%% file: PAPER_REVISED.tex
\documentclass[final,3p,times]{elsarticle}

\usepackage[breaklinks,colorlinks=true,linkcolor=blue,citecolor=red,backref=page]{hyperref}
\usepackage{amsfonts} 
\usepackage{amssymb,amsmath}
\usepackage[mathscr]{euscript} 
\usepackage{graphicx,epstopdf}
\usepackage{color}

\input{mathmacros.tex}

\journal{Wave Motion}

\begin{document}

\begin{frontmatter}

\title{Polarization effects for electromagnetic wave propagation in
  random media}

\author[label1]{Liliana Borcea}
\author[label2]{Josselin Garnier}

\address[label1]{Department of Mathematics,
University of Michigan,
Ann Arbor, MI 48109-1043
}
\address[label2]{Laboratoire de Probabilit\'es et Mod\`eles Al\'eatoires
\& Laboratoire Jacques-Louis Lions,
Universit{\'e} Paris Diderot,
75205 Paris Cedex 13,
France
}

\begin{abstract}
We study Maxwell's equations in random media with small fluctuations
of the electric permittivity. We consider a setup where the waves
propagate toward a preferred direction, called range. We decompose the
electromagnetic wave field in transverse electric and transverse
magnetic plane waves, called modes, with random amplitudes that model
cumulative scattering effects in the medium. Their evolution in range
is described by a coupled system of stochastic differential equations
driven by the random fluctuations of the electric permittivity. We
analyze the solution of this system with the Markov limit theorem and
obtain a detailed asymptotic characterization of the electromagnetic
wave field in the long range limit. In particular, we quantify the
loss of coherence of the waves due to scattering by calculating the
range scales (scattering mean free paths) on which the mean amplitudes
of the modes decay. We also quantify the energy exchange between the
modes, and consequently the loss of polarization induced by
scattering, by analyzing the Wigner transform (energy density) of the
electromagnetic wave field. This analysis involves the derivation of
transport equations with polarization.  We study in detail these
equations and connect the results with the existing literature in
radiative transport and paraxial wave propagation.
\end{abstract}


\end{frontmatter}

\section{Introduction}
Understanding the interaction of electromagnetic waves with complex
media through which they propagate is of great importance in
applications such as radar imaging and remote sensing
\cite{skolnik1970radar,elachi1990radar}, optical imaging
\cite{moscoso2008polarization}, laser beam propagation through the
atmosphere \cite{andrews2005laser}, and communications
\cite{goldsmith2005wireless}. As the microstructure (inhomogeneities)
of such media cannot be known in detail, we model it as a random
field, and thus study Maxwell's equations in random media. The goal is
to describe features of the solution, the electromagnetic wave field,
which do not depend on the particular realization of the random
medium, just on its statistics. Of particular interest in applications
are the statistical expectation of the solution, which describes the
coherent part of the waves, and the second moments which describe how
the waves decorrelate and depolarize, and how energy is transported in
the medium.

In most applications the inhomogeneities are weak scatterers, modeled
by small fluctuations of the wave speed. They have large cumulative
scattering effects at long distances of propagation. Among them are
the loss of coherence, manifested mathematically as an exponential
decay of the statistical expectation of the wave field, and thus
enhancement of the random fluctuations, and wave depolarization. The
quantification of these effects depends not only on the amplitude of
the fluctuations of the wave speed, but also on the relation between
the basic lengths scales: the wavelength, the scale of variations of
the medium (correlation length), the spatial support of the source,
and the distance of propagation.

Recent mathematical studies of electromagnetic waves in random media
are in \cite{kohler1991reflection} for layered media, in
\cite{alonso2015electromagnetic} for waveguides, and in
\cite{garnier2009paraxial} for beam propagation in open environments.
They decompose the electromagnetic field in plane wave components,
transverse to a preferred direction of propagation called range, and
then analyze the evolution of their amplitudes, which are frequency
and range dependent random fields.  The details and results differ
from one study to another, because the geometry and the scaling
regimes are different. For example, the decomposition in waveguides
leads to a countable set of waves, whereas in open environments, as we
consider here, there is a continuum of plane waves. The regime in
\cite{garnier2009paraxial} leads to statistical wave coupling by
scattering, however the waves form a narrow cone beam and 
they retain their initial linear polarization.

In this paper we build on the results in
\cite{alonso2015electromagnetic,garnier2009paraxial} to obtain a
detailed characterization of polarization effects in random open
environments. We consider a medium with random electric permittivity
that fluctuates on a scale (correlation length) $\ell$ that is larger
than the wavelength $\la$ by a factor $1/\gamma$, with $\gamma \in
(0,1)$, and the waves propagate over many correlation lengths.  
By
assuming that the support of the source is similar to the correlation
length and that the fluctuations of the medium are small and smooth
(i.e.  their standard deviation $\alpha$ is small and their power
spectral density decays fast enough), we identify an interesting
regime where the waves propagate along a preferred direction, the
range. It is between the paraxial regime studied in
\cite{garnier2009paraxial}, where the waves travel in the form of a
narrow cone beam, and the radiative transfer regime in
\cite{ryzhik96,bal00}, where the waves travel in all directions. In
our regime the waves propagate in a cone whose opening angle is
significantly larger than in the paraxial case, but smaller than $180$
degrees, so that the backscattered waves can be neglected.
The validity of this regime is controlled by the parameter $\gamma$ and the
distance $L$ of propagation.  We take $L \gg \ell = \la/\gamma$ so
that the scattering effects in the medium build up, but to ensure that
the waves remain forward propagating, we restrict the distance to $L
\sim \alpha^{-2} \ell$.

From the mathematical point of view, the advantage of having a
preferred direction of propagation is that we can reduce the analysis
of Maxwell's equations to the study of the range evolution of the random
amplitudes of the components of the wave field.  These amplitudes
satisfy a system of stochastic differential equations driven by the
fluctuations of the electric permittivity, and can be analyzed in
detail using the Markov limit theorem \cite{PK1974,PW1994}. Our main
results are:
\begin{enumerate}
\item The quantification of the scattering mean free
paths, the range scales on which the components of the electromagnetic
wave field lose coherence. 
\item The quantification of the statistical decorrelation of the waves
  over directions.
\item The derivation of the transport equations for the energy
  density, which allows us to quantify the depolarization of the waves
  and the diffusion of wave energy in direction.
 \item The connection of the results to the radiative transport theory with
  polarization given in \cite{chandra,ryzhik96,tsangkong}, and to the
  moment equations associated to the random paraxial wave
  equation given in
  \cite{Fannjiang,papa07,garnier2009paraxialderiv} (scalar case)
  and \cite{garnier2009paraxial} (electromagnetic case).
\end{enumerate}

The paper is organized as follows: We state up front, in Section
\ref{sect:stateres}, the mathematical problem and the main
results. The derivation of these results uses the formulation and
scaling described in Section \ref{sect:form}. The derivation of the
wave decomposition is in Section \ref{sect:WD}.  To give an intuitive
interpretation of the decomposition we consider first homogeneous
media. Then we give the decomposition in random media and derive the
system of stochastic differential equations satisfied by the random
wave amplitudes. The Markov limit of the solution of this system is
obtained in Section \ref{sect:DL}. We use it in Section \ref{sect:COH}
to quantify the loss of coherence of the waves. The analysis of the
second moments of the amplitudes and the derivation of the transport
equations is in Section \ref{sect:TEQ}. They are connected to the
radiative transport theory in \cite{chandra,ryzhik96,bal00} in
\ref{ap:RT}. We illustrate the results
with numerical simulations in statistically isotropic media in Section
\ref{sect:DEPOL}. We give a detailed analysis 
in the high-frequency limit $\gamma \to 0$ in Section \ref{sect:HF}
and connect these results to the white-noise paraxial wave equation
studied in \cite{garnier2009paraxial} in \ref{app:par}.  We
end with a summary in Section \ref{sect:sum}.

% ----------------------
\section{Statement of results}
\label{sect:stateres}
We begin in Section \ref{sect:MEq} with Maxwell's equations satisfied
by time-harmonic electric and magnetic fields in random media. The
scattering regime is defined in Section \ref{sect:scales}, by
identifying the important scales and describing their relations.  The
wave decomposition in transverse electric and magnetic modes is stated
in Section \ref{sect:WDec}. The characterization of the first and
second moments of the random amplitudes of these modes are our main
results, stated in Section \ref{sect:result}. They are proved in
Sections \ref{sect:COH} and \ref{sect:TEQ}.

\subsection{Maxwell's equations in random media}
\label{sect:MEq}
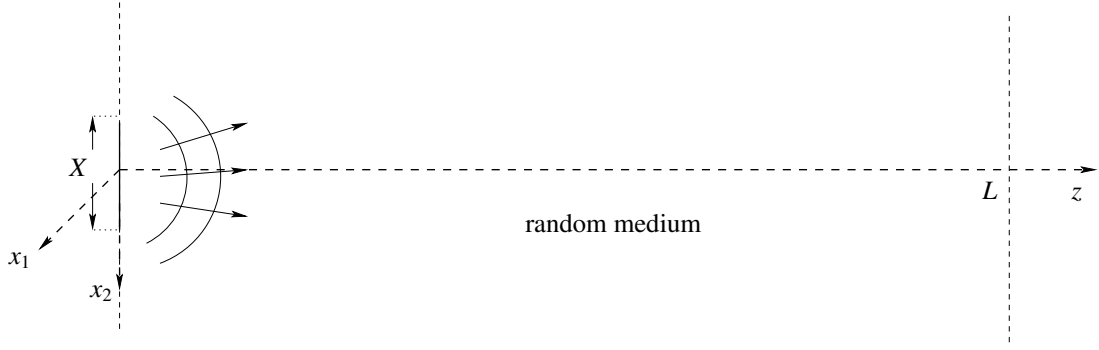
\begin{figure}[t]
  \centering
  \begin{picture}(0,100)%
\hspace{-2.7in}\includegraphics[width = 0.85\textwidth]{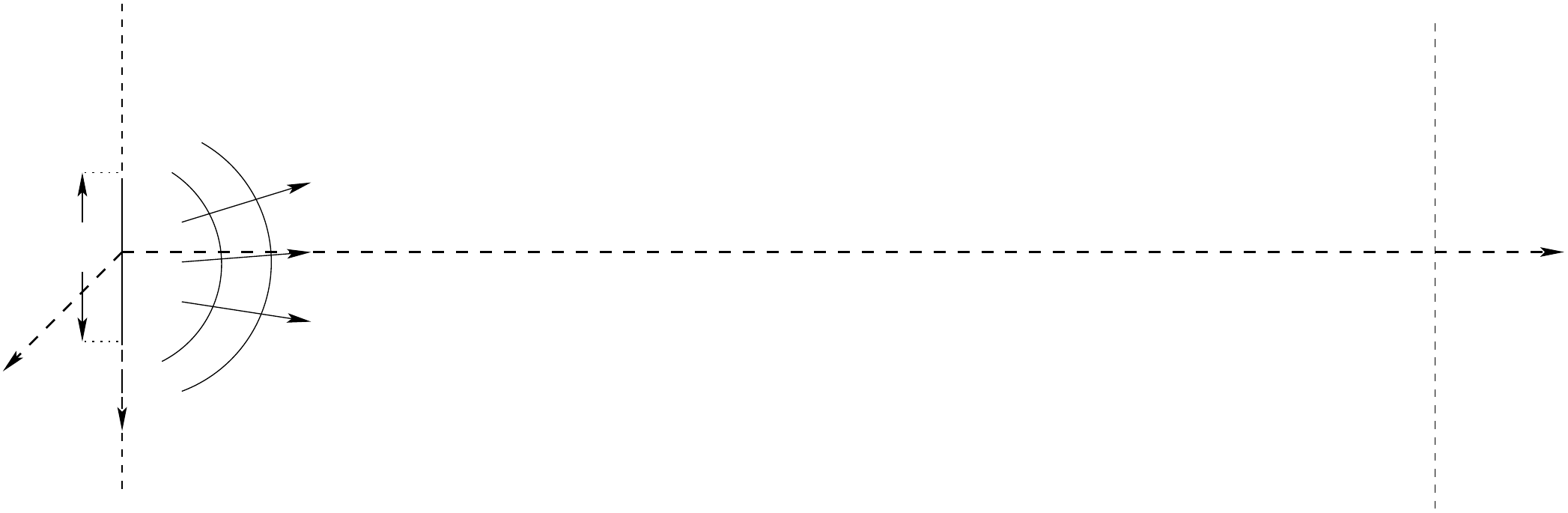}%
\end{picture}%
\setlength{\unitlength}{100sp}%
\begingroup\makeatletter\ifx\SetFigFont\undefined%
\gdef\SetFigFont#1#2#3#4#5{%
  \reset@font\fontsize{#1}{#2pt}%
  \fontfamily{#3}\fontseries{#4}\fontshape{#5}%
  \selectfont}%
\fi\endgroup%
\begin{picture}(0,100)(0,100)
\put(-135000,20000){\makebox(0,0)[lb]{\smash{{\SetFigFont{7}{8.4}{\familydefault}{\mddefault}{\updefault}{\color[rgb]{0,0,0}{\normalsize $x_1$}}%
}}}}
\put(-120000,41800){\makebox(0,0)[lb]{\smash{{\SetFigFont{7}{8.4}{\familydefault}{\mddefault}{\updefault}{\color[rgb]{0,0,0}{\normalsize $X$}}%
}}}}
\put(-115000,12000){\makebox(0,0)[lb]{\smash{{\SetFigFont{7}{8.4}{\familydefault}{\mddefault}{\updefault}{\color[rgb]{0,0,0}{\normalsize $x_2$}}%
}}}}
\put(126000,36000){\makebox(0,0)[lb]{\smash{{\SetFigFont{7}{8.4}{\familydefault}{\mddefault}{\updefault}{\color[rgb]{0,0,0}{\normalsize $z$}}%
}}}}
\put(104000,36000){\makebox(0,0)[lb]{\smash{{\SetFigFont{7}{8.4}{\familydefault}{\mddefault}{\updefault}{\color[rgb]{0,0,0}{\normalsize $L$}}%
}}}}
\put(-8000,28000){\makebox(0,0)[lb]{\smash{{\SetFigFont{7}{8.4}{\familydefault}{\mddefault}{\updefault}{\color[rgb]{0,0,0}{\normalsize random medium}}%
}}}}
\end{picture}%
\caption{Geometric setup. The source of diameter $X$ is localized at
  the origin of range $z$, and emits waves in the direction $z$, in
  the random medium extending from $z = 0$ to $z = L$.}
\label{fig:SETUP}
\end{figure}

The time-harmonic electric and magnetic fields $\vbE(\vx)$ and
$\vbH(\vx)$ satisfy Maxwell's equations
\begin{align}
\vna \times \vbE(\vx) & = i \om \mu_o 
\vbH(\vx), \label{eq:F1} \\ \vna \times \vbH(\vx) &=
\vec{\boldsymbol{\mathcal J}}(\vx) - i 
\om \eps(\vx)\vbE(\vx), \label{eq:F2} \\ \vna \cdot
\big[\eps(\vx) \vbE(\vx)\big] &= \rho(\vx), \label{eq:F3} \\ \vna
\cdot \big[\mu_o \vbH(\vx)\big] &= 0, \label{eq:F4}
\end{align}
for $\vx \in \mathbb{R}^3$ and frequency $\om \in \mathbb{R}$. The
setup is illustrated in Figure \ref{fig:SETUP}. The excitation is with
the current source density $\vec{\boldsymbol{\mathcal J}}$ localized
at $z= 0$, that emits waves in the direction $z$, called range. The
system of coordinates is $\vx = (\bx,z)$, with vector $\bx =
(x_1,x_2)$ in the cross-range plane.  The waves propagate in a linear
medium with random electric permittivity $\eps(\vx)$ and constant
magnetic permeability $\mu_o$. The medium is isotropic and
non-absorbing, meaning that $\eps(\vx)$ is scalar valued and
positive. Extensions to variable magnetic permeability and to
dissipative media with complex valued permittivity tensors are
possible, but for simplicity we do not consider them here.

We focus attention on equations (\ref{eq:F1}) and (\ref{eq:F2}),
because the other two equations are implied by them. Equation
(\ref{eq:F4}) follows by taking the divergence in (\ref{eq:F1}) and
the charge density $\rho(\vx)$ in equation (\ref{eq:F3}) is
related to the current source $\vec{\boldsymbol{\mathcal J}}(\vx)$
by the continuity of charge equation
\[
i \om \rho(\vx) = \vna \cdot \vec{\boldsymbol{\mathcal
    J}}(\vx),
\]
obtained by taking the divergence of equation (\ref{eq:F2}).

The random electric permittivity models small scale inhomogeneities in
the medium.  It fluctuates around the constant value $\eps_o$, as
described by
\begin{equation}
\eps_{\rm r}(\vx) = 1 + 1_{(0,L)}(z) \, {\alpha} \nu
\left(\frac{\vx}{\ell} \right),
\label{eq:F12}
\end{equation}
where $\eps_{\rm r}(\vx) = \eps(\vx)/\eps_o$ is the relative electric permittivity,
and $\nu$ is a dimensionless stationary random process of
dimensionless argument in $\mathbb{R}^3$. The process has zero mean
\begin{equation}
\EE \big[ \nu(\vr) \big] = 0, 
\label{eq:F13}
\end{equation}
and autocorrelation
\begin{equation}
 \EE\big[\nu(\vr)\nu(\vr')\big] = \cR(\vr-\vr') ,\quad 
\quad \forall \, \vr, \vr' \in \mathbb{R}^3.
\label{eq:F14}
\end{equation}
We assume that $\nu$ is bounded and differentiable, with bounded
derivative almost surely, and that $\cR(\vr)$ is integrable, with
Fourier transform (power spectral density) 
\begin{equation}
\widetilde \cR( \vec{\bq}) = \int_{\mathbb{R}^3} d \vr\, \cR(\vr)
e^{- i \vec{\bq} \cdot \vr}.
\label{eq:PSDensity}
\end{equation}
This is a non-negative function by Bochner's theorem, and we suppose
that it is is either compactly supported in a ball in
$\mathbb{R}^3$ of radius $O(1)$, or it is negligible outside this ball. The
autocorrelation is normalized so that
\begin{equation}
\int_{\mathbb{R}^3} d \vr \, \cR(\vr) =
O(1), \qquad \cR({\bf 0}) = O(1).
\label{eq:F15}
\end{equation}
The length scale $\ell$ in (\ref{eq:F12}) is the correlation length and the
positive and small dimensionless parameter ${\alpha}$ quantifies the
typical amplitude (standard deviation) of the fluctuations.

The indicator function $1_{(0,L)}(z)$ in (\ref{eq:F12}) limits the
support of the fluctuations to the range interval $z \in (0,L)$. This
allows us to state easily the boundary conditions satisfied by $\vbE$
and $\vbH$. The truncation at $z = L$ may be understood physically in
the time domain, where $L$ is determined by the duration of the
observation time. For a finite time the waves emitted by the source
cannot be affected by the medium beyond some range, called here $L$,
which is why we can truncate the random fluctuations there without
affecting the solution. The truncation at $z = 0$ is consistent with
the fact that in our scaling the backscattered waves are
negligible. Thus, the waves that propagate from the source in the
negative range direction have no effect on the waves at $z > 0$, which
is why we can truncate the fluctuations at $z = 0$.

\subsection{The scattering regime}
\label{sect:scales}
There are four important length scales that define the scattering
regime: the wavelength $\la$, the correlation length $\ell$, the
distance of propagation $\bar{L}$, and the support $X$ of the source in
cross-range. The wavelength is defined at the reference wave speed
$c_o = 1/\sqrt{\eps_o \mu_o}$ by $ \la = {2 \pi c_o}/{\om}$, and
$\bar{L}$ is of the same order as $L$. The support $X$ of the source
in cross-range determines the initial opening angle of the cone (beam)
of emitted waves.  To define it we let the current source be of the
form
\begin{equation}
\vec{\boldsymbol{\mathcal J}}(\vx) = \zeta_o^{-1}
\vbJ\left(\frac{\bx}{X} \right) \delta(z),
\label{eq:F16}
\end{equation}
where $\vbJ(\br) = \big(\bJ(\br),\mJ_z(\br)\big)$ is a vector-valued
function of the dimensionless vector $\br \in \mathbb{R}^2$. The
magnitude of $\vbJ$ is assumed negligible for $|\br| > 1$, so $X$
is the diameter of the spatial support of the source. The scaling by the constant
impedance $\zeta_o = \sqrt{\mu_o/\eps_o}$ is chosen for convenience.

Our scattering regime is defined by the scale ordering
\begin{equation}
\la   < \ell \sim X \ll \bar{L},
\label{eq:hypscal1}
\end{equation}
and the small standard deviation $\alpha$ of the random fluctuations
of the electric permittivity, satisfying
\begin{equation}
\alpha^2  \sim \ell / \bar{L}.
\label{eq:hypscal1b}
\end{equation}
Here the symbol $\sim$ denotes  the same order.
By taking a large distance of propagation $\bar{L}$, as in
(\ref{eq:hypscal1}), we ensure that the cumulative scattering effects
in the random medium have a significant effect on the electromagnetic
wave field. This effect is summarized in Section \ref{sect:result} by
the exponential decay in $z$ of the statistical expectation of the
wave field, which is a manifestation of its randomization, the
exchange of energy between the wave modes, and the diffusion of energy
over directions of propagation. Because of this diffusion we restrict
the propagation distance $\bar{L}$ as in (\ref{eq:hypscal1b}), so that
the forward going waves emitted by the source remain forward going for
all $z \in (0, \bar{L})$.

\subsection{Transverse electric and magnetic modes}
\label{sect:WDec}
The interaction of the waves with the random medium depends on the
direction of propagation, which is why we decompose the electric and
magnetic fields in plane waves with wave vector $\vka$.  The
decomposition stated here is derived in Section \ref{sect:WD3}, and in
our scattering regime it consists of two types of forward going waves
which we call modes: transverse electric and magnetic.

The decomposition of the electric field is
\begin{align}
  \vbE(\vx) =\int_{|\bka|<1} \frac{d(k \bka)}{(2 \pi)^2} \,&
  \beta^{-\frac{1}{2}}(\bka)\big[a(\bka,z) \, \vu(\bka) +
    a^{\perp}(\bka,z) \, \vu^\perp(\bka) \big] e^{i k
    \vec{\bka} \cdot \vx}, \label{eq:EField_SR}
\end{align}
and for the magnetic field we have 
\begin{align}
  \vbH(\vx) =\zeta_o^{-1} \int_{|\bka|<1} \frac{d(k \bka)}{(2 \pi)^2}
  &\beta^{-\frac{1}{2}}(\bka) \big[a(\bka,z) \, \vu^\perp(\bka) -
    a^{\perp}(\bka,z) \, \vu(\bka) \big] e^{ik \vec{\bka} \cdot
    \vx}, \label{eq:MField_SR}
\end{align}
where $k = 2 \pi /\la$ is the wavenumber and we use the notation $d(k
\bka) = k^2 d \bka$.  The three-dimensional wave vector is
\begin{equation} 
\vec{\bka} = \big( \bka, \beta(\bka)\big), \qquad
\bka = (\kappa_1,\kappa_2) \in \mathbb{R}^2.
\label{eq:WD3}
\end{equation}
It satisfies $|\vka| = 1$, and the third component is
\begin{equation}
\beta(\bka) = \sqrt{1-|\bka|^2}.
\label{eq:BetaDef}
\end{equation}
This is real valued for $|\bka| < 1$, meaning that the waves are
propagating and not evanescent. In general there are forward and
backward propagating waves, as explained in Section \ref{sect:WD3},
but in our scattering regime the waves move forward, which is why we
have $\beta(\bka)$ with positive sign in (\ref{eq:WD3}).

The vectors 
\begin{equation}
\vu(\bka) = \Big(\beta(\bka) \frac{\bka}{|\bka|}, 
   -|\bka| \Big),  \qquad
\vu^\perp(\bka)= \Big( \frac{\bka^\perp}{|\bka|}, 0 \Big),
\label{eq:ORT_SR}
\end{equation} 
where $\bka^\perp = (-\kappa_2,\kappa_1)$, distinguish the two modes
associated with wave vector $\vka$.  Note that the triplet $\{ 
\vu(\bka),\vu^\perp(\bka),\vka \}$ is an orthonormal basis of
$\mathbb{R}^3$. Thus, equations (\ref{eq:EField_SR}) and
(\ref{eq:MField_SR}) are decompositions in transverse waves, which are
orthogonal to the direction of propagation along $\vka$.  The
amplitudes $a(\bka,z)$ are for the transverse magnetic modes, because
they do not contribute to the longitudinal component of the magnetic
field. They multiply the vector $\vu^\perp(\bka)$ in equation
(\ref{eq:MField}). Similarly, $a^{\perp}(\bka,z)$ are the amplitudes
of transverse electric modes, because they do not contribute to the
longitudinal electric field.

For any wave vector $\vka$, the electric and magnetic plane waves
\begin{align*}
\vec{\boldsymbol{\mathscr{E}}}(\bka,z) &= \beta^{-\frac{1}{2}}(\bka)
\big[a(\bka,z) \vu(\bka) + a^\perp(\bka,z) \vu^\perp(\bka)\big], \\
\vec{\boldsymbol{\mathscr{H}}}(\bka,z) &= \zeta_o^{-1}\beta^{-\frac{1}{2}}(\bka)
\big[a(\bka,z) \vu^\perp(\bka) - a^\perp(\bka,z) \vu(\bka)\big],
\end{align*}
are orthogonal to $\vka$ and to each other. The amplitudes $a(\bka,z)$
and $a^\perp(\bka,z)$ are random, and encode the scattering effects in
the medium. In the system of coordinates with axes along $\vu(\bka)$, 
$\vu^\perp(\bka)$, and $\vka$, which is right-handed because 
\[
\vu(\bka) \times \vu^\perp(\bka) = \vka, \qquad \vka \times \vu(\bka) = \vu^\perp(\bka), \qquad 
\vu^\perp(\bka) \times \vka = \vu(\bka),
\]
we can quantify the evolution of the state of polarization of the waves using the  coherence matrix 
\begin{equation}
   \boldsymbol{\mathscr{P}}(\bka,z) = \EE \left[ \begin{pmatrix}
      a(\bka,z) \\ a^{\perp}(\bka,z) \end{pmatrix}
 \begin{pmatrix}
   a(\bka,z) \\ a^{\perp}(\bka,z) \end{pmatrix}^\dagger \right] = 
   \frac{1}{2} \begin{pmatrix} \mathfrak{S}_1(\bka,z) + \mathfrak{S}_2(\bka,z) & \mathfrak{S}_3(\bka,z) + i \mathfrak{S}_4(\bka,z) \\
   \mathfrak{S}_3(\bka,z) -i  \mathfrak{S}_4(\bka,z)  & \mathfrak{S}_1(\bka,z) - \mathfrak{S}_2(\bka,z) \end{pmatrix},
   \label{eq:Stokes}
\end{equation}
where $\dagger$ denotes complex conjugate and transpose, and 
\begin{align*}
\mathfrak{S}_1(\bka,z) &= \EE \big[ |a(\bka,z)|^2 +
|a^\perp(\bka,z)|^2 \big] , \qquad  \mathfrak{S}_2(\bka,z) =
\EE \big[ |a(\bka,z)|^2 -
|a^\perp(\bka,z)|^2 \big],\\ \mathfrak{S}_3(\bka,z) &= 2
{\rm Re} \big( \EE \big[a(\bka,z) \overline{a^\perp(\bka,z)}
  \big] \big), \qquad \mathfrak{S}_4(\bka,z) = 2 {\rm Im} \big(
\EE \big[a(\bka,z) \overline{a^\perp(\bka,z)} \big] \big),
\end{align*}
are the components of the Stokes vector $\boldsymbol{\mathfrak{S}}(\bka,z) = \big( \mathfrak{S}_1(\bka,z),
\mathfrak{S}_2(\bka,z) , \mathfrak{S}_3(\bka,z), \mathfrak{S}_4(\bka,z)
\big)$.  
We use henceforth the bar to denote complex conjugate. 
The degree of polarization is defined by 
\begin{equation}
{\rm Pol}(\bka,z) = \frac{\big( \mathfrak{S}_2(\bka,z)^2+\mathfrak{S}_3(\bka,z)^2+\mathfrak{S}_4(\bka,z)^2\big)^{1/2}}{\mathfrak{S}_1(\bka,z) }. 
\end{equation}
It is a number in the interval $[0,1]$, due to the Cauchy-Schwarz inequality 
\[\big|\EE \big[a(\bka,z) \overline{a^\perp(\bka,z)} \big]\big|^2 \leq 
\EE \big[|a(\bka,z)|^2\big]\EE\big[|  a^\perp(\bka,z)|^2 \big].\] 
The wave is fully polarized when  ${\rm Pol}(\bka,z) = 1$  and unpolarized when 
${\rm Pol}(\bka,z) = 0$. The mode amplitudes of fully polarized waves have a deterministic linear dependence so 
that 
\[
\EE[ |a(\kappa,z)|^2] \EE[ |a^\perp(\kappa,z)|^2] = \left| \EE [ a(\kappa,z) \overline{a^\perp(\bka,z)} ] \right|^2.
\]
The mode amplitudes of  unpolarized waves are uncorrelated and have the same mean power.

\subsection{Main results}
\label{sect:result}%
If the medium were homogeneous, the modes would be independent, with
constant amplitudes $a_o(\bka)$ and $a_o^\perp(\bka)$ defined by the
current source (\ref{eq:F16}),
\begin{align}
a_o(\bka) &= \frac{\beta^{-1/2}(\bka)\, |\bka|}{2} \int_{\mathbb{R}^2}
d \bx \, J_z\Big(\frac{\bx}{X} \Big) e^{- i k \bka \cdot \bx} -
\frac{\beta^{1/2}(\bka)}{2} \frac{\bka}{|\bka|} \cdot
\int_{\mathbb{R}^2} d \bx\, \bJ\Big(\frac{\bx}{X} \Big) e^{- i k
  \bka \cdot \bx}, \label{eq:MR4} \\ 
  a_o^\perp(\bka) &=
-\frac{\beta^{-1/2}(\bka)}{2} \frac{\bka^\perp}{|\bka|} \cdot
\int_{\mathbb{R}^2} d \bx\, \bJ \Big(\frac{\bx}{X} \Big) e^{- i k
  \bka \cdot \bx}. \label{eq:MR5}
\end{align}
In the random medium, in the scattering regime defined by
(\ref{eq:hypscal1}-\ref{eq:hypscal1b}), the mode amplitudes
$\{a(\bka,z),a^\perp(\bka,z)\}$ for $|\bka| < 1$ and $z \lesssim \bar{L}$
form a Markov process whose statistical moments can be characterized
explicitly. We are interested in the first two moments which describe
the loss of coherence of the waves and the transport of energy by the
modes.

\subsubsection{The mean mode amplitudes}
\label{sect:result.coh}
Let us denote the vector of the expectations of the mode amplitudes at
wave vector $\vka$ by
\begin{equation}
  \boldsymbol{\mathscr{A}}(\bka,z) = \begin{pmatrix} \EE[a(\bka,z)]
    \\ \EE[a^\perp(\bka,z)]
  \end{pmatrix}.
\end{equation}
Its initial value in the source plane $z=0$ is 
\begin{equation}
\label{eq:defAo}
\boldsymbol{\mathscr{A}}(\bka,0) = \boldsymbol{\mathscr{A}}_o(\bka)
= \begin{pmatrix} a_o(\bka) \\ a_o^\perp(\bka)
  \end{pmatrix},
\end{equation}
and for $z > 0$ we have 
\begin{equation}
 \boldsymbol{\mathscr{A}}(\bka,z)= \exp\Big[ {\bf Q}(\bka)z \Big]
 \boldsymbol{\mathscr{A}}_o(\bka).
\label{eq:MR3}
\end{equation}
Thus, the scattering effects in the random medium do not average out,
and the range evolution of the mean mode amplitudes is determined by
the matrix ${\bf Q}(\bka)$ in the exponential. This complex
symmetric matrix is given by 
\begin{align}
 {\bf Q}(\bka) = &-\frac{k^2 {\alpha}^2}{4} \int_{|\bka'|< 1} \frac{d
   (k \bka')}{(2 \pi)^2}\, \bGamma(\bka,\bka')
 \bGamma(\bka',\bka)\int_{\mathbb{R}^2} d \bx \int_0^\infty d z \,
 \cR\left(\frac{\vx}{\ell}\right) e^{ - i k
   (\vka-\vka') \cdot \vx } -\frac{i k {\alpha}^2}{2} \cR({\bf 0})
 \frac{|\bka|^2}{\beta(\bka)} \begin{pmatrix} 1 & 0 \\ 0
   &0 \end{pmatrix},
\label{eq:MR6}
\end{align}
for $\vka,\vka'$ defined as in (\ref{eq:WD3}),   $\vx = (\bx,z)$, and  
\begin{equation}
\bGamma(\bka,\bka') = \begin{pmatrix} \frac{|\bka|
    |\bka'|}{\sqrt{\beta(\bka) \beta(\bka')}} + \frac{\bka}{|\bka|}
  \cdot \frac{\bka'}{|\bka'|}  {\small \sqrt{\beta(\bka)\beta(\bka')}} &
  \frac{\bka}{|\bka|} \cdot \frac{\bka'^{\perp}}{|\bka'|}
  \sqrt{\frac{\beta(\bka)}{\beta(\bka')}} \\ \frac{\bka^\perp}{|\bka|}
  \cdot \frac{\bka'}{|\bka'|} \sqrt{\frac{\beta(\bka')}{\beta(\bka)}}
  & \frac{\bka}{|\bka|} \cdot \frac{\bka'}{|\bka'|}
  \frac{1}{\sqrt{\beta(\bka')\beta(\bka)}}\end{pmatrix}.
\label{eq:MR1}
\end{equation}
Note that $\bGamma(\bka',\bka) = \bGamma(\bka,\bka')^T$ and since the
autocorrelation $\cR$ is an even function, and the power spectral
density $\widetilde \cR$ is non-negative, the real part of ${\bf
  Q}(\bka)$ is a symmetric, negative definite matrix
\begin{equation}
{\rm Re} \big[{\bf Q}(\bka)\big] = - \frac{k^2
  \ell^3 \alpha^2}{8} \int_{|\bka'|< 1} \frac{d (k \bka')}{(2
  \pi)^2}\, \bGamma(\bka,\bka') \bGamma^T(\bka,\bka')
\, \widetilde{\cR}\big(k \ell (\vka-\vka') \big).
\end{equation}

In the case of transverse isotropic statistics of the media, where
$\cR(\vr)$ with $\vr = (\br,r_z)$ depends only on $|\br|$ and $r_z$,
we obtain by integrating in \eqref{eq:MR6} over $\bka'$ in polar
coordinates that $ {\bf Q}(\bka)$ is diagonal and moreover, it only
depends on $|\bka|$.  Therefore, the mean wave amplitudes decouple and
they decay exponentially in $z$ on the range scales (scattering mean
free paths)
\begin{align}
\cS(\bka) = -1/{\rm Re} \big[ {\rm Q}_{11}(\bka) \big], \qquad
\cS^\perp(\bka) = -1/{\rm Re} \big[  {\rm Q}_{22}(\bka) \big] ,
\label{eq:MR7}
\end{align}
which depend only on $|\bka|$.  We refer to Section \ref{sect:DEPOL}
for illustrations with numerical simulations in statistically 
isotropic media, and to Section \ref{sect:HF} for the analysis of the high-frequency limit
and connection to the white-noise paraxial regime described in \ref{app:par}.  We
also show in  \ref{ap:RT} that, if $\cR$ is isotropic i.e.,
$\cR(\vr)$ only depends on $|\vr|$, then ${\rm Re}({\bf
  Q}(\bka))$ is proportional to the identity matrix, so the scattering
mean free paths (\ref{eq:MR7}) are equal.

In general media the mean mode amplitudes are coupled, but they still
decay exponentially in $z$. Indeed, matrix 
\begin{equation}
{\bf S}(\bka) = - {\bf Q}(\bka) - {\bf Q}(\bka)^\dagger = \frac{k^2
  \ell^3 \alpha^2}{4} \int_{|\bka'|< 1} \frac{d (k \bka')}{(2
  \pi)^2}\, \bGamma(\bka,\bka') \bGamma^T(\bka,\bka')
\, \widetilde{\cR}\big(k \ell (\vka-\vka') \big)
\label{eq:defSM}
\end{equation}
is real symmetric and positive definite. Therefore, the mean amplitudes
satisfy
\begin{equation}
  \partial_z \|\boldsymbol{\mathscr{A}}(\bka,z)\|^2 =
  \boldsymbol{\mathscr{A}}(\bka,z)^\dagger \Big[ {\bf Q}(\bka)^\dagger
    + {\bf Q}(\bka) \Big] \boldsymbol{\mathscr{A}}(\bka,z) = -
  \boldsymbol{\mathscr{A}}(\bka,z)^\dagger {\bf S}(\bka)
  \boldsymbol{\mathscr{A}}(\bka,z), 
\end{equation}
and if we let $\Lambda_1(\bka) \ge \Lambda_2(\bka) > 0$ be the eigenvalues of
${\bf S}(\bka)$ and use Gronwall's lemma, we get
\begin{equation}
  e^{-z{\Lambda_1(\bka)}/{2}} \|\boldsymbol{\mathscr{A}}_o(\bka)\| \le
  \|\boldsymbol{\mathscr{A}}(\bka,z)\| \le e^{-z\Lambda_2(\bka)/2}
  \|\boldsymbol{\mathscr{A}}_o(\bka)\|.
\end{equation}
The scales $2/\Lambda_1(\bka)$ and $2/\Lambda_2(\bka)$, on which the
mean wave amplitudes decay exponentially in range, are the scattering
mean free paths.

The exponential decay of the mean mode amplitudes, and therefore of
the mean (coherent) fields $\EE \big[ \vbE(\vx)\big]$ and $\EE
\big[\vbH(\vx)\big]$, is a manifestation of the randomization of the
waves due to scattering. Energy is conserved, as shown in Section
\ref{sect:DiffA}, 
\begin{equation}
\int_{|\bka| < 1} \frac{d (k \bka)}{(2 \pi)^2} \Big[ |a(\bka,z)|^2 +
  |a^\perp(\bka,z)|^2 \Big] = ~ \mbox{constant},
\label{eq:MR2}
\end{equation}
so the second moments cannot all decay. Therefore, the mode amplitudes
have large random fluctuations. 

\vspace{0.05in}
\textbf{Remark 1.}
The exponential decay of the coherent field is a general phenomenon in
random media, and can be characterized under many scattering regimes,
including when the forward scattering approximation is not valid. We
refer to \cite{vanrossum} for an overview.

\subsubsection{The transport equations}
\label{sect:result.mp}
The energy density (Wigner transform) of the modes at wave vector
$\vka$ is defined by 
\begin{equation}
\bcW(\bka,\vx) = \int \frac{d(k {\itbf q})}{(2 \pi)^2} \, \exp \big[ i
  k {\itbf q} \cdot \big(\bx + \nabla \beta({\itbf q}) z\big)\big]\,
\EE \left[ \begin{pmatrix} a\big(\bka + \frac{{\itbf q}}{2},z\big)
    \\ a^{\perp}\big(\bka + \frac{ {\itbf q}}{2},z\big) \end{pmatrix}
 \begin{pmatrix}
   a\big(\bka-\frac{{\itbf q}}{2},z\big) \\ a^{\perp}\big(\bka-\frac{
     {\itbf q}}{2},z\big) \end{pmatrix}^\dagger \right], \qquad \vx =
(\bx,z),
\label{eq:MR8}
\end{equation}
where the integral is over all ${\itbf q}$ such that $|\bka \pm {\itbf
  q}/2| < 1$. It satisfies the transport equation
\begin{align}
\partial_z \bcW(\bka,\vx) - \nabla \beta(\bka) \cdot \nabla_\bx
\bcW(\bka, \vx) = \mathbf{Q}(\bka) \bcW(\bka,\vx) + \bcW(\bka,\vx)
\mathbf{Q}(\bka)^\dagger \nonumber \\  + \frac{k^2 \ell^3 {\alpha}^2}{4
}\int_{|\bka'| < 1} \frac{d (k \bka')}{(2 \pi)^2}\,
\bGamma(\bka,\bka') \bcW(\bka',\vx) \bGamma(\bka',\bka) \,
\widetilde \cR\left(k \ell (\vka-\vka')\right),
  \label{eq:MR9}
\end{align}
for $z>0$. The derivation of this equation, given in Section
\ref{sect:TEQ}, is one of the main results of this paper. We
show some numerical illustrations in Section \ref{sect:DEPOL} 
in statistically isotropic media and
relate it to the radiative transport theory in
\cite{chandra,ryzhik96,bal00} in \ref{ap:RT}. We consider
the high frequency limit $\la/\ell \to 0$ in Section \ref{sect:HF}, where we can connect the
results to those in the white-noise paraxial regime
\cite{garnier2009paraxial} as described in \ref{app:par}.

The integral over $\bx$ of the Wigner transform is the Hermitian,
positive definite coherence matrix
\begin{equation}
   \boldsymbol{\mathscr{P}}(\bka,z) = \EE \left[ \begin{pmatrix}
      a(\bka,z) \\ a^{\perp}(\bka,z) \end{pmatrix}
 \begin{pmatrix}
   a(\bka,z) \\ a^{\perp}(\bka,z) \end{pmatrix}^\dagger \right] =
   \int_{\mathbb{R}^2} d \bx \, \bcW(\bka,\vx).
  \label{eq:MR10}
\end{equation}
Its diagonal entries are the mean mode powers, and the relation with the Stokes vector is as in 
\eqref{eq:Stokes}. The
coherence matrix evolves from its initial value
$\boldsymbol{\mathscr{P}}(\bka,0) = \boldsymbol{\mathscr{A}}(\bka,0)
\boldsymbol{\mathscr{A}}(\bka,0)^\dagger$ at the source according to
equation
\begin{align}
 \partial_z \boldsymbol{\mathscr{P}}(\bka,z) =& \mathbf{Q}(\bka)
 \boldsymbol{\mathscr{P}}(\bka,z) + \boldsymbol{\mathscr{P}}(\bka,z)
 \mathbf{Q}(\bka)^\dagger + \frac{k^2 \ell^3 {\alpha}^2}{4
 }\int_{|\bka'| < 1} \frac{d (k \bka')}{(2 \pi)^2}\,
 \bGamma(\bka,\bka') \boldsymbol{\mathscr{P}}(\bka',z)
 \bGamma(\bka',\bka) \, \widetilde \cR\left(k \ell
 (\vka-\vka')\right).
\label{eq:MR13}
\end{align}
It follows by direct calculation, using the commutation properties of
the trace and \eqref{eq:defSM}, that
\begin{equation}
  \partial_z \int_{|\bka|<1} \frac{d (k\bka)}{(2 \pi)^2}\, {\rm Tr}
  \Big[ \boldsymbol{\mathscr{P}}(\bka,z)\Big] = \int_{|\bka|<1}\frac{d
    (k\bka)}{(2 \pi)^2}\, {\rm Tr} \Big[ \Big(\mathbf{Q}(\bka) +
    \mathbf{Q}(\bka)^\dagger + {\bf S}(\bka) \Big)
    \boldsymbol{\mathscr{P}}(\bka,z)\Big] = 0, 
\label{eq:MR12}
\end{equation}
which is consistent with the conservation identity
(\ref{eq:MR2}). Moreover,  (\ref{eq:MR13}) can be written in integral form as
\begin{align*}
   \boldsymbol{\mathscr{P}}(\bka,z) = e^{\mathbf{Q}(\bka) z}
   \boldsymbol{\mathscr{P}}(\bka,0) e^{\mathbf{Q}(\bka)^\dagger z} +
   \frac{k^2 \ell^3{\alpha}^2}{4} \int_0^z dz' \int_{|\bka'| < 1}
   \frac{d (k \bka')}{(2 \pi)^2}\, \widetilde
   \cR\left( k \ell (\vka-\vka')\right) \times 
 \nonumber \\ \times
   \Big[e^{\mathbf{Q}(\bka) \, (z-z')}\bGamma(\bka,\bka')\Big]
   \boldsymbol{\mathscr{P}}(\bka',z')\Big[e^{\mathbf{Q}(\bka) \,
       (z-z')} \bGamma(\bka,\bka)\Big]^\dagger,
\end{align*}
where we used the symmetry relation $\bGamma(\bka',\bka) =
\bGamma(\bka,\bka')^T$. Since $ \boldsymbol{\mathscr{P}}(\bka,0)$ is
Hermitian positive definite by definition, and the power spectral
density $\widetilde \cR$ is non-negative, we conclude that 
equation (\ref{eq:MR13}) has a Hermitian positive definite solution 
for all $z$, which is consistent with (\ref{eq:MR10}).

% --------------------------
\section{Mathematical formulation of the problem}
\label{sect:form}

In this section we give the mathematical formulation for the analysis
of the electromagnetic wave field in the random medium with electric
permittivity (\ref{eq:F12}) . We begin in Section \ref{sect:form.1}
with the derivation of the $4\times4$ system of partial differential
equations satisfied by the transverse components of the electric and
magnetic field. Our goal is to analyze the solution of this system in
the scattering regime (\ref{eq:hypscal1}-\ref{eq:hypscal1b}),
made precise  in Section \ref{sect:form.3}. This scaling
defines the asymptotic regime in which we can characterize the
statistics of the mode amplitudes using a Markov limit.

\subsection{The equations for the transverse waves}
\label{sect:form.1}
Let us eliminate from Maxwell's equations (\ref{eq:F1}) and
(\ref{eq:F2}) the longitudinal components $\mE_z$ and $\mH_z$ of $\vbE
= (\bE,\mE_z)$ and $\vbH = (\bH,\mH_z)$, and derive a closed system of
partial differential equations for the transverse electric and
magnetic fields $\bE$ and $\bH$. We obtain that
\begin{align}
\mH_z(\vx) &= -\frac{i}{\om \mu_o} \nabla^\perp \cdot
\bE(\vx), \label{eq:F5} \\ \mE_z(\vx) &= \frac{i}{\om \eps(\vx)} \big[
  \nabla^\perp \cdot \bH(\vx) -\vec{\itbf e}_z \cdot \boldsymbol{\mathcal
    J}(\vx)\big], \label{eq:F6}
\end{align}
where $\vec{\itbf e}_z$ is the unit vector along the range axis, $\nabla
= (\partial_{x_1},\partial_{x_2})$ is the gradient in the transverse
coordinates and $\nabla^\perp = (-\partial_{x_2},\partial_{x_1})$ is
its rotation by $90$ degrees, counterclockwise. The transverse components satisfy
\begin{align}
\partial_z \bE(\vx) &= -i \om \mu_o \bH^\perp(\vx) - \frac{i}{\om}
\nabla \left[\frac{\nabla \cdot \bH^\perp(\vx)}{\eps(\vx)} \right] -
\frac{i \zeta_o^{-1}}{\om} \nabla_{\bx} \left[\frac{
    {J}_z(\bx/X)}{\eps(\vx)}\right] \delta(z), \label{eq:F7}
\\ \partial_z \bH^\perp(\vx) &= - i \om \eps(\vx) \bE(\vx) -
\frac{i}{\om \mu_o} \nabla^\perp \big[ \nabla^\perp \cdot
  \bE(\vx)\big] + \zeta_o^{-1}\bJ(\bx/X) \delta(z) , \label{eq:F8}
\end{align}
with $\bH^\perp = (-\mH_2,\mH_1)$, the rotation of $\bH =
(\mH_1,\mH_2)$ by $90$ degrees, counterclockwise.  We study the
$4\times 4$ system (\ref{eq:F7}-\ref{eq:F8}), written in more
convenient form in terms of the scaled, rotated magnetic field
\begin{equation}
\bU(\vx) = - \zeta_o \bH^\perp(\vx).
\label{eq:F9}
\end{equation}
The scaling by the impedance $\zeta_o$ ensures that $\bE$ and $\bU$
have the same physical units. Equations (\ref{eq:F7}-\ref{eq:F8})
become
\begin{align} 
\partial_z \bE(\vx) &= i k \bU(\vx) + \frac{i}{k} \nabla
\left[\frac{\nabla \cdot \bU(\vx)}{\eps_{\rm r}(\vx)} \right] -
\frac{i}{k} \nabla_\bx \left[\frac{ {J}_z(\bx/X)}{\eps_{\rm
      r}(\vx)}\right]\delta(z), \label{eq:F10} \\ \partial_z \bU(\vx)
&= i k \eps_{\rm r}(\vx) \bE(\vx) + \frac{i}{k} \nabla^\perp \left[
  \nabla^\perp \cdot \bE(\vx) \right] -\bJ(\bx/X)
\delta(z). \label{eq:F11}
\end{align}

\subsection{Asymptotic regime}
\label{sect:form.3}%
We model the separation of scales in
(\ref{eq:hypscal1}-\ref{eq:hypscal1b}) with two dimensionless
parameters $\ep$ and $\gamma$, ordered as
\begin{equation}
0 < \ep \ll \gamma < 1,
\label{eq:hypscal2}
\end{equation}
and defined by 
\begin{equation}
\ep  = \frac{\la}{\bar{L}}, \qquad \gamma = \frac{\la}{\ell}.
\label{eq:F17}
\end{equation}
Our primary asymptotic analysis is for $\ep \to 0$. The parameter
$\gamma$ remains finite, except in Section \ref{sect:HF} where we
study the high-frequency limit $\gamma \to 0$ in which the transport
equations simplify and become equivalent to those in the white-noise paraxial regime.

We let $\bar{L}$ be the reference length scale, and introduce the scaled
length variables
\begin{equation}
\bx' = \bx/(\ep \bar{L}), \quad z' = z/\bar{L}, \quad L' = L/\bar{L}, \quad
X' = X/\bar{L} = \ep /\gamma_{_\mJ} .
\label{eq:F18}
\end{equation} 
The scaled wavenumber is $k' = k \bar{L} \ep = 2\pi$ and the
normalized standard deviation of the fluctuations is 
\begin{equation}
\label{def:alpha}
{\alpha}' =
{\alpha} / \ep^{1/2} = O(1).
\end{equation}
As we will see below, the dimensionless parameter
\begin{equation}
\gamma_{_\mJ} = \gamma \ell/X = O(\gamma),
\label{eq:defGJ}
\end{equation}
controls the opening angle of the cone (beam) emitted by the
source.

Substituting (\ref{eq:F18}) into (\ref{eq:F10}-\ref{eq:F11}) and
dropping all the primes, since all the variables are scaled
henceforth, we obtain
\begin{align}
\partial_z \bE^\ep(\bx,z) = \frac{ik}{\ep} \left[ {\bf I} +
  \frac{1}{k^2} \nabla\big(\nabla \cdot \big)\right]\bU^\ep(\bx,z) -
\frac{i {\alpha}}{k \ep^{1/2}} \nabla \left[ \nu \left(\gamma
  \bx,\gamma \frac{z}{\ep}\right) \nabla \cdot \bU^\ep(\bx,z)\right]
\nonumber \\ + \frac{i {\alpha}^2}{k} \nabla \left[\nu^2\left(\gamma
  \bx,\gamma\frac{z}{\ep}\right) \nabla \cdot \bU^\ep(\bx,z)\right] -
\frac{i}{k} \nabla_\bx \mJ_z\big(\gamma_{_\mJ} {\bx} \big)
\delta(z), \label{eq:F19}
\end{align}
and 
\begin{align}
\partial_z \bU^\ep(\bx,z) = \frac{i k}{\ep} \left[ {\bf I} +
  \frac{1}{k^2} \nabla^\perp\big(\nabla^\perp \cdot
  \big)\right]\bE^\ep(\bx,z) + \frac{ik  {\alpha}}{\ep^{1/2}} \nu
\left(\gamma\bx,\gamma\frac{z}{\ep}\right) \bE^\ep(\bx,z) 
 - \bJ\big(\gamma_{_\mJ}  {\bx} \big) \delta(z), \label{eq:F20}
\end{align}
for $0 \le z \le L$ and ${\bf I}$ the $2 \times 2$ identity matrix. For ranges
$z < 0$ and $z > L$ the equations are simpler, as all the terms
involving the process $\nu$ vanish. The fields $\bE^\ep$ and $\bU^\ep$
are approximations of $\bE$ and $\bU$ for $\ep \ll 1$, and equations
(\ref{eq:F19}-\ref{eq:F20}) follow by neglecting terms of order
$\ep^{1/2}$ and higher in (\ref{eq:F10}-\ref{eq:F11}). 

\section{Wave decomposition}
\label{sect:WD}
Because the interaction of the waves with the random medium depends on
the direction of propagation, we decompose $\bE^\ep$ and $\bU^\ep$
over plane waves, using the Fourier transform with respect to the
transverse coordinates $\bx$.  We denote this transform with
a hat, to distinguish it from the Fourier transform with respect to the three dimensional 
vector $\vx = (\bx,z)$. We have 
\begin{equation}
\what \bE^\ep(k \bka,z) = \int_{\mathbb{R}^2} d \bx \, \bE^\ep(\bx,z) e^{-
    i k \bka \cdot \bx}, 
\label{eq:F21} 
\end{equation}
and similar for $\what \bU^\ep(k \bka,z)$, where $\bka$ is the scaled
wave vector.  The inverse transform is
\begin{equation}
\bE^\ep(\bx,z) = \int_{\mathbb{R}^2} \frac{d(k \bka)}{(2
  \pi)^2}\, \what \bE^\ep (k\bka,z) e^{i k \bka \cdot \bx}.
\label{eq:F22}
\end{equation}
We begin in Section \ref{sect:WDH} with the decomposition in
homogeneous media and then consider random media in Section
\ref{sect:WDR}. We state the energy conservation in Section
\ref{sect:WD2EC} and interpret the wave decomposition in terms of
transverse electric and magnetic plane waves in Section
\ref{sect:WD3}.

\subsection{Homogeneous media}
\label{sect:WDH}
The transformed fields in homogeneous media are $\what
\bE_o^\ep(k\bka,z)$ and $\what \bU_o^\ep(k\bka,z)$. They satisfy the
system of ordinary differential equations
\begin{align}
\partial_z \begin{pmatrix} \what \bE_o^\ep(k\bka,z) \\ \what
  \bU_o^\ep(k\bka,z) \end{pmatrix} = \frac{i k}{\ep} {\bf
  M}(\bka) \begin{pmatrix} \what \bE_o^\ep(k\bka,z) \\ \what
  \bU_o^\ep(k\bka,z) \end{pmatrix} + \frac{\delta(z)}{\gamma_{_\mJ} ^2}
 \begin{pmatrix} \bka \what \mJ_z\Big(\frac{k\bka}{\gamma_{_\mJ} }
  \Big) \\- \what \bJ\Big(\frac{k\bka}{\gamma_{_\mJ} }
  \Big) \end{pmatrix}, \label{eq:WD1}
\end{align}
with $4 \times 4$ matrix 
\begin{equation}
{\bf M}(\bka) = \begin{pmatrix} {\bf 0} & {\bf I} - \bka \otimes \bka \\
{\bf I} - \bka^\perp \otimes \bka^\perp & {\bf 0} \end{pmatrix},
\label{eq:WD2}
\end{equation}
and symbol $\otimes$ denoting vector outer product. 
The Fourier transform of the current source $\vbJ = (\bJ,J_z)$
is 
\begin{equation}
\what \bJ(k\bka) = 
\int_{\mathbb{R}^2} d \bx \, \bJ(\bx) e^{ - i k \bka \cdot \bx}, \qquad 
\what J_z(k\bka) =
\int_{\mathbb{R}^2} d \bx \, J_z(\bx) e^{ - i k \bka \cdot \bx}.
\end{equation}
Naturally, the
fields are outgoing and bounded away from the source.
Moreover, the Fourier transform $(\what \bJ(k\bka), \what J_z(k\bka))$ of the current source 
is supported at $ |\bka| < 1 / k$. Therefore, the excitation in equations
(\ref{eq:WD1}) is supported at transverse wave vectors satisfying
\begin{equation}
|\bka| \le \frac{\gamma_{_\mJ}}{k} < 1 .
\label{eq:asGJSC}
\end{equation}

The plane wave decomposition is based on the diagonalization of the
matrix ${\bf M}(\bka)$. This has two double eigenvalues $\pm
\beta(\bka)$ and the eigenvectors are
\begin{equation}
\bpsi_\pm(\bka) = \begin{pmatrix} \pm \beta^{1/2}(\bka)
  \frac{\bka}{|\bka|} \\ \beta^{-1/2}(\bka) \frac{\bka}{|\bka|} 
\end{pmatrix} \quad \mbox{and} \quad  
\bpsi^\perp_\pm(\bka) = \begin{pmatrix} \beta^{-1/2}(\bka)
  \frac{\bka^\perp}{|\bka|} \\ \pm \beta^{1/2}(\bka) \frac{\bka^\perp}{|\bka|}
\end{pmatrix}. \label{eq:WD4}
\end{equation}
They are linearly independent when $|\bka| \neq 1$, so they form a
basis of $\mathbb{R}^4$ in which we can expand the solution of
(\ref{eq:WD1}),
\begin{align}
\begin{pmatrix} \what \bE_o^\ep(k\bka,z) \\ \what
  \bU_o^\ep(k\bka,z) \end{pmatrix} = \upsilon_o(\bka,z) \psi_+(\bka) +
\upsilon_o^\perp(\bka,z) \psi_+^\perp(\bka) + \chi_o(\bka,z) \psi_{-}(\bka)
+ \chi_o^\perp(\bka,z) \psi_{-}^\perp(\bka). \label{eq:WD5}
\end{align}
Because the frequency is fixed, we do not include the wavenumber $k$ in the argument 
of the coefficients in the expansion. 
Equation \eqref{eq:WD5} is the Helmholtz decomposition of the electric field and the
rotated magnetic field, because the components along $\bka$ correspond
to curl free vector fields in the $\bx$ domain, and the components
along $\bka^\perp$ to divergence free fields. 

The expressions of the coefficients in (\ref{eq:WD5}) are obtained by
substituting into (\ref{eq:WD1}) and using the linear independence of
the eigenvectors. We obtain that
\[
\upsilon_o(\bka,z) = a_o(\bka) e^{\frac{ik}{\ep} \beta(\bka)z}, \qquad
\chi_o(\bka,z) = b_o(\bka) e^{-\frac{ik}{\ep} \beta(\bka)z},
\]
and similar for $\upsilon_o^\perp$ and $\chi_o^\perp$.  Consequently, the
electric field is given by
\begin{align}
\what \bE_o^\ep(k\bka,z) = 1_{(0,\infty)}(z) \Big[a_o(\bka)
  \beta^{1/2}(\bka) \frac{\bka}{|\bka|} + a_o^\perp(\bka)
  \beta^{-1/2}(\bka) \frac{\bka^\perp}{|\bka|}\Big] e^{\frac{i k}{\ep}
  \beta(\bka) z}   \nonumber \\ 
  + 1_{(-\infty,0)}(z)
\Big[-b_o(\bka) \beta^{1/2}(\bka) \frac{\bka}{|\bka|} +
  b_o^\perp(\bka) \beta^{-1/2}(\bka)
  \frac{\bka^\perp}{|\bka|}\Big] e^{-\frac{i k}{\ep} \beta(\bka) z},
\label{eq:WD6}
\end{align}
and the rotated magnetic field is 
\begin{align}
\what \bU_o^\ep(k\bka,z) = 1_{(0,\infty)}(z) \Big[a_o(\bka)
  \beta^{-1/2}(\bka) \frac{\bka}{|\bka|} + a_o^\perp(\bka)
  \beta^{1/2}(\bka) \frac{\bka^\perp}{|\bka|}\Big] e^{\frac{i k}{\ep}
  \beta(\bka) z} \nonumber \\ 
  + 1_{(-\infty,0)}(z)
\Big[b_o(\bka) \beta^{-1/2}(\bka) \frac{\bka}{|\bka|} -
  b_o^\perp(\bka) \beta^{1/2}(\bka) \frac{\bka^\perp}{|\bka|}\Big]
e^{-\frac{i k}{\ep} \beta(\bka) z}. \label{eq:WD7}
\end{align}
Recalling the Fourier transform (\ref{eq:F21}), we see that this is a
plane wave decomposition with wave vectors $\big(\bka,\pm
\beta(\bka)\big)$ multiplying $(\bx,z/\ep)$ in the phases of the
exponentials. The plus sign corresponds to forward going waves, along
the positive range axis, and the negative sign to backward going
waves.

In (\ref{eq:WD6}) we used the radiation conditions, so that the waves
are outgoing and bounded, and the amplitudes are determined by the
jump conditions at the source
\begin{align*}
\what \bE_o^\ep(k\bka, 0+)-\what \bE_o^\ep(k\bka,0-) =
\frac{\bka}{\gamma_{_\mJ} ^2} \what \mJ_z\left(\frac{k\bka}{\gamma_{_\mJ} }\right), \qquad 
\what \bU_o^\ep(k\bka, 0+)-\what \bU_o^\ep(k\bka,0-) &=
-\frac{1}{\gamma_{_\mJ} ^2} \what \bJ\left(\frac{k\bka}{\gamma_{_\mJ} }\right).
\end{align*}
We obtain that
\begin{align}
a_o(\bka) &= \frac{1}{2 \gamma_{_\mJ} ^2} \left[ \beta^{-1/2}(\bka)
  |\bka| \what \mJ_z \left(\frac{k\bka}{\gamma_{_\mJ} }\right) -
  \beta^{1/2}(\bka) \frac{\bka}{|\bka|}\cdot \what
  \bJ\left(\frac{k\bka}{\gamma_{_\mJ} }\right)
  \right], \label{eq:WD8} \qquad  a_o^\perp(\bka) = -\frac{1}{2
  \gamma_{_\mJ} ^2} \beta^{-1/2}(\bka) \frac{\bka^\perp}{|\bka|} \cdot
\what \bJ\left(\frac{k\bka}{\gamma_{_\mJ} }\right)
,\\ b_o(\bka) &= \frac{1}{2 \gamma_{_\mJ} ^2} \left[
  \beta^{-1/2}(\bka) |\bka| \what \mJ_z
  \left(\frac{k\bka}{\gamma_{_\mJ} }\right) + \beta^{1/2}(\bka)
  \frac{\bka}{|\bka|}\cdot \what \bJ\left(\frac{k\bka}{\gamma_{_\mJ}
  }\right) \right], \qquad  b_o^\perp(\bka) =
a_o^\perp(\bka). \label{eq:WD8bp}
\end{align}
Note that equations (\ref{eq:WD8}) written in dimensional units are
the same as (\ref{eq:MR4})-(\ref{eq:MR5}). Moreover, 
as noted above, the excitation in equations
(\ref{eq:WD1}) is supported at transverse wave vectors satisfying (\ref{eq:asGJSC}).
Thus, there are no evanescent waves in the decompositions 
(\ref{eq:WD6}-\ref{eq:WD7}).

\subsection{Random media}
\label{sect:WDR}
The Fourier transforms of the wave fields in the random medium satisfy
the system of equations
\begin{align}
\partial_z \begin{pmatrix} \what \bE^\ep(k\bka,z) \\ \what
  \bU^\ep(k\bka,z) \end{pmatrix}= \frac{i k}{\ep} {\bf M}(\bka)
\begin{pmatrix} \what \bE^\ep(k\bka,z) \\ \what
  \bU^\ep(k\bka,z) \end{pmatrix} +  1_{(0,L)}(z)
\, \left[ \mathcal{M}^\ep
\begin{pmatrix} \what \bE^\ep \\ \what \bU^\ep \end{pmatrix} \right]  (k\bka,z),
\label{eq:WD10}
\end{align}
derived from (\ref{eq:F19}-\ref{eq:F20}), with radiation conditions
at $z < 0$ and $z > L$, and source conditions at $z = 0$. The leading
order term in the right hand side involves the same matrix ${\bf
  M}(\bka)$ as in the homogeneous medium. The perturbation due to the
random medium is given by the integral operator $\mathcal{M}^\ep$
defined by
\begin{align}
 \left[\mathcal{M}^\ep  \begin{pmatrix}
      \what \bE^\ep \\ \what \bU^\ep \end{pmatrix}
\right](k\bka,z) =&
  \frac{ik {\alpha}}{\ep^{1/2}\gamma^2} \int \frac{d (k
    \bka')}{(2 \pi)^2} \, \what \nu \left(\frac{k(\bka
    -\bka')}{\gamma},\gamma\frac{z}{\ep}\right) 
    \begin{pmatrix} {\bf 0} & \bka \otimes \bka' \\ {\bf I} &
    {\bf 0} \end{pmatrix}  \begin{pmatrix} \what
    \bE^\ep(k\bka',z) \\ \what \bU^\ep(k\bka',z) \end{pmatrix} 
  \nonumber \\& 
  -\frac{i k {\alpha}^2}{\gamma^2} \int \frac{d (k
    \bka')}{(2 \pi)^2}\, \what{\nu^2} \left(\frac{k(\bka -
    \bka')}{\gamma},\gamma\frac{z}{\ep}\right)
 \begin{pmatrix} {\bf 0} & \bka \otimes \bka' \\ {\bf 0} &
    {\bf 0} \end{pmatrix} \begin{pmatrix} \what
    \bE^\ep(k\bka',z) \\ \what \bU^\ep(k\bka',z) \end{pmatrix},
  \label{eq:WD11}
\end{align}
where $\what \nu$ and $\what{\nu^2}$ are the Fourier transforms of $\nu$ and $\nu^2$ with respect to the
first argument.

We use a similar decomposition of the waves as in the homogeneous
medium, based on the same eigenvector basis
$\{\bpsi_\pm(\bka),\bpsi_\pm^\perp(\bka)\}$ that diagonalizes the
matrix ${\bf M}(\bka)$ in the leading term of \eqref{eq:WD10}. We
obtain that
\begin{align}
\begin{pmatrix} \what \bE^\ep(k\bka,z) \\ \what
  \bU^\ep(k\bka,z) \end{pmatrix} =& \Big[ a^\ep(\bka,z) \bpsi_+(\bka)
  + a^{\ep,\perp}(\bka,z) \bpsi^\perp_+(\bka) \Big] e^{\frac{i k}{\ep}
  \beta(\bka) z} + \Big[b^\ep(\bka,z) \bpsi_{-}(\bka) +
  b^{\ep,\perp}(\bka,z) \bpsi^\perp_{-}(\bka) \Big] e^{-\frac{i k}{\ep}
  \beta(\bka) z},
\label{eq:WD14}
\end{align}
where the amplitudes $a^\ep(\bka,z)$, $a^{\ep,\perp}(\bka,z)$,
$b^\ep(\bka,z)$, and $b^{\ep,\perp}(\bka,z)$ are no longer constants
determined by the current source density, but random fields. We 
keep all the components of the waves, forward and backward going,
because of scattering in the random medium.

The radiation conditions are
\begin{align}
  a^\ep(\bka,z) &= a^{\ep,\perp}(\bka,z) = 0, \quad \mbox{if} ~~ z <
  0, \quad \mbox{and} \quad b^\ep(\bka,z) = b^{\ep,\perp}(\bka,z) = 0, \quad
  \mbox{if}~~ z\ge L, \label{eq:WD15}
\end{align}
because the medium is homogeneous outside the range interval $(0,L)$.
This also implies that
\begin{align}
  a^\ep(\bka,z) = a^\ep(\bka,L), ~~\mbox{and} ~ ~ 
  a^{\ep,\perp}(\bka,z) = a^{\ep,\perp}(\bka,L), \quad \mbox{if} ~ z > L,
\end{align}
and
\begin{align}
  b^\ep(\bka,z) = b^\ep(\bka,0-) ~~\mbox{and}~~
  b^{\ep,\perp}(\bka,z) = b^{\ep,\perp}(\bka,0-), \quad \mbox{if}~ z < 0.
\end{align}
The excitation comes from the jump conditions at the source
\begin{align*}
\what \bE^\ep(k\bka, 0+)-\what \bE^\ep(k\bka,0-) &=
\frac{\bka}{\gamma_{_\mJ}^2} \what \mJ_z
\left(\frac{k\bka}{\gamma_{_\mJ} } \right), \qquad \what \bU^\ep(k\bka,
0+)-\what \bU^\ep(k\bka,0-) = -\frac{1}{\gamma_{_\mJ}^2} \what \bJ
\left( \frac{k\bka}{\gamma_{_\mJ} } \right),
\end{align*}
which give the following initial conditions at $z = 0$,
\begin{align}
  a^\ep(\bka,0+) &= a_o(\bka), \quad a^{\ep,\perp}(\bka,0+) =
  a_o^\perp(\bka), \label{eq:WD16}
\end{align}
and
\begin{align}
  b^\ep(\bka,0-) &= b_o(\bka) + b^\ep(\bka,0+), \quad b^{\ep,\perp}(\bka,0-)=
  b_o^\perp(\bka) + b^{\ep,\perp}(\bka,0+). \label{eq:WD17}
\end{align}
Equations (\ref{eq:WD16}) say that the forward going waves leaving the
source are the same as in the homogeneous medium. This is physical,
because the waves must travel a long distance before they are affected
by the small fluctuations of the electric permittivity.  Equations
(\ref{eq:WD17}) say that the waves at $z < 0$ are given by the
superposition of those emitted by the source, modeled by $b_o(\bka)$
and $b_o^\perp(\bka)$, and the waves backscattered by the random
medium, modeled by $b^\ep(\bka,0+)$ and $b^{\ep,\perp}(\bka,0+)$. 

To determine the amplitudes in the random medium, we substitute
equations (\ref{eq:WD14}) into (\ref{eq:WD10}), and use the linear
independence of the eigenvectors $\bpsi_\pm(\bka)$ and
$\bpsi_\pm^\perp(\bka)$. If we let 
\begin{eqnarray}
 {\itbf Y}^\ep (\bka,z)=
\begin{pmatrix}
{a}^\ep(\bka,z) \\
{a}^{\ep,\perp}(\bka,z) \\
{b}^\ep(\bka,z) \\
{b}^{\ep,\perp}(\bka,z)
\end{pmatrix},
\end{eqnarray}
we obtain that 
\begin{align}
\nonumber \frac{\partial {\itbf Y}^\ep}{\partial z} (\bka,z)= \frac{i
  k {\alpha}}{2 \gamma^2 \ep^{1/2}} \int \frac{d (k
  \bka')}{(2\pi)^2}\, \widehat{\nu} \left(
\frac{k(\bka-\bka')}{\gamma},\gamma \frac{z}{\ep}\right) {\bf
  F}\left(\bka,\bka',\frac{z}{\ep} \right) {\itbf Y}^\ep(\bka',z)  \\
+ \frac{i k {\alpha}^2}{2\gamma^2 } \int \frac{d (k
  \bka')}{(2\pi)^2} \, \widehat{\nu^2} \left(
\frac{k(\bka-\bka')}{\gamma},\gamma\frac{z}{\ep}\right) {\bf
  G}\left(\bka,\bka',\frac{z}{\ep} \right) {\itbf Y}^\ep(\bka',z),
\label{eq:twopointbvp}
\end{align}
in $z \in (0,L)$, with boundary conditions (\ref{eq:WD15}) and
(\ref{eq:WD16}). We are interested in the propagating waves,
corresponding to $|\bka| < 1$ in (\ref{eq:twopointbvp}), and we
explain in Section \ref{sect:DL} that in our regime the evanescent
waves may be neglected. This means that we can restrict the integral
in (\ref{eq:twopointbvp}) to vectors $\bka'$ satisfying $|\bka'| < 1$.

The mode amplitudes are coupled in equations (\ref{eq:twopointbvp}) by
the $4 \times 4$ complex matrices ${\bf F}(\bka,\bka',\zeta)$ and
${\bf G}(\bka,\bka',\zeta)$, with block structure
\begin{equation}
{\bf F} (\bka,\bka',\zeta ) =
\begin{pmatrix}
{\bf F}^{aa} (\bka,\bka',\zeta )  & {\bf F}^{ab} (\bka,\bka',\zeta ) \\
{\bf F}^{ba} (\bka,\bka',\zeta )  & {\bf F}^{bb} (\bka,\bka',\zeta ) 
\end{pmatrix},
\end{equation}
and 
\begin{equation}
{\bf G} (\bka,\bka',\zeta ) =
\begin{pmatrix}
{\bf G}^{aa} (\bka,\bka',\zeta )  & {\bf G}^{ab} (\bka,\bka',\zeta ) \\
{\bf G}^{ba} (\bka,\bka',\zeta )  & {\bf G}^{bb} (\bka,\bka',\zeta )
\end{pmatrix},
\end{equation}
where the superscripts of the $2\times2$ blocks indicate which types
of waves they couple.  The dependence on $\zeta$ of the blocks in
${\bf F}$ can be factored as
\begin{align}
{\bf F}^{aa} (\bka,\bka',\zeta ) &= \bGamma^{aa}(\bka,\bka') e^{ i k
  \big[\beta(\bka')-\beta(\bka) \big] \zeta} ,
  \label{eq:Fzetaaa}
   \\ {\bf F}^{bb}
(\bka,\bka',\zeta ) &= \bGamma^{bb}(\bka,\bka') e^{i
  k\big[-\beta(\bka')+\beta(\bka) \big] \zeta} ,\
  \label{eq:Fzetabb}
   \\ {\bf  F}^{ab} (\bka,\bka',\zeta ) &= \bGamma^{ab}(\bka,\bka') e^{-i
  k\big[\beta(\bka')+\beta(\bka) \big] \zeta} , \label{eq:Fzetaab} \\ 
  {\bf F}^{ba} (\bka,\bka',\zeta ) &= \bGamma^{ba}(\bka,\bka') e^{ i
  k\big[\beta(\bka')+\beta(\bka) \big] \zeta} ,\ \label{eq:Fzetaba}
\end{align}
with $2\times 2$ real-valued matrices 
\begin{equation}
\bGamma^{aa}(\bka,\bka') = \begin{pmatrix} \frac{|\bka|
    |\bka'|}{\sqrt{\beta(\bka) \beta(\bka')}} + \frac{\bka}{|\bka|}
  \cdot \frac{\bka'}{|\bka'|}  \sqrt{\beta(\bka)\beta(\bka')} &
  \frac{\bka}{|\bka|} \cdot \frac{\bka'^{\perp}}{|\bka'|}
  \sqrt{\frac{\beta(\bka)}{\beta(\bka')}} \\ \frac{\bka^\perp}{|\bka|}
  \cdot \frac{\bka'}{|\bka'|} \sqrt{\frac{\beta(\bka')}{\beta(\bka)}}
  & \frac{\bka}{|\bka|} \cdot \frac{\bka'}{|\bka'|}
  \frac{1}{\sqrt{\beta(\bka')\beta(\bka)}}
\label{eq:DefGamma}
\end{pmatrix},
\end{equation}
and 
\begin{equation}
\bGamma^{bb}(\bka,\bka') = \begin{pmatrix} -\frac{|\bka|
    |\bka'|}{\sqrt{\beta(\bka) \beta(\bka')}} - \frac{\bka}{|\bka|}
  \cdot \frac{\bka'}{|\bka'|}  \sqrt{\beta(\bka)\beta(\bka')} &
  \frac{\bka}{|\bka|} \cdot \frac{\bka'^{\perp}}{|\bka'|}
  \sqrt{\frac{\beta(\bka)}{\beta(\bka')}} \\ \frac{\bka^\perp}{|\bka|}
  \cdot \frac{\bka'}{|\bka'|} \sqrt{\frac{\beta(\bka')}{\beta(\bka)}}
  & -\frac{\bka}{|\bka|} \cdot \frac{\bka'}{|\bka'|}
  \frac{1}{\sqrt{\beta(\bka')\beta(\bka)}}
\label{eq:DefGammab}
\end{pmatrix},
\end{equation}
on the diagonal. Note that $\bGamma^{aa}$ is the same as the matrix
$\bGamma$ in (\ref{eq:MR1}). The off-diagonal matrices are
\begin{equation}
\bGamma^{ab}(\bka,\bka') = \begin{pmatrix} \frac{|\bka|
    |\bka'|}{\sqrt{\beta(\bka) \beta(\bka')}} - \frac{\bka}{|\bka|}
  \cdot \frac{\bka'}{|\bka'|}  \sqrt{\beta(\bka)\beta(\bka')} &
  \frac{\bka}{|\bka|} \cdot \frac{\bka'^{\perp}}{|\bka'|}
  \sqrt{\frac{\beta(\bka)}{\beta(\bka')}} \\ -\frac{\bka^\perp}{|\bka|}
  \cdot \frac{\bka'}{|\bka'|} \sqrt{\frac{\beta(\bka')}{\beta(\bka)}}
  & \frac{\bka}{|\bka|} \cdot \frac{\bka'}{|\bka'|}
  \frac{1}{\sqrt{\beta(\bka')\beta(\bka)}}
\label{eq:DefGammaAB}
\end{pmatrix},
\end{equation}
and 
\begin{equation}
\bGamma^{ba}(\bka,\bka') = \begin{pmatrix} -\frac{|\bka|
    |\bka'|}{\sqrt{\beta(\bka) \beta(\bka')}} + \frac{\bka}{|\bka|}
  \cdot \frac{\bka'}{|\bka'|} \sqrt{\beta(\bka)\beta(\bka')} &
  \frac{\bka}{|\bka|} \cdot \frac{\bka'^{\perp}}{|\bka'|}
  \sqrt{\frac{\beta(\bka)}{\beta(\bka')}}
  \\ -\frac{\bka^\perp}{|\bka|} \cdot \frac{\bka'}{|\bka'|}
  \sqrt{\frac{\beta(\bka')}{\beta(\bka)}} & -\frac{\bka}{|\bka|} \cdot
  \frac{\bka'}{|\bka'|} \frac{1}{\sqrt{\beta(\bka')\beta(\bka)}}
\label{eq:DefGammaBA}
\end{pmatrix}.
\end{equation}
The diagonal blocks in ${\bf G}$ are
\begin{eqnarray}
  \label{eq:DefGaa}
{\bf G}^{aa}(\bka,\bka',\zeta) &=& \begin{pmatrix} -\frac{|\bka|
    |\bka'|}{\sqrt{\beta(\bka) \beta(\bka')}} & 0 \\ 0 &
  0 \end{pmatrix} e^{ i k \big[\beta(\bka')-\beta(\bka) \big] \zeta},\\ 
  {\bf G}^{bb}(\bka,\bka',\zeta) &=& - \overline{{\bf G}^{aa}(\bka,\bka',\zeta)}, 
  \label{eq:DefGbb}
\end{eqnarray}
and the off-diagonal blocks are similar
\begin{eqnarray}
{\bf G}^{ab}(\bka,\bka',\zeta) &=& \begin{pmatrix} -\frac{|\bka|
    |\bka'|}{\sqrt{\beta(\bka) \beta(\bka')}} & 0 \\ 0 &
  0 \end{pmatrix} e^{- i k \big[\beta(\bka')+\beta(\bka) \big] \zeta},\\ 
  {\bf G}^{ba}(\bka,\bka',\zeta) &=& - \overline{{\bf
    G}^{ab}(\bka,\bka',\zeta)}.
\end{eqnarray}

\subsection{Energy conservation}
\label{sect:WD2EC}
The diagonal blocks of coupling matrices ${\bf F}$ and ${\bf G}$
satisfy the symmetry relations
\begin{align}
\label{eq:SYM1}
{\bf F}^{aa} (\bka,\bka',\zeta ) & = {\bf
  F}^{aa}(\bka',\bka,\zeta)^\dagger, \qquad {\bf G}^{aa}
(\bka,\bka',\zeta ) = {\bf G}^{aa}(\bka',\bka,\zeta)^\dagger , \\ 
{\bf F}^{bb} (\bka,\bka',\zeta ) & = {\bf
  F}^{bb}(\bka',\bka,\zeta)^\dagger, \qquad {\bf G}^{bb}
(\bka,\bka',\zeta ) = {\bf G}^{bb}(\bka',\bka,\zeta)^\dagger.
\end{align}
Using these in (\ref{eq:twopointbvp}) we obtain the conservation
identity
\begin{equation}
\int_{|\bka|<1} \frac{d (k \bka)}{(2 \pi)^2} \Big[ |a^\ep(\bka,z)|^2 +
  |a^{\ep,\perp}(\bka,z)|^2 - |b^\ep(\bka,z)|^2 -
  |b^{\ep,\perp}(\bka,z)|^2\Big] = ~\mbox{constant},
\label{eq:ECId}
\end{equation}
independent of $z$. 

Now let us recall that $\frac{1}{2} \mbox{Re}\big[ \vbE^\ep \times
  \vbH^\ep\big] $ is the time average of the Poynting vector of the
time harmonic electromagnetic wave \cite{jackson1999classical}. Its
$z$-component is
\begin{equation*}
 -\frac{1}{2} \mbox{Re}\Big[ \bE^\ep(\bx,z) \cdot
   \bH^{\ep,\perp}(\bx,z) \Big] = \frac{\zeta_o^{-1}}{2}
 \mbox{Re}\Big[ \bE^\ep(\bx,z) \cdot \bU^\ep(\bx,z) \Big],
\end{equation*}
where we used (\ref{eq:F9}). The integral of this is the energy
flux which must be conserved
\begin{equation}
\int_{\mathbb{R}^2} d \bx \,\mbox{Re}\Big[ \bE^\ep(\bx,z) \cdot
  \bU^\ep(\bx,z) \Big]  = \mbox{constant}.
\label{eq:ECId2}
\end{equation}
Our identity (\ref{eq:ECId}) is the statement of this energy
conservation, written in terms of the wave amplitudes. It follows from
(\ref{eq:ECId2}) by substituting equations (\ref{eq:WD14}) and using
Plancherel's theorem.

\subsection{Transverse electric and transverse magnetic waves}
\label{sect:WD3}
To interpret the wave decomposition (\ref{eq:WD14}) as an expansion of
the electromagnetic wave field in transverse electric and magnetic
waves, we recall the definitions (\ref{eq:F5})-(\ref{eq:F6}) of
the longitudinal electric and magnetic field. After the scaling, and
in the Fourier domain, these become
\begin{align}
  \what H_z^\ep(k\bka,z) &= \zeta_o^{-1} \bka^\perp \cdot \what
  \bE^\ep(k\bka,z) , \label{eq:HZeq} \\ \what E_z^\ep(k\bka,z) & = -
  \zeta_o \bka^\perp \cdot \what \bH^\ep(k\bka,z) = - \bka \cdot \what
  \bU^\ep(k\bka,z). \label{eq:Ezeq}
\end{align}
Then, using (\ref{eq:WD14}), the definition (\ref{eq:WD4}) of the
eigenfunctions, and the inverse Fourier transform (\ref{eq:F22}), we
obtain that the electric field $\vbE^\ep = ( \bE^\ep,\mE_z^\ep)$ is given by 
\begin{align}
\nonumber
  \vbE^\ep \big(\frac{\bx}{\ep},z\big) =  \int_{|\bka|<1} \frac{d(k \bka)}{(2 \pi)^2} \,&
  \beta^{-\frac{1}{2}}(\bka)\Big\{ \big[a^\ep(\bka,z) \, \vu(\bka) +
    a^{\ep,\perp}(\bka,z) \, \vu^\perp(\bka) \big] e^{i \frac{k}{\ep}
    \vec{\bka} \cdot \vx} \\ & +\big[-b^\ep(\bka,z) \,
    \vu^{\,-}(\bka)  + b^{\ep,\perp}(\bka,z) \, \vu^\perp(\bka) \big] e^{i
    \frac{k}{\ep} \vec{\bka}^{\,-} \cdot \vx}\Big\}, \label{eq:EField}
\end{align}
and the magnetic field is 
\begin{align}
\nonumber
  \vbH^\ep \big( \frac{\bx}{\ep},z\big)
 =  \zeta_o^{-1} \int_{|\bka|<1} \frac{d(k \bka)}{(2
   \pi)^2} &\beta^{-\frac{1}{2}}(\bka) \Big\{ \big[a^\ep(\bka,z) \,
   \vu^\perp(\bka) - a^{\ep,\perp}(\bka,z) \, \vu(\bka) \big] e^{i
   \frac{k}{\ep} \vec{\bka} \cdot \vx} \\
   & +\big[b^\ep(\bka,z) \,
   \vu^\perp(\bka) + b^{\ep,\perp}(\bka,z) \, \vu^{\,-}(\bka) \big] e^{i
   \frac{k}{\ep} \vec{\bka}^{\,-} \cdot \vx}\Big\}, \label{eq:MField}
\end{align}
for $\vx = (\bx,z)$.  Here we used the wave vectors $\vec{\bka} =
\big(\bka, \beta(\bka)\big)$ and $ \vec{\bka}^{\,-} = \big(\bka,
-\beta(\bka) \big)$ of the forward and backward going waves, and
\begin{equation*}
\vu(\bka) = \left( \beta(\bka) \frac{\bka}{|\bka|}, -|\bka| \right),
\quad \vu^{\,-}(\bka) = \left( \beta(\bka) \frac{\bka}{|\bka|},
|\bka|\right), \quad \vu^\perp(\bka)= \left(
  \frac{\bka^\perp}{|\bka|}, 0\right).
\end{equation*} 
The triplets $\{ \vu(\bka),\vu^\perp(\bka),\vka \}$ and
$\{ \vu^{\,-}(\bka), \vu^\perp(\bka),\vka^{\,-}\}$ are orthonormal
bases of $\mathbb{R}^3$, when $|\bka| < 1$.  Thus, (\ref{eq:EField})
and (\ref{eq:MField}) are decompositions in transverse waves, which
are orthogonal to the direction of propagation along the wave vectors
$\vka$ and $\vec{\bka}^-$, respectively. The amplitudes
$a^\ep(\bka,z)$ and $b^\ep(\bka,z)$ are for the forward and backward
transverse magnetic plane waves, and $a^{\ep,\perp}(\bka,z)$ and
$b^{\ep,\perp}(\bka,z)$ are the amplitudes of transverse electric
plane waves. Note that the forward going part, written in dimensional units,
is the same as in (\ref{eq:EField_SR}) and (\ref{eq:MField_SR}).

\section{The Markov limit}
\label{sect:DL}
Here we describe the $\ep \to 0$ limit of the random mode amplitudes
satisfying the two-point linear boundary value problem
(\ref{eq:twopointbvp}) with boundary conditions
(\ref{eq:WD15})-(\ref{eq:WD16}). We begin in Section \ref{sect:propag}
by writing the solution in terms of propagator matrices. Then we
explain in Section \ref{sect:FScA} why we can neglect the backward and
evanescent waves. The limit of the forward going wave amplitudes is
given in Section \ref{sect:DiffA}.

\subsection{Propagator matrices}
\label{sect:propag}
The $4 \times 4$ propagator matrices ${\bf P}^\ep(\bka, z;\bka_o)$ are
solutions of the initial value problem
\begin{align}
\frac{\partial {\bf P}^\ep(\bka,z;\bka_o)}{\partial z} \nonumber
=\frac{i k {\alpha}}{2 \gamma^2 \ep^{1/2}} \int_{|\bka'|<1}
\frac{d(k\bka')}{(2\pi)^2}\, \widehat{\nu} \left(
\frac{k(\bka-\bka')}{\gamma},\gamma \frac{z}{\ep}\right) {\bf
  F}\Big(\bka,\bka',\frac{z}{\ep} \Big) {\bf P}^\ep(\bka',z;\bka_o)
\\ + \frac{i k {\alpha}^2}{2 \gamma^2} \int_{|\bka'|<1}
\frac{d(k\bka')}{(2 \pi)^2} \, \widehat{\nu^2} \left(
\frac{k(\bka-\bka')}{\gamma}, \gamma \frac{z}{\ep}\right) {\bf G}\Big(
\bka,\bka', \frac{z}{\ep} \Big) {\bf P}^\ep(\bka',z;\bka_o) ,
\label{eq:ivp}
\end{align}
 with initial condition ${\bf P}^\ep(\bka,z=0;\bka_o) =
 \delta(\bka-\bka_o) {\bf I}$, where ${\bf I}$ is the $4\times 4$
 identity matrix.  The solution of (\ref{eq:twopointbvp}) satisfies
\[ 
{\itbf Y}^\ep (\bka,z)= \int_{|\bka_o|<1} d \bka_o\, {\bf
  P}^\ep(\bka,z;\bka_o) {\itbf Y}^\ep (\bka_o,0) ,
\]
for any $z \in [0,L]$. In particular, at $z=L$ and using
boundary conditions (\ref{eq:WD15})-(\ref{eq:WD16}), 
\begin{equation}
\begin{pmatrix}
{a}^\ep(\bka,L) \\
{a}^{\ep,\perp}(\bka,L) \\
0\\ 0
\end{pmatrix}
= \int_{|\bka_o|<1} d \bka_o\, {\bf P}^\ep(\bka,L;\bka_o)
\begin{pmatrix}
{a}_o(\bka_o) \\
{a}^{\perp}_o(\bka_o) \\
{b}^\ep(\bka_o,0) \\
{b}^{\ep,\perp}(\bka_o,0) 
\end{pmatrix}
. \label{eq:ivp1}
\end{equation}

\subsection{The forward scattering approximation}
\label{sect:FScA}
The $\ep \to 0$ limit ${\bf P}$ of ${\bf P}^\ep$ can be obtained and
identified as a Markov process, meaning that it satisfies a system
of stochastic differential equations. We refer to \cite{PK1974,PW1994} and 
%appendix
 \ref{ap:DL} for details, and state here the result.

In the limit $\ep \rightarrow 0$, the terms of order one in
(\ref{eq:ivp}) i.e., those involving the kernel ${\bf G}$, become
equal to their average with respect to the distribution of $\nu$ and
to the rapid $O(1/\ep)$ phase. We denote the average with respect to
this phase by $\left< \cdot \right>$, and obtain
\begin{align}
\frac{ik {\alpha}^2}{2 \gamma^2} \int_{|\bka'|<1} \frac{d(k
  \bka')}{(2\pi)^2}\, \EE
\left[\widehat{\nu^2}\left(\frac{k(\bka-\bka')}{\gamma},\cdot\right)\right]
\left< {\bf G}(\bka,\bka',\cdot) \right> {\bf P}(\bka',z;\bka_o)
 %\nonumber \\ 
 = \frac{i k {\alpha}^2}{2} \cR({\bf 0}) \left< {\bf
  G}(\bka,\bka,\cdot) \right> {\bf P}(\bka,z;\bka_o).
\end{align}
Due to the large phases of the off-diagonal blocks in ${\bf G}$ we
have
\[
\left< {\bf G}(\bka,\bka,\cdot) \right> =
\frac{|\bka|^2}{\beta(\bka)}
\begin{pmatrix}
-1 & 0 & 0& 0\\
0 & 0& 0& 0\\
0 & 0& 1& 0\\
0 & 0& 0& 0
\end{pmatrix},
\]
a diagonal matrix, so the order one terms in (\ref{eq:ivp}) do not
introduce any wave coupling.

The order $O(\ep^{-1/2})$ part in (\ref{eq:ivp}) i.e., involving ${\bf
  F}$, gives diffusive terms.  Let us split the propagator matrix into
four blocks:
\[{\bf P}^{\ep} =\begin{pmatrix} {\bf P}^{aa,\ep} & {\bf   P}^{ab,\ep}\\ 
{\bf P}^{ba,\ep} & {\bf P}^{bb,\ep}
\end{pmatrix} 
\]
and consider first only the propagating waves.  The stochastic
differential equations for the limit entries of ${\bf
  P}^{ab,\ep}(\bka,z;\bka_o)$ and ${\bf
  P}^{ba,\ep}(\bka,z;\bka_o)$ are coupled to the limit entries of
${\bf P}^{aa,\ep}(\bka,z;\bka_o)$  and ${\bf P}^{bb,\ep}(\bka,z;\bka_o)$ through the coefficients
\begin{align}
\int_{\mathbb{R}} d r_z \, \what \cR
\left(\frac{k(\bka-\bka')}{\gamma},\gamma r_z\right) e^{ -i k \big[
    \beta(\bka)+\beta(\bka')\big] r_z} = \int_{\RR^3} d \vr \,
     {\mathscr R} (\gamma \vr) e^{ - i {k} \big[ ( \bka-\bka')\cdot \br
         + (\beta(\bka)+\beta(\bka') ) r_z \big]} = \frac{1}{\gamma^3}
     \widetilde \cR\left(\frac{k (\vka-\vec{\bka'}^{-})}{\gamma}\right).
\label{eq:couplingterm1}
\end{align}
Here $\what \cR$ is the Fourier transform of $\cR$ with respect to the
cross-range variable, and
$\widetilde \cR$ is the three-dimensional Fourier transform \eqref{eq:PSDensity}.  
The sum $\beta(\bka) + \beta(\bka')$ in the phase of (\ref{eq:couplingterm1})
is because the phase factors in the matrices ${\bf F}^{ab}$ and  ${\bf F}^{ba}$  are
 $\pm (\beta(\bka)+\beta(\bka'))\zeta$.  The
coupling between the entries of ${\bf P}^{aa,\ep}(\bka',z;\bka_o)$ and between the entries of
${\bf P}^{bb,\ep}(\bka',z;\bka_o)$   is
through the coefficients
\begin{align}
\int_{\mathbb{R}} d r_z \, \what \cR
\left(\frac{k(\bka-\bka')}{\gamma},\gamma r_z\right) e^{- i k \big[
    \beta(\bka)-\beta(\bka')\big] r_z} = \int_{\RR^3} d \vr \,
     {\mathscr R} (\gamma \vr) e^{- i {k} \big[ ( \bka-\bka')\cdot \br
         + (\beta(\bka)-\beta(\bka') ) r_z \big]} = \frac{1}{\gamma^3}
     \widetilde \cR\left( \frac{k(\vka - \vec{\bka'})}{\gamma}\right),
\label{eq:couplingterm2}
\end{align}
because the phase factors in matrices ${\bf F}^{aa}$ and ${\bf F}^{bb}$  are $\pm (
\beta(\bka)-\beta(\bka') )\zeta $. 

We conclude that the interaction of the forward and backward going
waves depends on the decay of the power spectral density.  Recall that
 the source gives wave amplitudes supported at $|\bka| \le
\gamma_{_\mJ} /k < 1$ and $\widetilde\cR$ is negligible outside the ball 
with center at ${\bf 0}$ and radius one in $\mathbb{R}^3$. 
From \eqref{eq:defGJ} and \eqref{eq:asGJSC} it
is possible to choose some $\kappa_{_M} \in (\gamma_{_\mJ} /k,1)$ such
that $\gamma$ satisfies
\begin{equation}
\label{eq:ASG2}
\frac{2 k \beta(\kappa_{_M})}{\gamma} >  1. 
\end{equation}
Under these circumstances, for all  $\bka'$ satisfying $|\bka'| < \kappa_{_M}$,
the coupling coefficients (\ref{eq:couplingterm1})
vanish, because 
\[
\frac{k \big|\vka-\vec{\bka'}^- \big|}{\gamma} = \frac{k}{\gamma}
\left[|\bka-\bka'|^2 + \left( \beta(\bka) +
  \beta(\bka')\right)^2\right]^{1/2} \ge \frac{k}{\gamma} \left[
  \beta(\bka) + \beta(\bka')\right] \ge \frac{2k
  \beta(\kappa_{_M})}{\gamma} >  1.  
\]
This means that the
forward and backward wave amplitudes are asymptotically decoupled as
long as the energy of the wave is supported at the transverse
wave vectors $\bka$ satisfying $|\bka| \le \kappa_{_M}$.

Nevertheless, the forward going amplitudes are coupled with each
other, because 
\[
 \frac{k \big|\vka-\vec{\bka'}\big|}{\gamma} = \frac{k}{\gamma}
\left[|\bka-\bka'|^2 + \left( \beta(\bka) - \beta(\bka')\right)^2\right]^{1/2}
\]
is small for at least a subset of transverse wave vectors satisfying
$|\bka|, |\bka'| \le \kappa_{_M}$. Due to this coupling there is
diffusion of energy from the waves emitted by the source with $|\bka|
< \gamma_{_\mJ}/ k$, to waves at larger values of $|\bka|$, as
explained in Section \ref{sect:DEPOL}. This is why we take
$\kappa_M >\gamma_{_\mJ} /k $ in (\ref{eq:ASG2}). By assuming that
$a^\ep(\bka,z)$ and $a^{\ep,\perp}(\bka,z)$ are supported at $|\bka|
\le \kappa_M $ we essentially restrict $z$ by $Z_{_M}$, so that the
energy does not diffuse to waves with $|\bka| > \kappa_{_M}$ for $z\le
Z_{_M}$. Physically, it means that the three-dimensional wave vectors
$\vka$ of the forward going waves remain within a cone with opening
angle smaller than $180$ degrees.  We will see that the evolution of
the $\bka$-distribution of the wave energy is described by a radiative
transfer equation, which means that the wave energy undergoes a random
walk (or diffusion). We will also see that the diffusion coefficient is 
of the order $\gamma \alpha^2$. Thus, the 
$\bka$-distribution of the wave energy reaches $\kappa_{_M}$ after a
propagation distance of the order of $Z_{_M}$ such that $ \alpha^2
\gamma Z_{_M} = \kappa_{_M}^2$. Recalling the scaling
\eqref{def:alpha} of the standard deviation of the fluctuations and
\eqref{eq:hypscal1b}, we see that it is possible to choose $Z_{_M} =
O(1)$ (i.e., similar to $\bar{L}$ in dimensional units).

The evanescent waves can only couple with the propagating waves with
wave vectors of magnitude close to $1$. Thus, as long as the energy of
the wave is supported at wave vectors  satisfying $|\bka|< \kappa_M$,
assumption (\ref{eq:ASG2}) implies that the evanescent waves do not
get excited.

The forward scattering approximation amounts to neglecting all the
backward and evanescent waves in the system of stochastic differential
equations. It is justified in the limit $\ep \to 0$ by the decoupling
of the waves, the zero initial condition of the evanescent waves, and
the zero condition at $z=L$ of the left going wave amplitudes.

\subsection{Markov limit of the forward going mode amplitudes}
\label{sect:DiffA}
Under the forward scattering approximation, the wave amplitudes
$a^\ep(\bka,z)$ and $a^{\ep,\perp}(\bka,z)$ satisfy the initial value
problem
\begin{align}
\frac{\partial }{\partial z} \begin{pmatrix} a^\ep (\bka,z) \nonumber
    \\ a^{\ep,\perp} (\bka,z) \end{pmatrix}  =&
  \frac{i k {\alpha}}{2 \gamma^2 \ep^{1/2}} \int_{|\bka'| < 1} \frac{d (k
    \bka')}{(2 \pi)^2}\, \widehat{\nu} \left( \frac{k(\bka-\bka')}{\gamma}
  , \gamma \frac{z}{\ep}\right) {\bf
    F}^{aa}\left(\bka,\bka',\frac{z}{\ep} \right)
\begin{pmatrix}
a^\ep (\bka',z) \\ a^{\ep,\perp} (\bka',z) \end{pmatrix}
\\ &
+  \frac{i
  k {\alpha}^2}{2 \gamma^2} \int_{|\bka'| <1}\frac{d (k
  \bka')}{(2 \pi)^2}\, \widehat{\nu^2} \left(\frac{ k(\bka-\bka')}{\gamma}
, \gamma \frac{z}{\ep}\right) {\bf G}^{aa}\left(\bka,\bka',\frac{z}{\ep}
\right)
\begin{pmatrix}
a^\ep (\bka',z) \\ a^{\ep,\perp} (\bka',z) \end{pmatrix}, 
\label{eq:ivteps}
\end{align}
for $z > 0$, and the initial condition
\begin{equation}
\begin{pmatrix} a^\ep (\bka,0) 
    \\ a^{\ep,\perp} (\bka,0) \end{pmatrix} = \boldsymbol{\mathscr{A}}_o(\bka),
\label{eq:ivtepso}
\end{equation}
with $\boldsymbol{\mathscr{A}}_o$ defined in \eqref{eq:defAo}.  These
equations conserve energy, as the special structure of the kernels
${\bf F}^{aa}$ and ${\bf G}^{aa}$ described in equations
(\ref{eq:Fzetaaa}), (\ref{eq:DefGaa}), and (\ref{eq:SYM1}) implies
that for all $\ep > 0$, for all $z \geq0$,
\begin{equation}
\int_{|\bka| < 1} \frac{d (k \bka)}{(2\pi)^2} \Big[ |a^\ep(\bka,z)|^2 +
  |a^{\ep,\perp}(\bka,z)|^2 \Big] = \int_{|\bka| < 1} \frac{d (k \bka)}{(2\pi)^2}
\Big[ |a_o(\bka)|^2 + |a^\perp_o(\bka)|^2 \Big] \, .\label{eq:CENA}
\end{equation}

We refer to \ref{ap:DL} for the details on the Markov limit $\ep \to
0$. In particular, we explain there that the process
\[
\bX^\ep(z) = \begin{pmatrix} {\rm Re}\big( {a}^\ep(\bka,z) )\\ {\rm
    Im}\big( {a}^\ep(\bka,z) )\\ {\rm Re} \big(
  {a}^{\ep,\perp}(\bka,z) \\ {\rm Im} \big( {a}^{\ep,\perp}(\bka,z)
\end{pmatrix}_{\bka \in {\cal O}}, \quad \mbox{for} ~  \mathscr{O} = \big\{\bka \in \RR^2, \,
|\bka|<1\big\},
\]
converges weakly in ${\cal C}([0,L],{\cal D}')$ to a Markov process
$\bX(z)$, where ${\cal D}'$ is the space of distributions, dual to the
space $\mathcal{D}(\mathscr{O},\mathbb{R}^4)$ of infinitely
differentiable vector valued functions in $\mathbb{R}^4$, with compact
support. The generator of $\bX(z)$ is given in 
%appendix 
\ref{ap:DL}, and we denote henceforth the limit amplitudes by $a(\bka,z)$ and
$a^\perp(\bka,z)$. Their first and second moments are described in the
next two sections and are summarized in Section \ref{sect:result}.

\section{The coherent field}
\label{sect:COH}
The vector
\begin{equation}
  \boldsymbol{\mathscr{A}}(\bka,z) = \begin{pmatrix} \EE[a(\bka,z)]
    \\ \EE[a^\perp(\bka,z)] 
  \end{pmatrix}= \lim_{\ep \to 0} \EE
    \left[ \begin{pmatrix} a^\ep(\bka,z)
        \\ a^{\ep,\perp}(\bka,z) \end{pmatrix} \right],
  \label{eq:defAmean}
\end{equation}
whose components are the expectations of the limit amplitudes of the
modes at wave vector $\vka$, satisfies as explained in \ref{ap:DL} the
initial value problem
\begin{equation}
  \partial_z \boldsymbol{\mathscr{A}}(\bka,z) = {\bf Q}(\bka)
  \boldsymbol{\mathscr{A}}(\bka,z), \quad z > 0,
  \label{eq:Amean}
\end{equation}
with initial condition $\boldsymbol{\mathscr{A}}(\bka,0) =
\boldsymbol{\mathscr{A}}_o(\bka)$, and $2\times 2$ complex matrix $
           {\bf Q}(\bka)$ given by
\begin{align}
 {\bf Q}(\bka) = &-\frac{k^2 {\alpha}^2}{4 \gamma^3} \int_{|\bka'|< 1}
 \frac{d (k \bka')}{(2 \pi)^2}\, \bGamma(\bka,\bka')
 \bGamma(\bka',\bka)\int_0^\infty d \zeta \, \what
 \cR\left(\frac{k(\bka-\bka')}{\gamma}, \zeta\right) e^{
 -  i\frac{k}{\gamma} \big(\beta(\bka)-\beta(\bka')\big)\zeta}
 -\frac{i k {\alpha}^2}{2} \cR({\bf 0})
 \frac{|\bka|^2}{\beta(\bka)} \begin{pmatrix} 1 & 0 \\ 0
   &0 \end{pmatrix}.
\label{eq:C5}
\end{align}
In dimensional units this is the same as \eqref{eq:MR6}, with
$\bGamma$ given by \eqref{eq:MR1}, which equals $\bGamma^{aa}$ defined
in \eqref{eq:DefGamma}.

As stated in Section \ref{sect:result.coh}, the real part of ${\bf
  Q}(\bka)$ is a symmetric, negative definite matrix
\begin{align}
{\rm Re}[ {\bf Q}(\bka)] = -\frac{k^2 {\alpha}^2}{8 \gamma^3}\int_{|\bka'|<
  1}  \frac{d (k \bka')}{(2 \pi)^2} \, \widetilde
\cR\left(\frac{k(\vka-\vec{\bka'})}{\gamma}\right) \bGamma(\bka,\bka')
\bGamma(\bka',\bka),
  \label{eq:C5pa}
\end{align}
because $\widetilde \cR$ is non-negative and $\bGamma(\bka',\bka) =
\bGamma(\bka,\bka')^T$.  It is an effective diffusion term in
\eqref{eq:Amean}, which removes energy from the mean field and gives
it to the incoherent fluctuations. This causes the randomization or
loss of coherence of the waves.  The imaginary part of $ {\bf
  Q}(\bka)$ is the sum of two terms. The first, given by
\[
\frac{k^2 {\alpha}^2}{4 \gamma^3} \int_{|\bka'|< 1} \frac{d (k
  \bka')}{(2 \pi)^2} \int_0^\infty d \zeta \, \what
\cR\left(\frac{k(\bka-\bka')}{\gamma}, \zeta\right) \sin
\left[\frac{k}{\gamma} \big(\beta(\bka)-\beta(\bka')\big)
  \zeta\right]\bGamma(\bka,\bka') \bGamma(\bka',\bka),
\]
is an effective dispersion term, which does not remove energy from the
mean field and ensures causality\footnote{ If we write the coherent
  wave fields in the time domain, using the inverse Fourier transform
  with respect to the frequency $\om$, we obtain a causal result.}.
The second term in \eqref{eq:C5} is a homogenization effect due to the
$\what{\nu^2}$ part in \eqref{eq:ivteps}.  It accounts for the
slightly different velocities of the transverse electric and
transverse magnetic waves, with discrepancy of the order of $O(\ep)$.

\section{The transport equation}
\label{sect:TEQ}
To describe the transport of energy, we study the $z$-evolution of the
Hermitian positive definite coherence matrix
\begin{equation}
   \boldsymbol{\mathscr{P}}(\bka,z) = \lim_{\ep \to 0} \EE
   \left[ \begin{pmatrix} a^\ep(\bka,z)
       \\ a^{\ep,\perp}(\bka,z) \end{pmatrix}
 \begin{pmatrix}
   a^\ep(\bka,z) \\ a^{\ep,\perp}(\bka,z) \end{pmatrix}^\dagger \right].
  \label{eq:T1}
\end{equation}
The diagonal entries of $\boldsymbol{\mathscr{P}}$ are the mode
powers.  The evolution equation is derived as explained in
\ref{ap:DL},
\begin{align}
 \partial_z \boldsymbol{\mathscr{P}}(\bka,z) =& \mathbf{Q}(\bka)
 \boldsymbol{\mathscr{P}}(\bka,z) + \boldsymbol{\mathscr{P}}(\bka,z)
 \mathbf{Q}(\bka)^\dagger +\frac{k^2 {\alpha}^2}{4
   \gamma^3}\int_{|\bka'| < 1} \frac{d (k \bka')}{(2 \pi)^2}\,
 \bGamma(\bka,\bka') \boldsymbol{\mathscr{P}}(\bka',z)
 \bGamma(\bka',\bka) \widetilde \cR\left(\frac{k
   (\vka-\vec{\bka'})}{\gamma}\right),
  \label{eq:T2}
\end{align}
for $z>0$, with initial condition $ \boldsymbol{\mathscr{P}}(\bka,0) =
\boldsymbol{\mathscr{A}}_o(\bka)
\boldsymbol{\mathscr{A}}_o(\bka)^\dagger$. In dimensional units this
is the same as \eqref{eq:MR13}.

We have from \eqref{eq:CENA} that
\begin{equation*}
  \int_{|\bka|<1} \frac{d (k \bka)}{(2\pi)^2} \Big[
    \EE\big[|a^\ep(\bka,z)|^2\big] +
    \EE\big[|a^{\ep,\perp}(\bka,z)|^2\big] \Big] = \int_{|\bka|<1} \frac{d (k
    \bka)}{(2\pi)^2} \Big[ |a_o(\bka)|^2 + |a^\perp_o(\bka)|^2 \Big],
\end{equation*} 
for all $z>0$ and $0 \le \ep \ll 1$, so energy conservation must hold
in the limit $\ep \to 0$. This is derived from \eqref{eq:T2} using the
commutation properties of the trace, as stated in \eqref{eq:MR12}.

By studying the second moments of the wave amplitudes, we can also
quantify the decorrelation of the waves due to scattering in the
random medium. The calculation of these moments uses the generator in
\ref{ap:DL} and the result is that the wave amplitudes at different
transverse wave vectors $\bka \ne \bka'$ are decorrelated
\begin{align*}
 \lim_{\ep \to 0} \EE \left[ \begin{pmatrix}
      a^\ep(\bka,z) \\ a^{\ep,\perp}(\bka,z) \end{pmatrix} 
 \begin{pmatrix}
   a^\ep(\bka',z) \\ a^{\ep,\perp}(\bka',z) \end{pmatrix}^\dagger \right] =
 \boldsymbol{\mathscr{A}}  (\bka,z)  \boldsymbol{\mathscr{A}}  (\bka',z)^\dagger
.
\end{align*}
Intuitively, this is because these waves travel in different
directions and see a different region of the random medium. It is only
when $|\bka - \bka'| = O(\ep)$ that the waves are correlated, and we
can calculate the energy density (Wigner transform)
\begin{align}
\bcW(\bka,\bx,z) =& \lim_{\ep \to 0} \int \frac{d(k {\itbf q})}{(2
  \pi)^2} \, \exp\Big[ i k {\itbf q} \cdot \big(\nabla \beta(\bka) z +
  \bx\big) \Big]
  \EE \left[ \begin{pmatrix} a^\ep\big(\bka
    + \frac{\ep {\itbf q}}{2},z\big) \\ a^{\ep,\perp}\big(\bka +
    \frac{\ep {\itbf q}}{2},z\big) \end{pmatrix}
 \begin{pmatrix}
   a^\ep\big(\bka-\frac{\ep {\itbf q}}{2},z\big)
   \\ a^{\ep,\perp}\big(\bka-\frac{\ep {\itbf
       q}}{2},z\big) \end{pmatrix}^\dagger \right].
\label{eq:T6}
\end{align}
It satisfies the transport equation 
\begin{align}
 \partial_z \bcW(\bka,\bx,z) - \nabla \beta(\bka) \cdot \nabla_\bx
  \bcW(\bka, \bx,z) = \mathbf{Q}(\bka) \bcW(\bka,\bx,z) +
  \bcW(\bka,\bx,z) \mathbf{Q}(\bka)^\dagger \nonumber \\
  +\frac{k^2 {\alpha}^2}{4 \gamma^3}\int_{|\bka'| < 1} \frac{d (k
    \bka')}{(2 \pi)^2}\, \bGamma(\bka,\bka') \bcW(\bka',\bx,z)
  \bGamma(\bka',\bka) \widetilde
  \cR\left(\frac{k(\vka-\vec{\bka'})}{\gamma} \right),
  \label{eq:T7}
\end{align}
for $z>0$. In dimensional units this is the same as \eqref{eq:MR9}.
When the initial condition is smooth, we have from \eqref{eq:T6} that
\begin{equation}
\label{eq:T8}
\bcW(\bka,\bx,0) = \delta(\bx)
\boldsymbol{\mathscr{A}}_o(\bka)\boldsymbol{\mathscr{A}}_o(\bka)^\dagger,
\end{equation}
and therefore at $z >0$
\begin{equation}
\label{eq:T8b}
\bcW(\bka,\bx,z) = \delta\left(\bx + \nabla \beta(\bka) z\right)
\boldsymbol{\mathscr{P}}(\bka ,z ),
\end{equation}
so energy is transported on the characteristic
\[
\bx = - \nabla \beta(\bka) z = \frac{\bka}{\beta(\bka)} z.
\]

\textbf{Remark 2:} It can also be shown, with an analysis similar to
that in \ref{ap:DL}, that the waves decorrelate over frequency offsets
larger than $\ep$. Thus, one can study the energy density resolved
over both time and space i.e., the space-time Wigner transform. We do
not consider this here, as we limit our study to a single frequency.

\section{Illustration of the scattering effects in statistically isotropic media}
\label{sect:DEPOL}
In this section we give numerical illustrations of the net scattering
effects in isotropic media. 
We already stated  that in transverse
isotropic media, where $ \cR$ depends only on $|\bx|$ and $z$, the
matrix $\mathbf{Q}(\bka)$ defined in \eqref{eq:C5} is diagonal, and
depends only on $|\bka|$.  We also defined in \eqref{eq:MR7} the
scattering mean free paths $\cS(\bka)$ and $\cS^\perp(\bka)$, which
depend only on $|\bka|$.  The illustrations in this section are for
statistically isotropic media, where $\cR (\vx)= \cR_{\rm iso}(|\vx|)$
and as shown in \ref{ap:RT}, the real part of matrix ${\bf Q}(\bka)$
is a multiple of the identity. Therefore $\cS(\bka) =
\cS^\perp(\bka)$, meaning that both components of the waves lose
coherence on the same $|\bka|$-dependent range scale.  

%we consider transverse isotropic media, where
%$\cR(\vx)$ depends only on $|\bx|$ and $z$, and isotropic media, where
%$\cR (\vx)$ depends only on $|\vx|$.  We give numerical illustrations of the net scattering
%effects in isotropic media. 
%Then we present the high-frequency
%analysis of the transport equations \eqref{eq:T2} in section
%\ref{sect:HF}.
%\subsection{Illustration of net scattering effects in isotropic media}
%\label{sect:num}

\begin{figure}[t]
 \begin{center}
    \includegraphics[width=7.3cm]{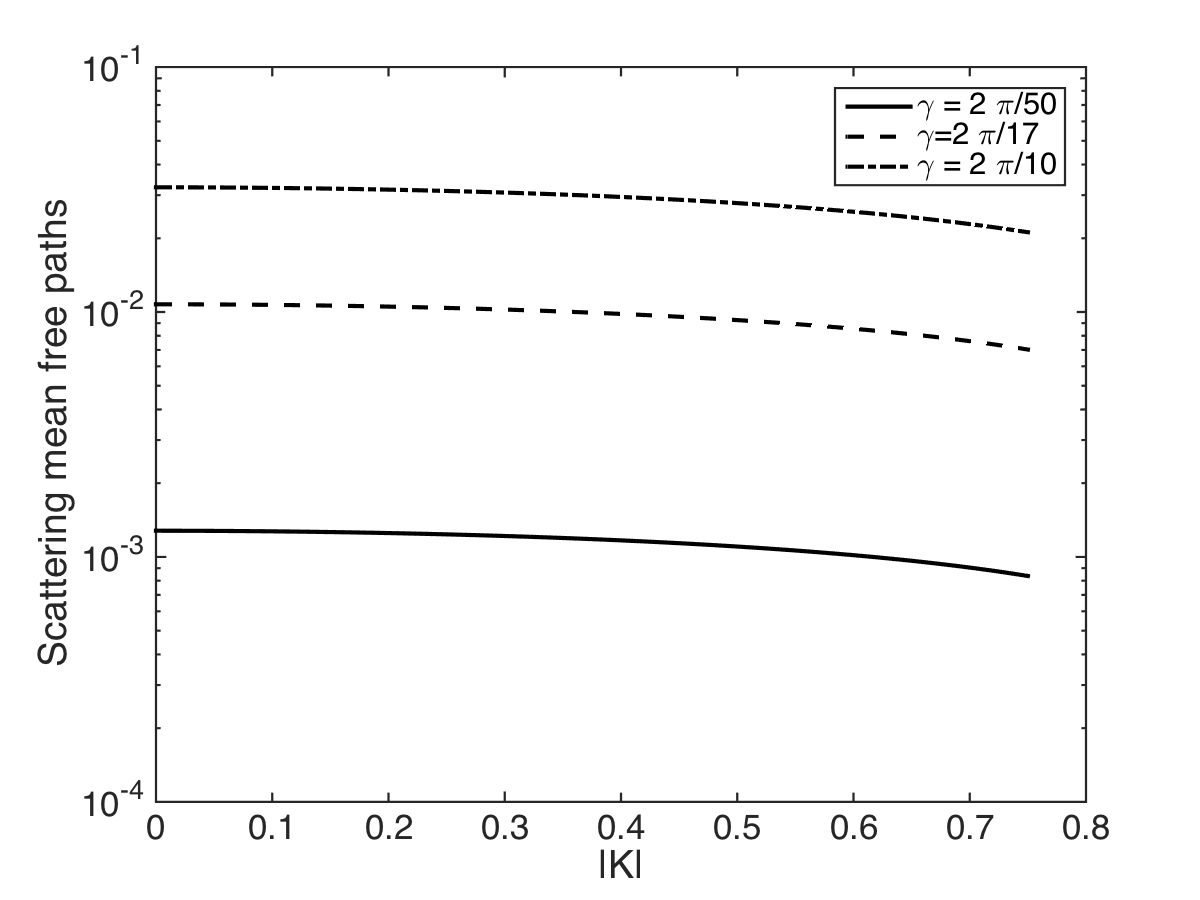}
    \end{center}
\caption{Illustration of the scattering mean free paths in an
  isotropic medium with Gaussian autocorrelation, for three different values of 
  $\gamma = \lambda/\ell$, as indicated in the legend. The abscissa is the magnitude
  $|\bka|$ of the transverse wave number. }
\label{fig:SCMF}
\end{figure}

In Figure \ref{fig:SCMF} we plot the scattering mean free paths for a
Gaussian autocorrelation $ \cR(\vx) = \exp (-{|\vx|^2}/{2} ),
$ and for three different values of $\gamma = \la/\ell$ equal to $2
\pi/50$, $2\pi/17$ and $2 \pi/10$.  We observe that the scattering
mean free paths decrease with $\gamma$, meaning that higher 
frequency waves lose coherence faster. The scattering mean free paths
also decrease monotonically with $|\bka|$. We may
understand this intuitively by noting that at a given range $z$, the
length of the path of the waves with transverse vector $|\bka|$ is
$z/\beta(\bka)$. This increases with $|\bka|$, and the waves lose
coherence faster because they travel a longer distance in the random
medium. The monotonicity of the scattering mean free paths can also be
seen clearly in Section \ref{sect:HF1}, where we study the
high-frequency limit $\gamma \to 0$.

To illustrate the evolution of the coherence matrix $
\boldsymbol{\mathscr{P}}(\bka,z)$, satisfying the evolution equation
\eqref{eq:T2}, we consider first the case of isotropic initial
conditions, where $\boldsymbol{\mathscr{A}}_o(\bka)$ depends only on
$|\bka|$, and we excite either the transverse magnetic 
or electric waves. From \eqref{eq:WD8bp} we see that this corresponds 
to a current source with $\what J_z = 0$ and $\what \bJ$ along $\bka$ or 
$\bka^\perp$. In the spatial coordinates, the latter means that $\bJ$ is the gradient 
or the perpendicular gradient of a scalar valued function. When the initial conditions are isotropic, 
the diagonal and off-diagonal terms of the coherence matrix
$ \boldsymbol{\mathscr{P}}(\bka,z)$ decouple. Indeed,
there is no coupling in the first two terms of equation \eqref{eq:T2},
because ${\bf Q}(\bka)$ is diagonal.  The integral in \eqref{eq:T2}
can be analyzed with a Picard iteration. For each iterate we obtain by
explicit calculation, using integration in polar coordinates, that the
coupling terms of the kernel vanish, and that the result preserves the
isotropic dependence on $\bka$. We obtain two autonomous systems of
equations for the power vector
\begin{equation}
\begin{pmatrix} \mathscr{P}_{11}(\bka,z) \\ \mathscr{P}_{22}(\bka,z) 
\end{pmatrix}= \lim_{\ep \to 0} \begin{pmatrix} \EE \big[
    |a^\ep(\bka,z)|^2 \big] \\ \EE\big[|a^{\ep,\perp}(\bka,z)|^2
    \big] \end{pmatrix} ,
\label{eq:EVect}
\end{equation}
and the off-diagonal terms
\begin{equation}
  \mathscr{P}_{12}(\bka,z) = \lim_{\ep \to 0} \EE\big[ a^\ep(\bka,z)
    \overline{a^{\ep,\perp}(\bka,z)}\big].
\label{eq:EVect12}
\end{equation}
This means in particular that if the wave is transverse electric or
magnetic initially, so that ${\mathscr{P}}_{12}(\bka,0) = 0$, the
coherence matrix $ \boldsymbol{\mathscr{P}}(\bka,z)$ remains diagonal
for all $z$.

\begin{figure}[!h]
 \begin{center}
    \includegraphics[width=0.49\textwidth]{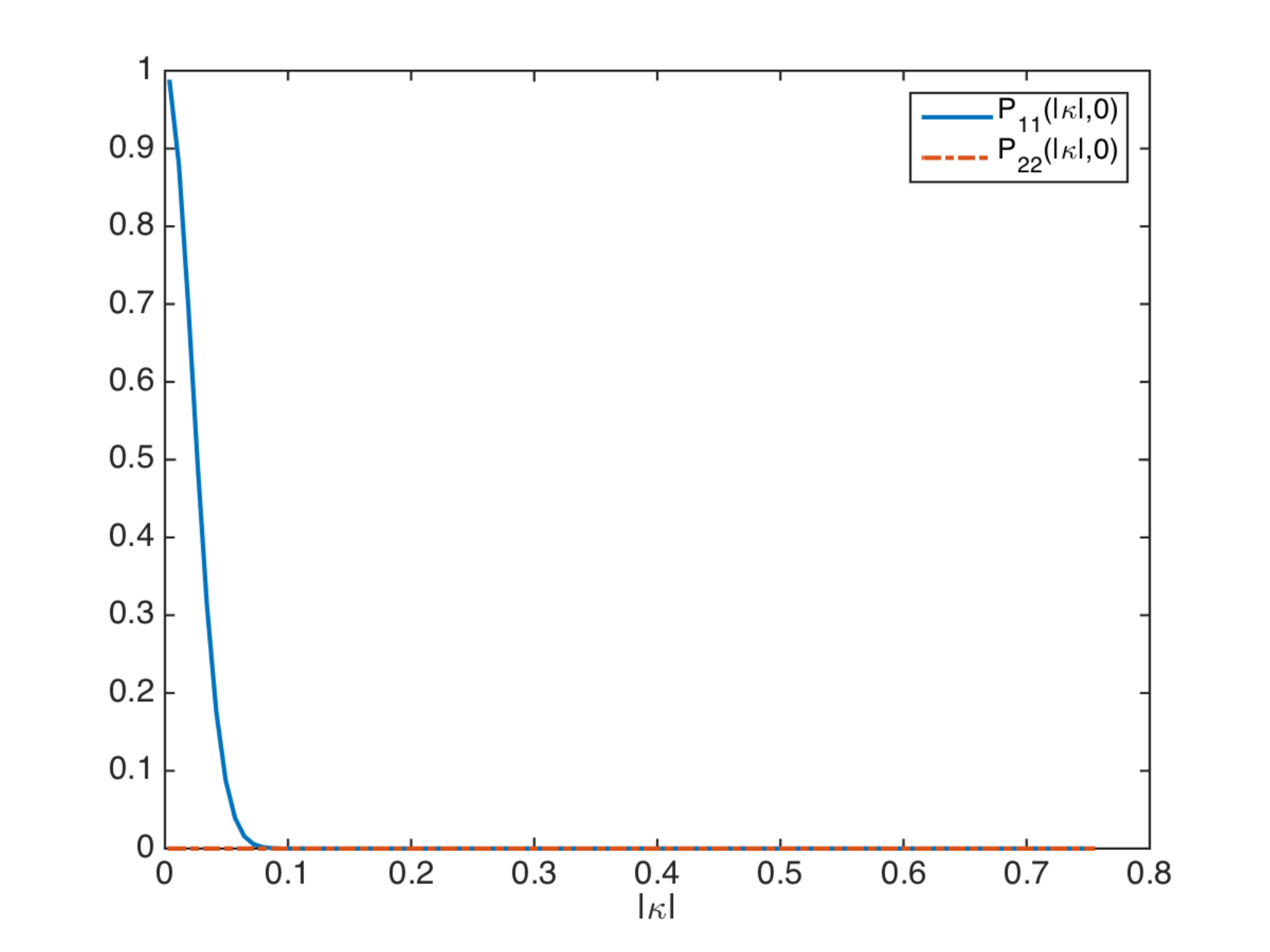}
    \includegraphics[width=0.49\textwidth]{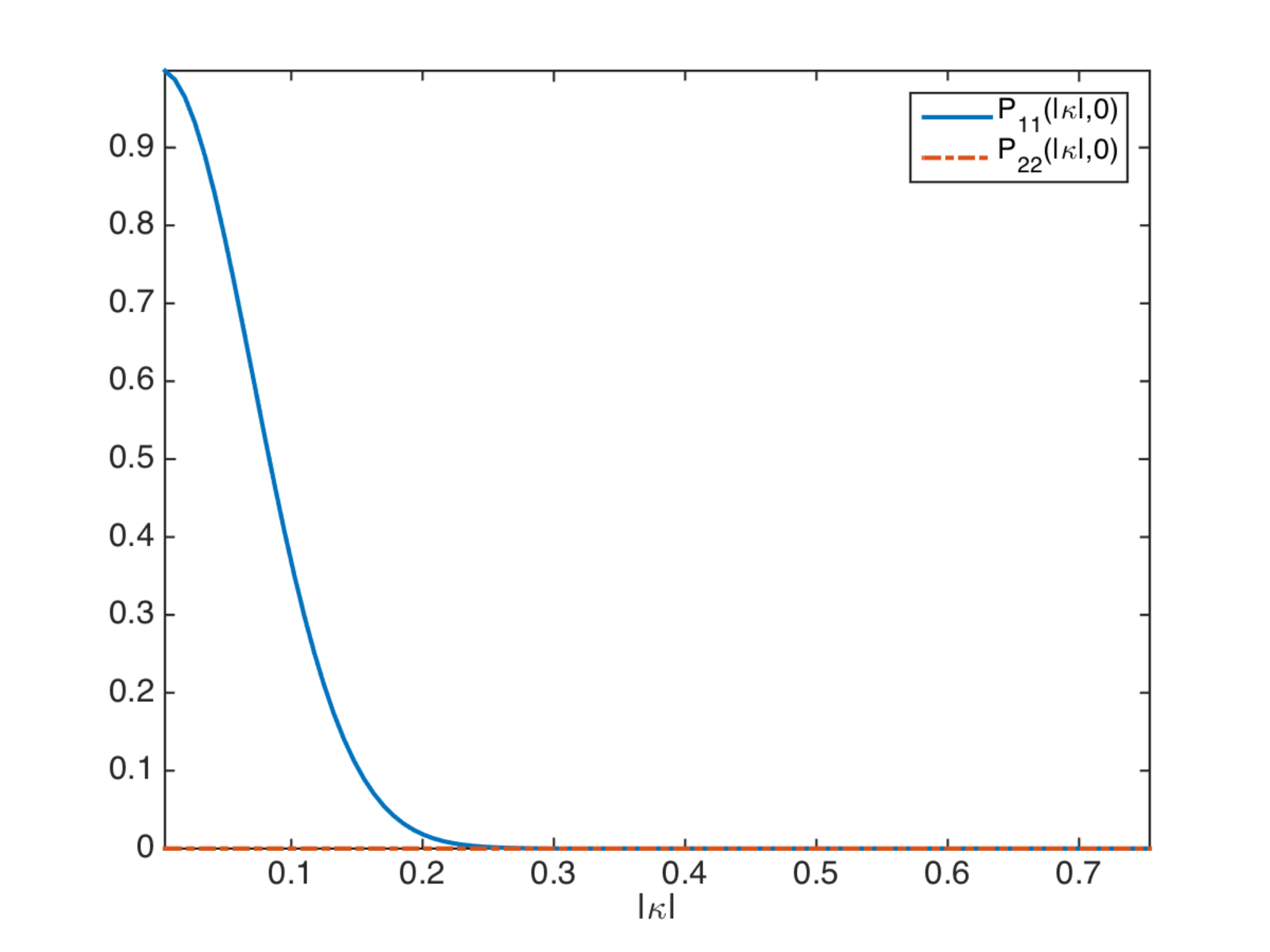} \\
     \includegraphics[width=0.49\textwidth]{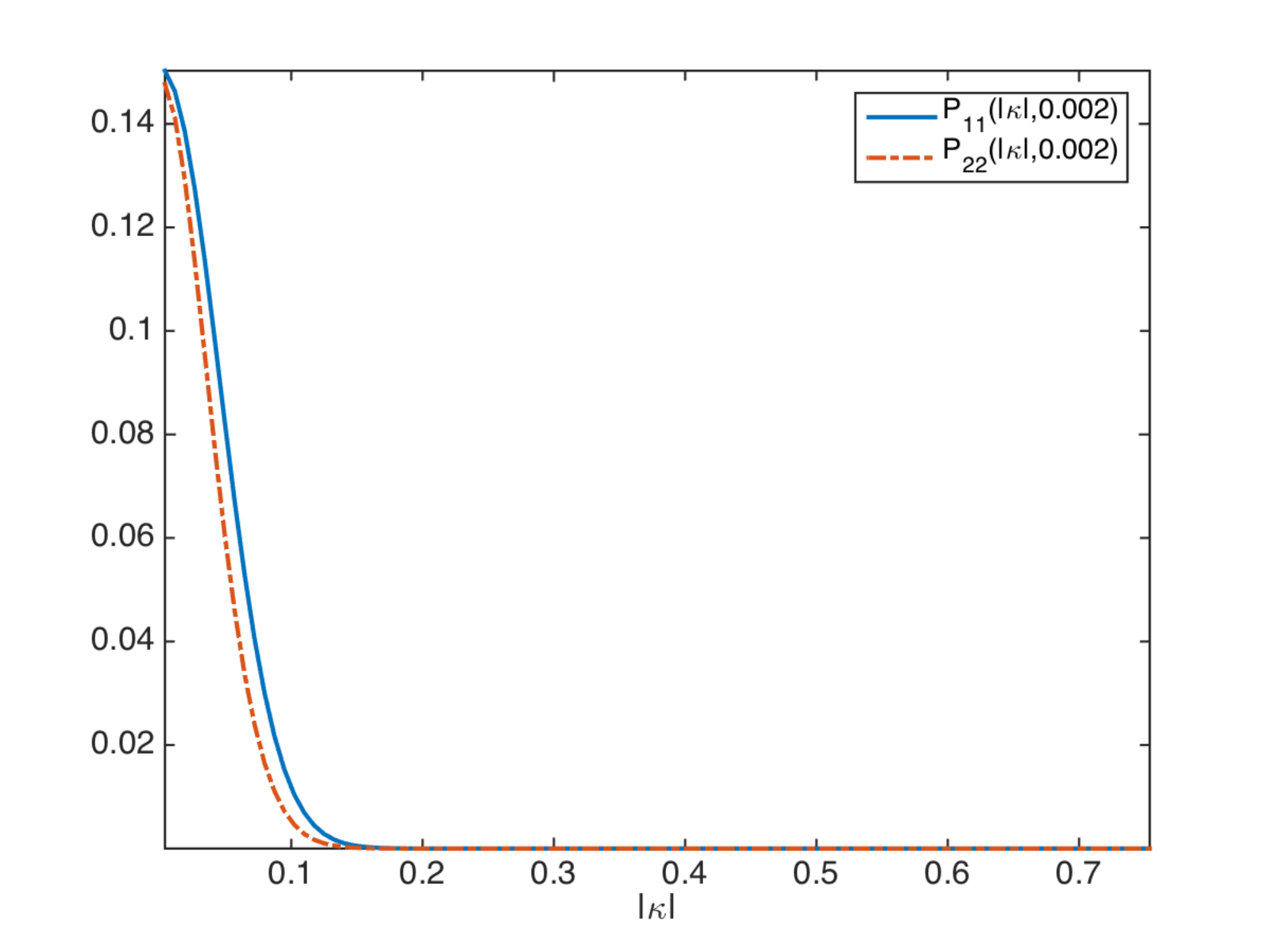}
    \includegraphics[width=0.49\textwidth]{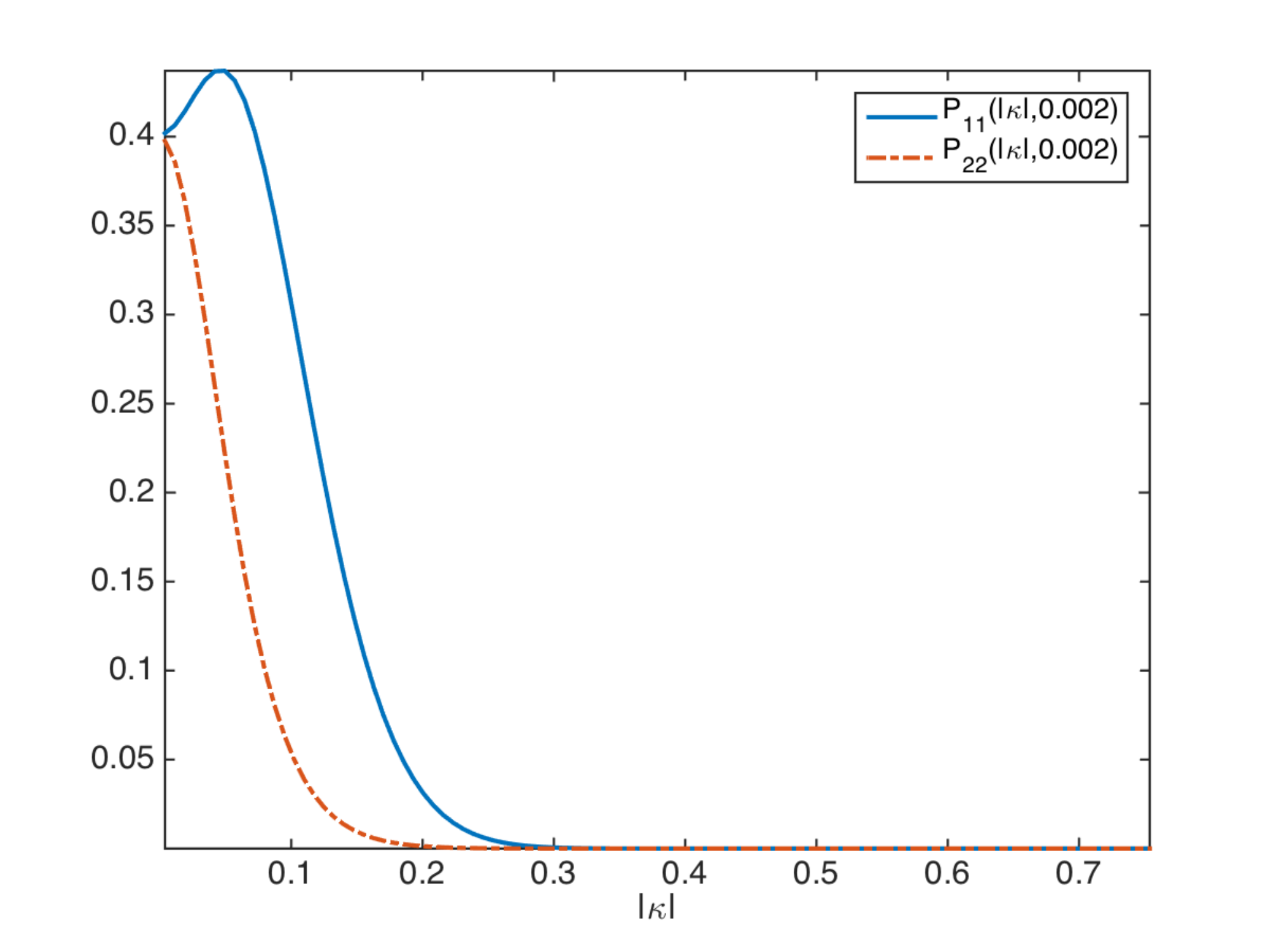} \\ 
\includegraphics[width=0.49\textwidth]{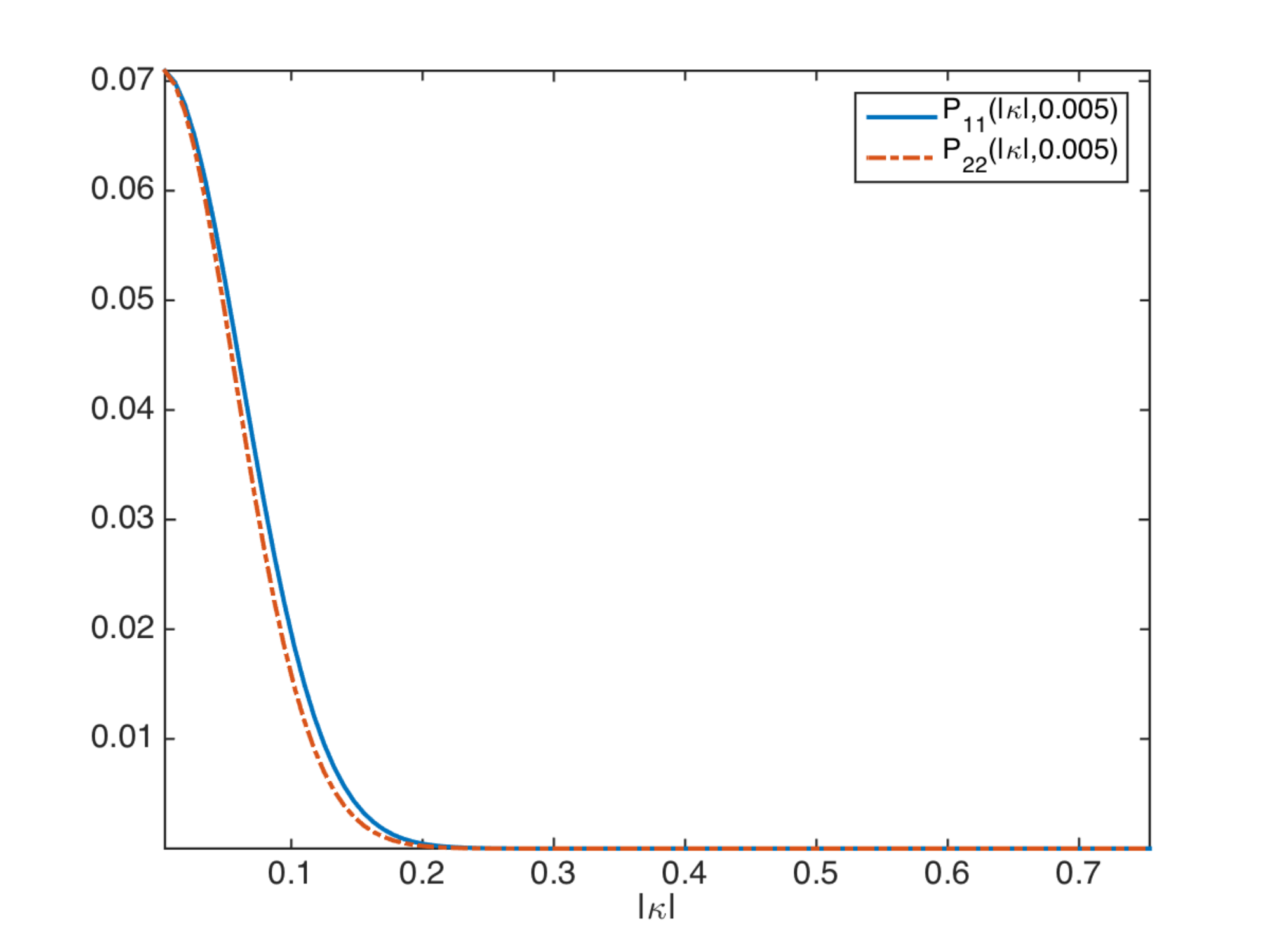}
    \includegraphics[width=0.49\textwidth]{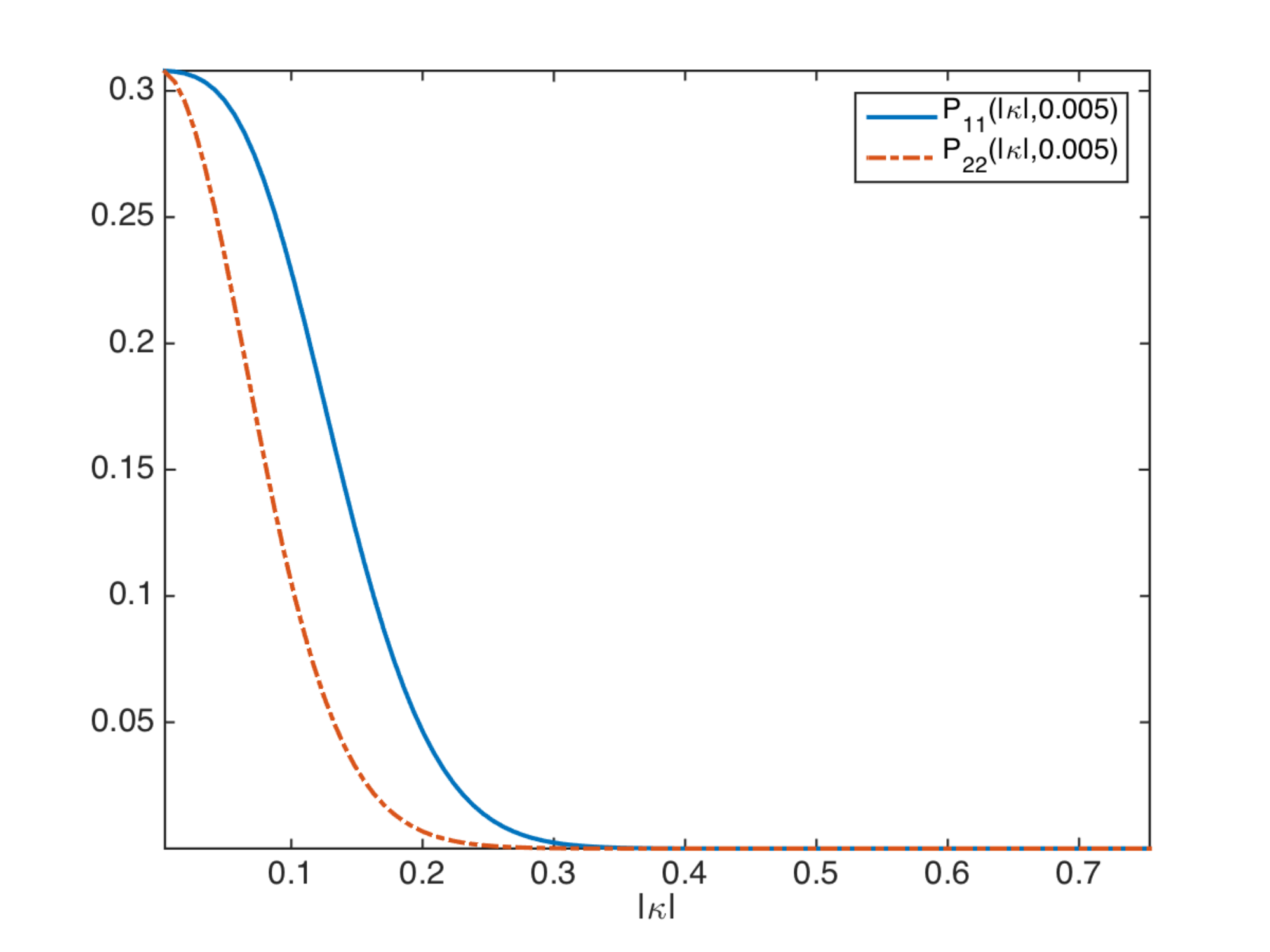}
    \end{center}
\caption{Display of ${\mathscr{P}}_{11}(\bka,z)$ (full blue line) and
  $ {\mathscr{P}}_{22}(\bka,z)$ (dashed red line) as a function of
  $|\bka|$, for three different values of $z$: Top line is for $z =
  0$, middle line is for $z = 2 \times 10^{-3}$ and the bottom line
  for $z = 5 \times 10^{-3}$. We take $\gamma = 2\pi/50$,
  $\gamma_{_\mJ} =\gamma$ in the left column and $\gamma_{_\mJ} = 3.5
  \gamma$ in the right column.  }
\label{fig:P_ISO_EVOL}
\end{figure}

We illustrate the transfer of power between the modes due to
scattering by displaying in Figure \ref{fig:P_ISO_EVOL} the
$z$-evolution of the power vector \eqref{eq:EVect} for the initial
condition
\[
\begin{pmatrix} \mathscr{P}_{11}(\bka,z) \\ \mathscr{P}_{22}(\bka,z) 
\end{pmatrix}= \begin{pmatrix} |a_o(\bka)|^2 \\ |a_o^\perp(\bka)|^2 \end{pmatrix}, 
\qquad |a_o(\bka)|^2 = \exp \left(-\frac{|\bka|^2}{2 \gamma_{_\mJ} ^2}\right),
\quad a_o^\perp(\bka) = 0.
\]
The electromagnetic plane wave at wave vector $\vka$ is transverse magnetic initially
and we have in terms of the Stokes vector defined in \eqref{eq:Stokes} 
\[
\boldsymbol{\mathfrak{S}}(\bka,0) = \exp \left(-\frac{|\bka|^2}{2
  \gamma_{_\mJ} ^2}\right) \left(1,1,0,0\right).
\]
In Figure \ref{fig:P_ISO_EVOL} we show the
results for $\gamma = 2 \pi/10$, $\gamma_{_\mJ} = \gamma$ and
$\gamma_{_\mJ} = 3.5 \gamma$. We plot ${\mathscr{P}}_{11}$ and
${\mathscr{P}}_{22}$ as functions of $|\bka|$, for $z = 0$, $z = 2
\times 10^{-3}$ and $z = 5 \times 10^{-3}$. Note from Figure
\ref{fig:SCMF} that the scattering mean free paths at $\bka = 0$ are
approximately $1.3 \times 10^{-3}$. Thus, in the last line of the
plots in Figure \ref{fig:P_ISO_EVOL} the waves have traveled
approximately five scattering mean free paths. We observe from the
figure that energy from the transverse magnetic waves, which are the
only ones excited initially, is transmitted to the transverse electric
modes.  This energy transfer is faster at smaller $|\bka|$, so the
modes have equal power at five scattering mean free paths in the case
$\gamma_{_\mJ} = \gamma$, whereas for $\gamma_{_\mJ} = 3.5 \gamma$,
the transverse magnetic modes still carry more energy. 
Recalling the
Stokes vector \eqref{eq:Stokes}, we note that when the transverse
electric and magnetic waves have equal power, we have
$\boldsymbol{\mathfrak{S}}(\bka,z) = 2 \mathscr{P}_{11}(\bka,z)
\left(1,0,0,0\right)$
since $\mathscr{P}_{12} = 0$ in the simulation,
and the electromagnetic plane wave at wave vector $\vka$ is
unpolarized.  
Aside from the transfer of energy between
the transverse electric and magnetic components of the waves, we note
in Figure \ref{fig:P_ISO_EVOL} the diffusion of energy at higher
$|\bka|$. For example, in the case of $\gamma_{_\mJ} = \gamma$, the
support of the energy is at $|\bka| < 0.1$ initially (top left plot)
but it extends to $|\bka|<0.2$ at five scattering mean free paths
(bottom left plot).

\begin{figure}[t]
 \begin{center}
    \includegraphics[width=7.3cm]{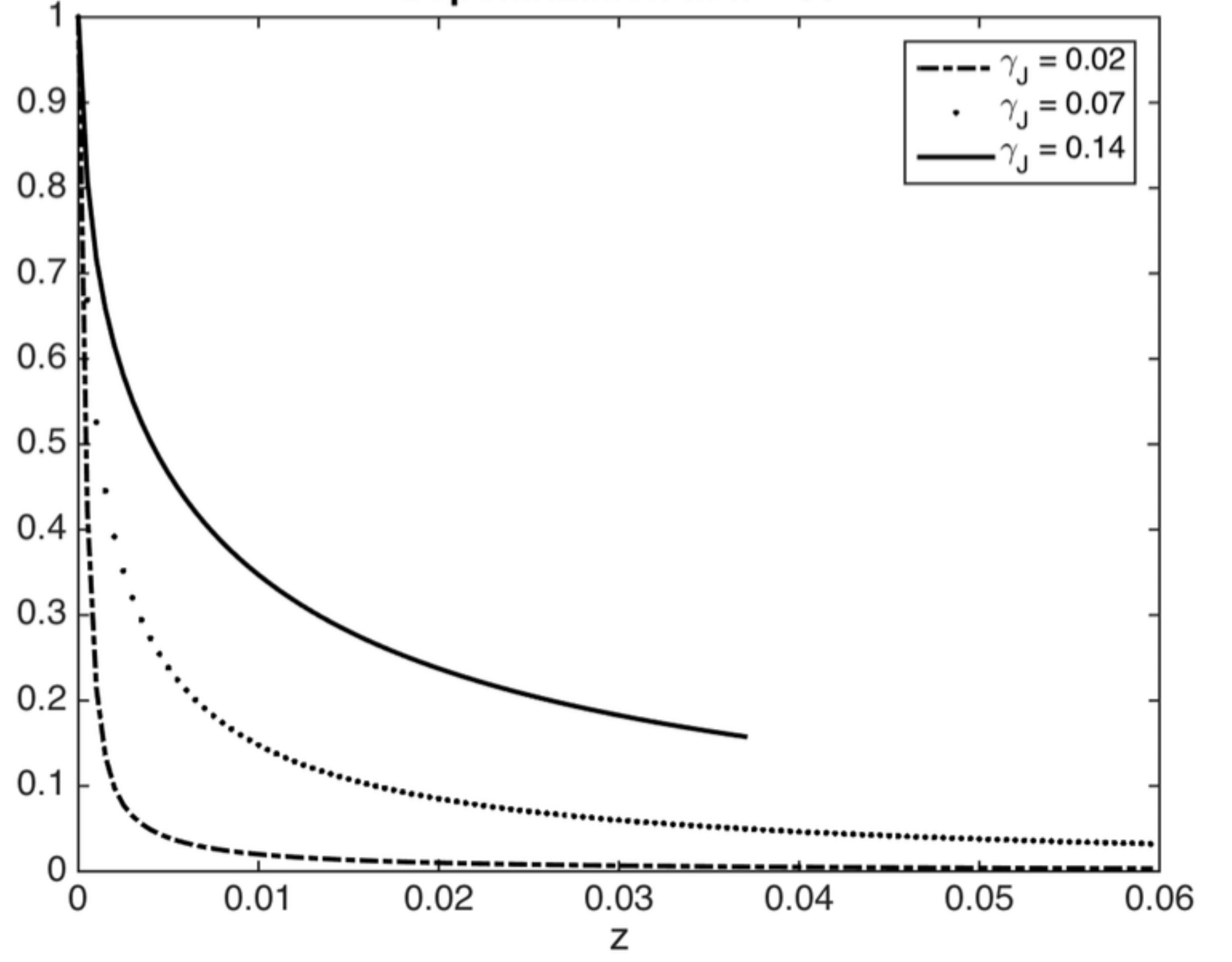}
    \end{center}
\caption{Plot of the coefficient of power distribution
  ${C}_{_\mathcal{P}}(z)$. At $z = 0$ all the power is in the
  transverse magnetic mode. As $z$ grows the other mode gains power,
  and eventually the energy becomes equally distributed between the modes. }
\label{fig:DEPOL}
\end{figure}

Figure \ref{fig:DEPOL} is another illustration of the effect of the
initial conditions on the power exchange between the modes. We
consider there, in addition to $\gamma_{_\mJ} = \gamma$ and $3.5
\gamma$, the case of an even larger $\gamma_{_\mJ} $, equal to $7
\gamma$. We plot in Figure \ref{fig:DEPOL} the coefficient of power
distribution ${C}_{_\mathcal{P}}(z)$, calculated as the difference of
the mode powers, normalized by the total power, which is constant in
$z$,
\begin{align}
\label{def:degp}
C_{_\mathcal{P}}(z) = \frac{\int_{|\bka| < 1} d \bka \, \big[
    {\mathscr{P}}_{11}(\bka,z) - {\mathscr{P}}_{22}(\bka,z)\big]}{
  \int_{|\bka| < 1} d \bka\,  \big[ {\mathscr{P}}_{11}(\bka,z) +
    {\mathscr{P}}_{22}(\bka,z)\big]}.
\end{align}
Consistent with the left bottom plot in Figure \ref{fig:P_ISO_EVOL},
$C_{_\mathcal{P}}$ is almost zero at $z = 5 \times 10^{-3}$ for
$\gamma_{_\mJ} = \gamma$. However, when $\gamma_{_\mJ} = 7 \gamma$
i.e., when the diameter $X$ of the source is $7$ times smaller and the
opening angle of the emitted wave cone is $7$ times larger,
$C_{_\mathcal{P}}$ is still large, so the waves still remember the
initial excitation. We terminate the plot for $\gamma_{_\mJ} = 7
\gamma$ at $ z < 0.04$, where the opening angle of the cone comes
close to $180$ degrees and the evanescent waves start to get excited.

\begin{figure}[!h]
 \begin{center}
    \hspace{-0.12in}\includegraphics[width=0.34\textwidth]{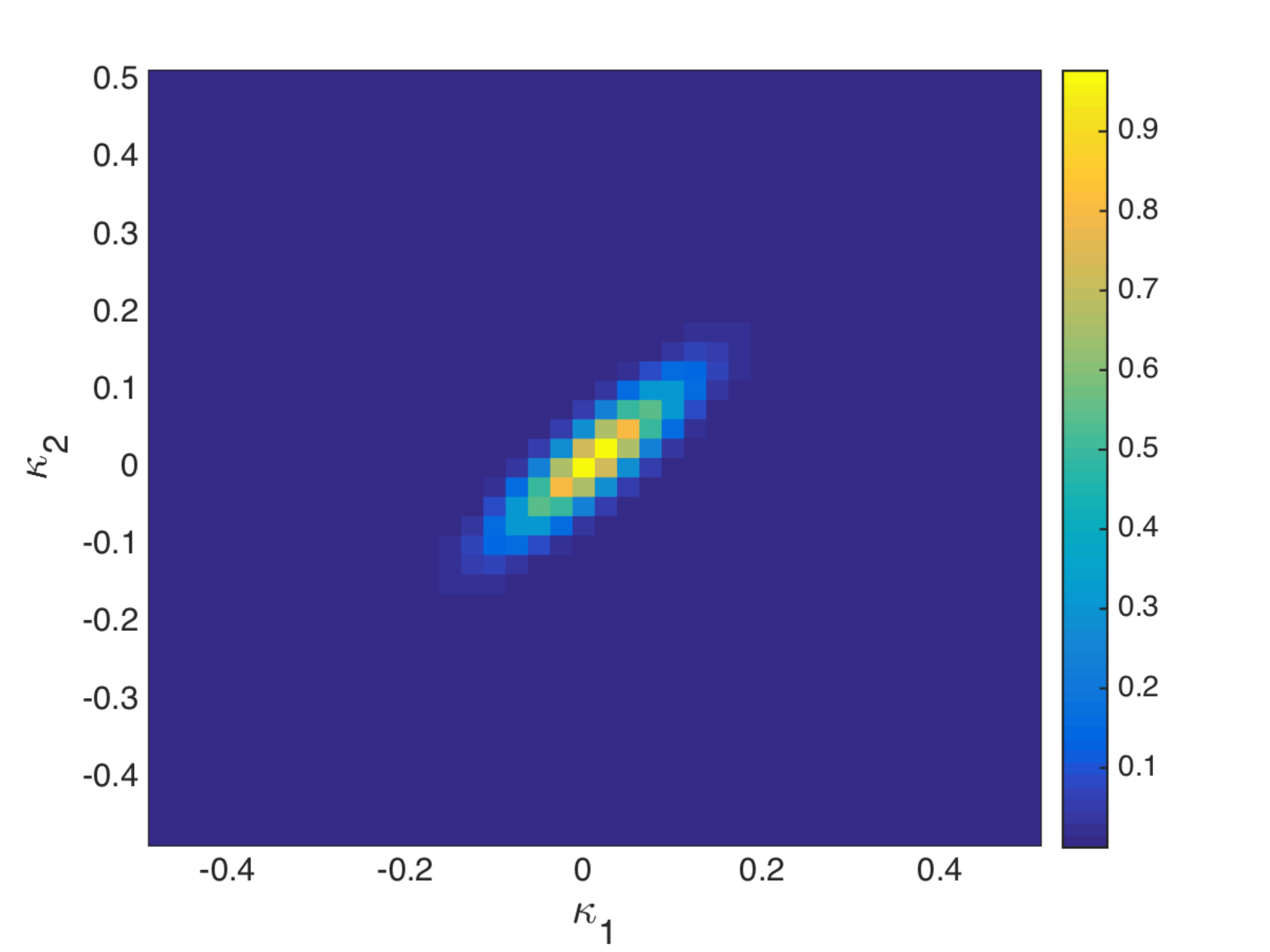}
    \hspace{-0.15in}\includegraphics[width=0.34\textwidth]{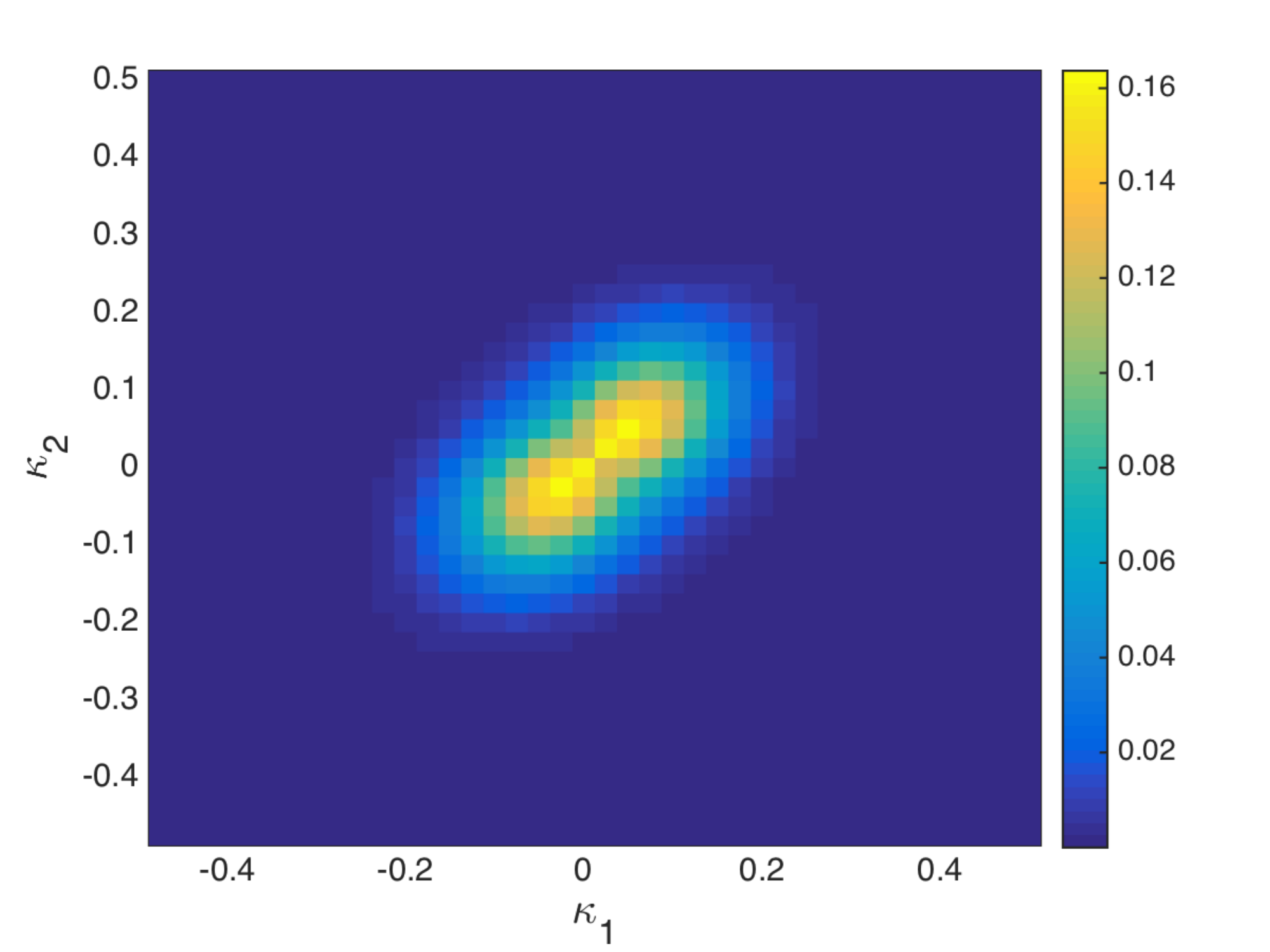}
    \hspace{-0.15in}\includegraphics[width=0.34\textwidth]{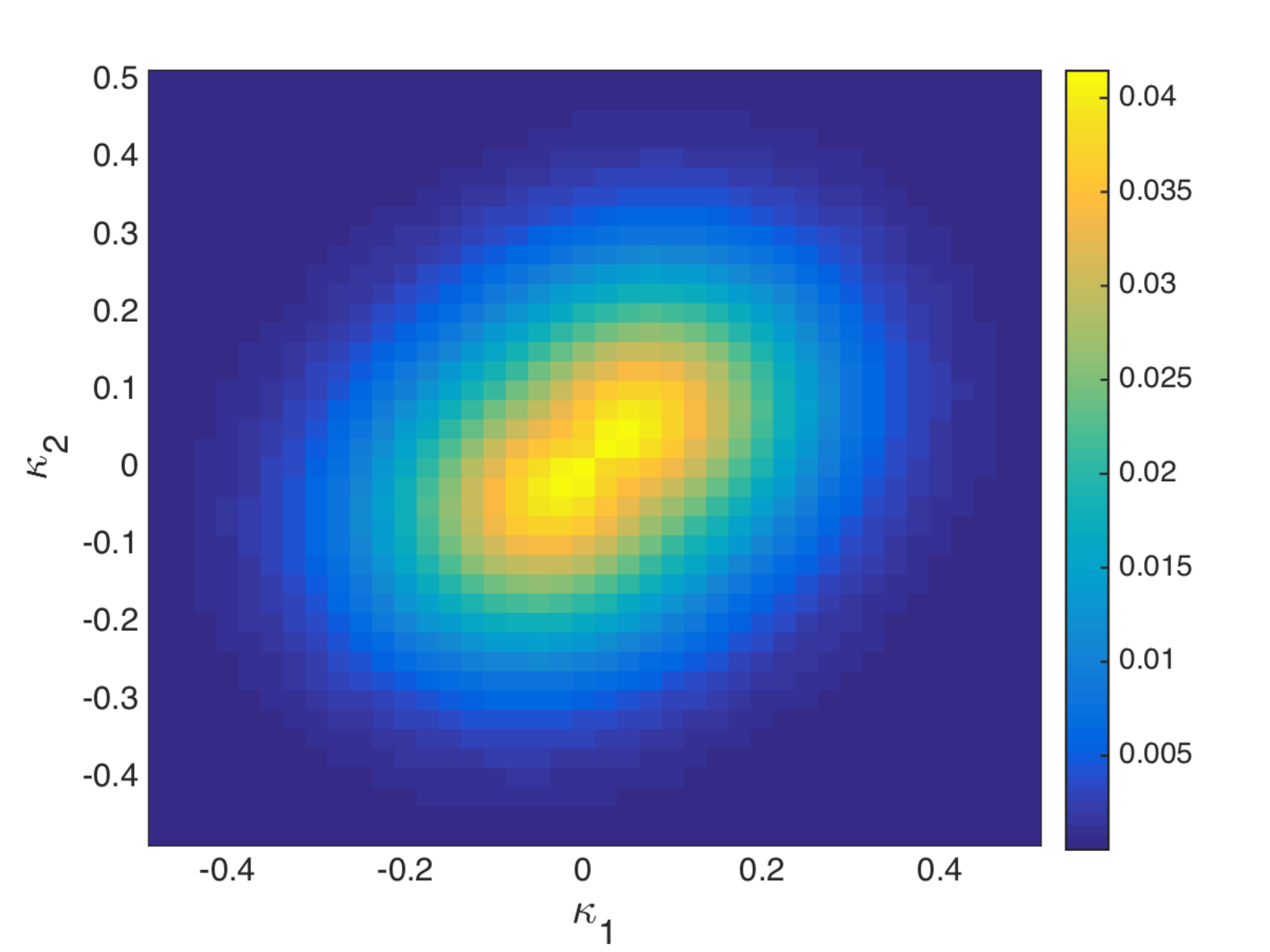}\\
    \hspace{-0.12in}\includegraphics[width=0.34\textwidth]{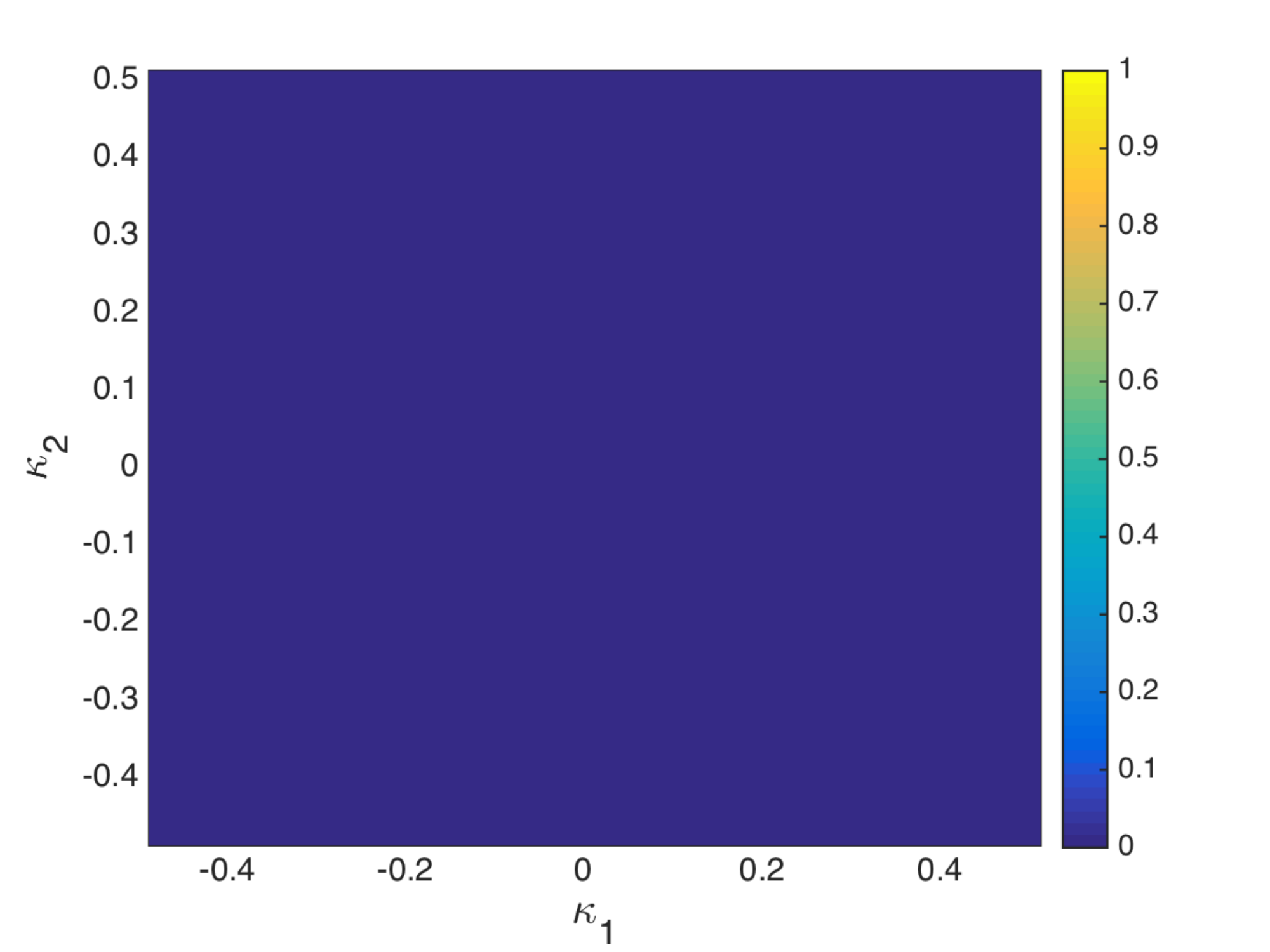}
    \hspace{-0.15in}\includegraphics[width=0.34\textwidth]{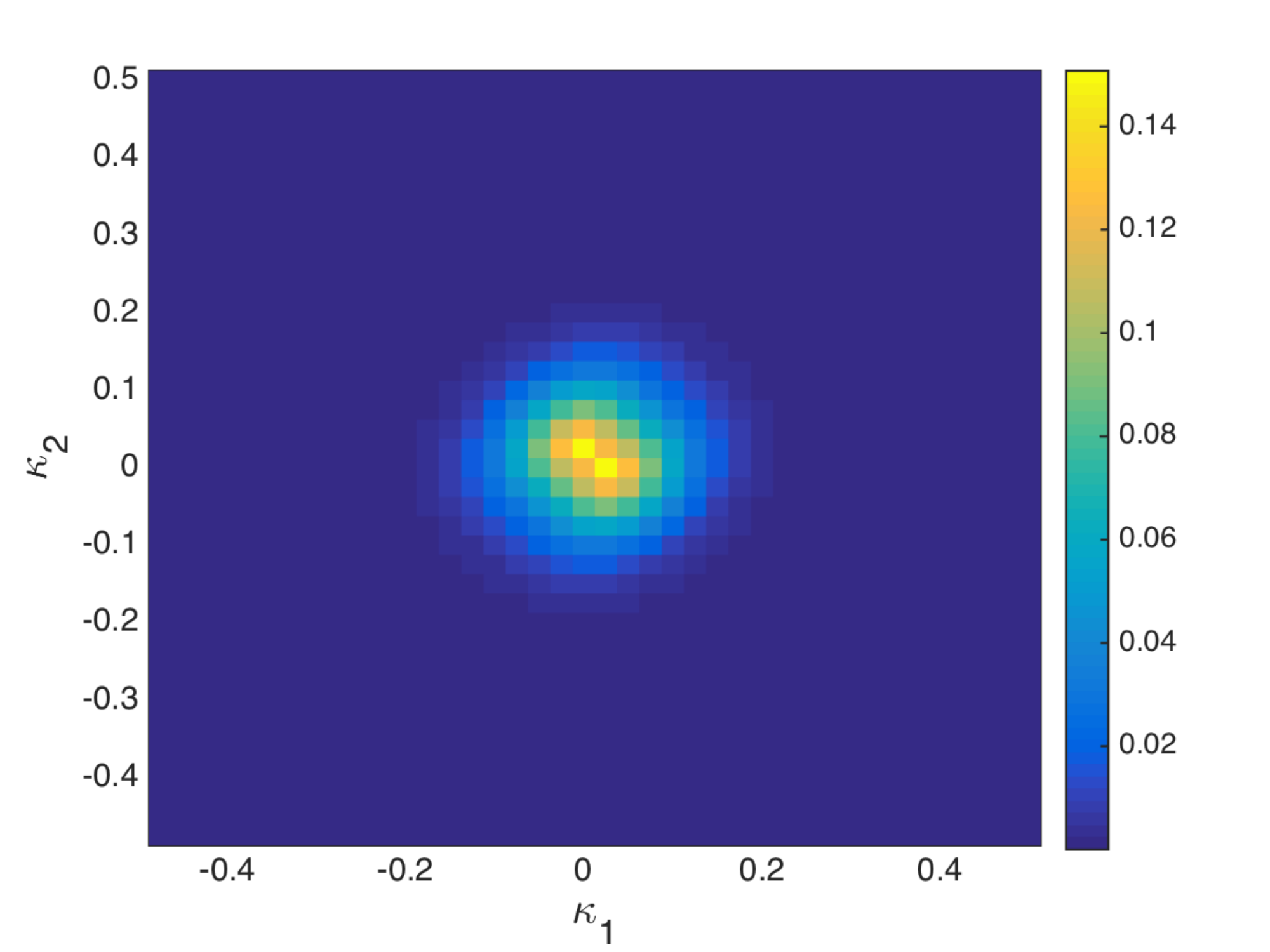}
    \hspace{-0.15in}\includegraphics[width=0.34\textwidth]{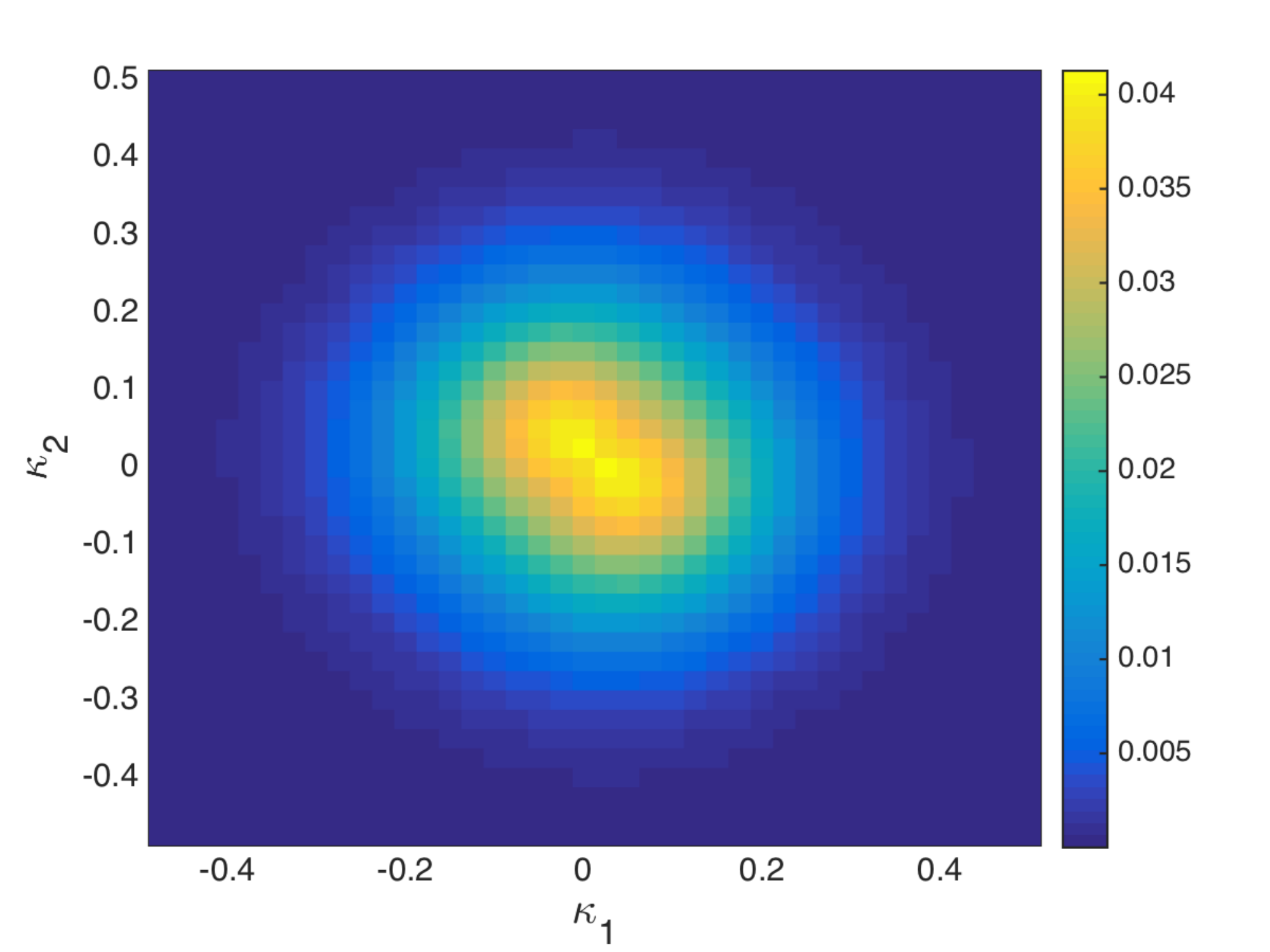}\\ 
    \hspace{-0.12in}\includegraphics[width=0.34\textwidth]{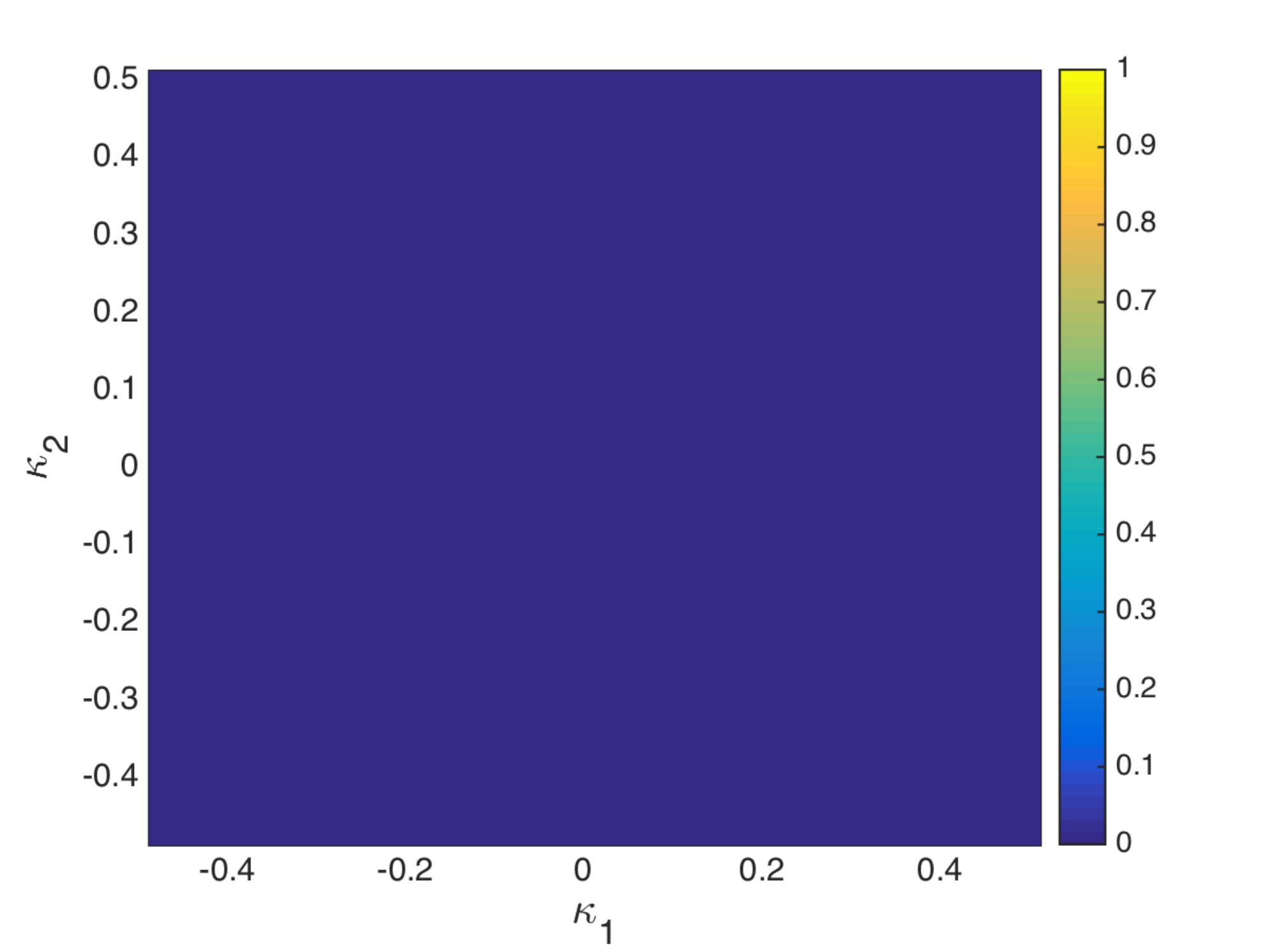}
   \hspace{-0.15in} \includegraphics[width=0.34\textwidth]{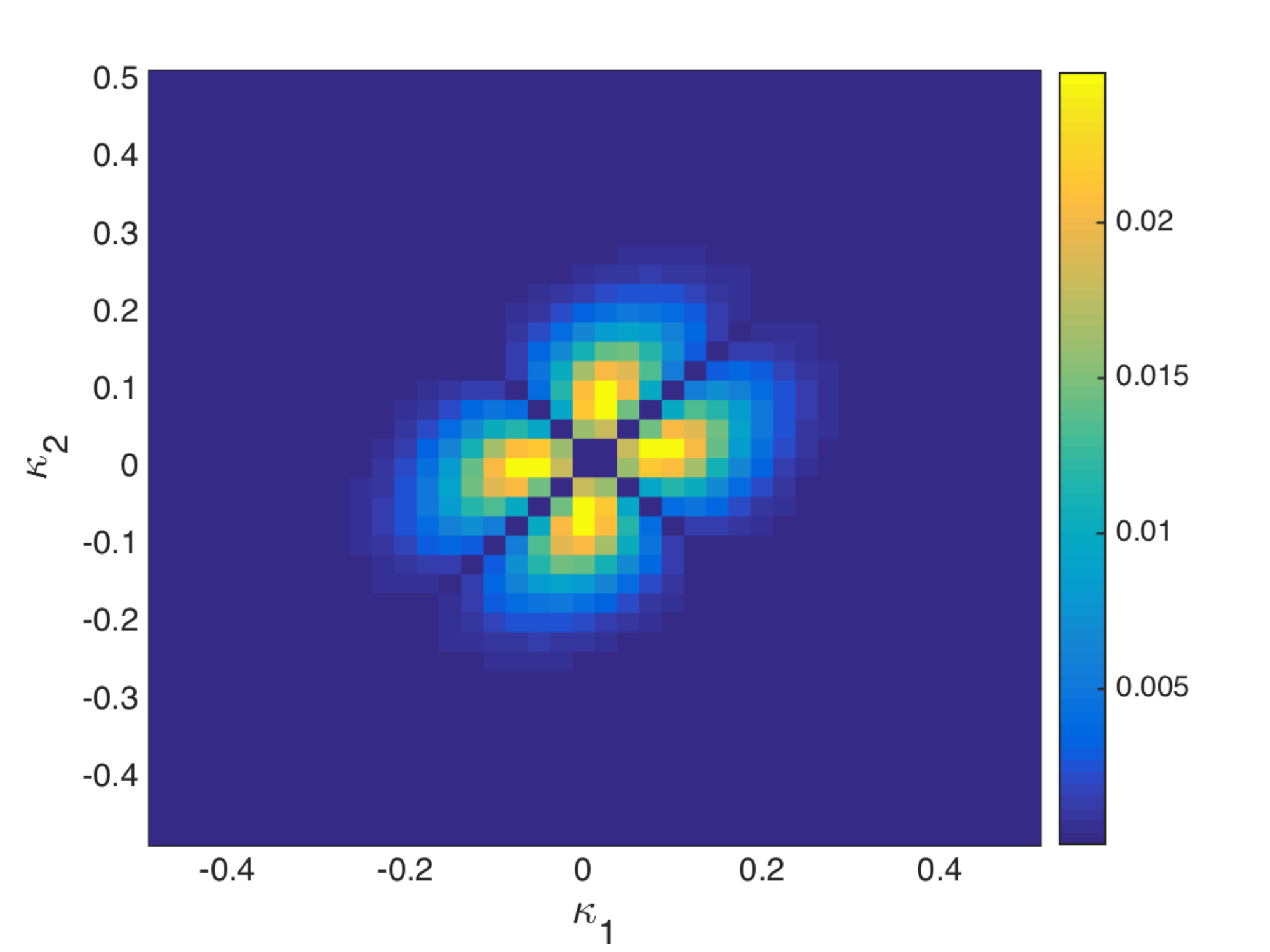}
    \hspace{-0.15in}\includegraphics[width=0.34\textwidth]{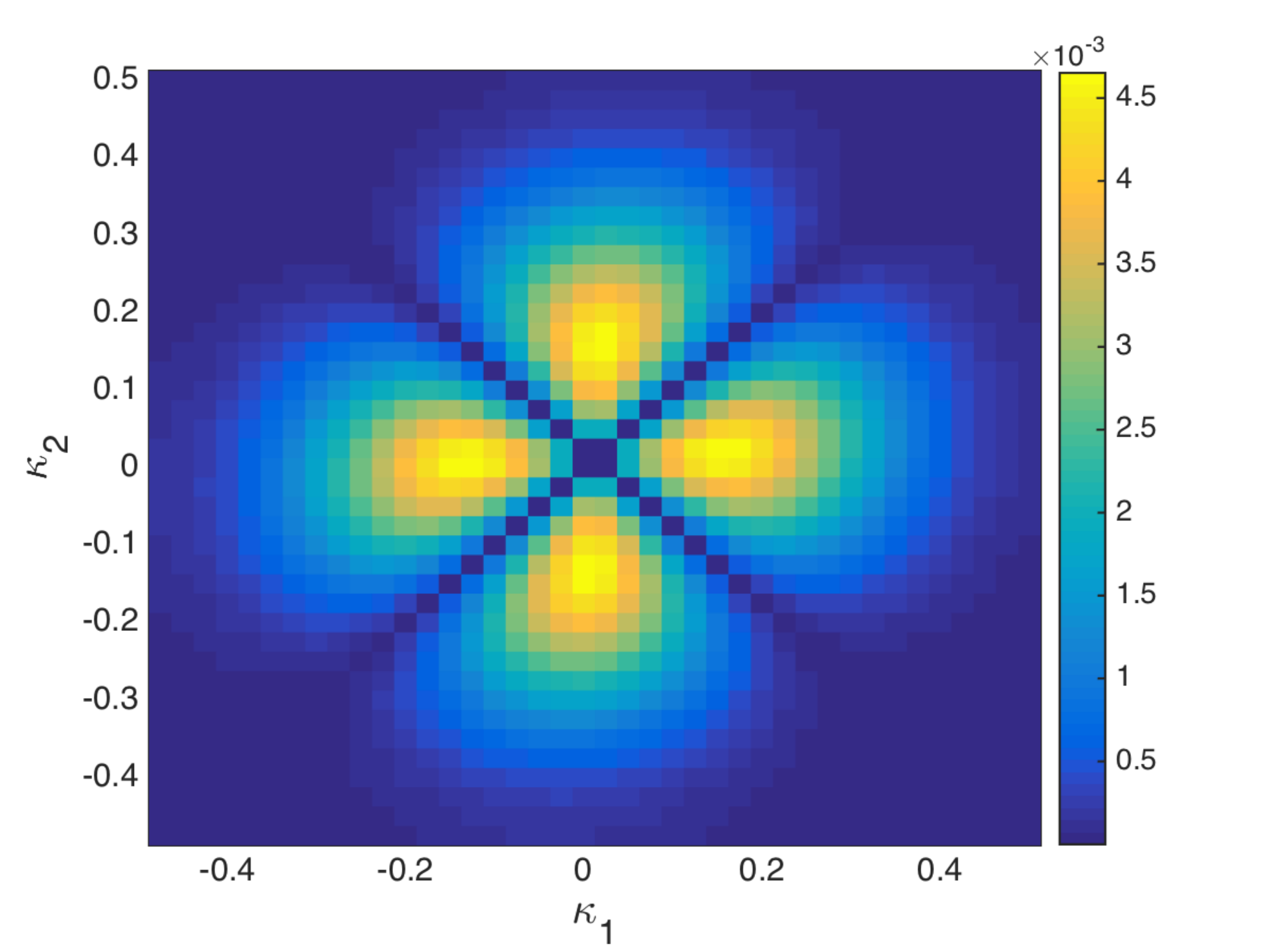}
    \end{center}
\caption{Display of $ {\mathscr{P}}_{11}(\bka,z)$ (top line),
  ${\mathscr{P}}_{22}(\bka,z)$ (middle line) and $|
  {\mathscr{P}}_{12}(\bka,z)|$ (bottom line) as a function of $\bka =
  (\kappa_1,\kappa_2)$, for $\gamma = 2 \pi/50$ and three different
  values of $z$: Left column $z = 0$, middle column $z = 5 \times
  10^{-3}$ and right column $z = 2.5 \times 10^{-2}$. }
\label{fig:P_ANISO}
\end{figure}

The results in Figure \ref{fig:P_ANISO} display the $z$-evolution of
the entries $ {\mathscr{P}}_{jl}(\bka,z)$ of the coherence matrix, for the
anisotropic initial condition 
\[
|a_o(\bka)|^2 = \exp \left[ -\frac{1}{2}
  \left(\frac{\kappa_1-\kappa_2}{0.03}\right)^2 -\frac{1}{2}
  \left(\frac{\kappa_1 + \kappa_2}{0.1}\right)^2 \right], \quad
a_o^\perp(\bka) = 0.
\]
Again, we observe the transfer of energy from the transverse magnetic
components to the transverse electric ones, which  is
more efficient at $|\bka| \approx 0$. We also note the diffusion of
energy to waves with larger wave vectors, which means physically that
the focused (beam-like) field emitted by the source spreads out as it
propagates along $z$. The off-diagonal components of the coherence
matrix are no longer zero, due to the anisotropic initial
condition. However, the solution approaches an isotropic one as $z$
grows, and the off-diagonal components decay.

\section{High-frequency analysis}
\label{sect:HF}%
In this section we analyze in detail the net scattering effects of the
random medium by considering the high-frequency limit $\gamma \to 0$.
%in which the transport equations simplify. 
%We assume a transverse isotropic medium, meaning that $\cR(\vec\bx) $ depends only on $|\bx|$ and $z$. 
We begin with the quantification of the scattering mean free paths and then analyze the
coherence matrix. 

\subsection{Scattering mean free paths}
\label{sect:HF1}
In transverse isotropic random media the matrix ${\bf Q}(\bka)$ is
diagonal. We are interested in its real part which defines the
scattering mean free paths \eqref{eq:MR7}. We have
\begin{align}
{\rm Re} \big[ {\rm Q}_{11}(\bka)\big] = - \frac{k^2 {\alpha}^2}{8 \gamma^3}
\int_{|\bka'|<1} \frac{d (k \bka')}{(2 \pi)^2} \widetilde \cR
\left(\frac{k(\vka-\vec{\bka'})}{\gamma}\right)
   \Big[  \Gamma_{11}^2(\bka',\bka) + \Gamma_{21}^2(\bka',\bka) \Big],
\label{eq:HF1}
\end{align}
and 
\begin{align}
{\rm Re} \big[ {\rm Q}_{22}(\bka)\big] = - \frac{k^2 {\alpha}^2}{8 \gamma^3} \int_{|\bka'|<1} \frac{d
  (k \bka')}{(2 \pi)^2} \widetilde \cR
\left(\frac{k(\vka-\vec{\bka'})}{\gamma}\right)
  \Big[ \Gamma_{22}^2(\bka',\bka) + \Gamma_{12}^2(\bka',\bka) \Big],
\label{eq:HF2}
\end{align}
where we recall that 
$
\vka - \vec{\bka'} = \left( \bka-\bka',
\beta(\bka)-\beta(\bka')\right)$.
Let us change variables $\bka' = \bka - \gamma \bu$, and expand
in powers of $\gamma$
\begin{align} 
\beta(\bka) - \beta(\bka-\gamma \bu) =  \gamma \bu \cdot \nabla \beta(\bka)
+ O(\gamma^2) = -\gamma \frac{\bu \cdot \bka}{\beta(\bka)} +
O(\gamma^2).
\end{align}
For the diagonal entries of $\bGamma$ we have
\begin{align}
\label{eq:EXPGD}
\Gamma_{jj}(\bka-\gamma \bu,\bka) = \frac{1}{\beta(\bka)} - \gamma
\frac{\bu \cdot \bka}{2 \beta^3(\bka)} + O(\gamma^2), \quad j = 1,2 ,
\end{align}
and for the off-diagonal ones
\begin{align}
\Gamma_{21}(\bka-\gamma \bu,\bka) = \gamma \frac{\bu \cdot
  \bka^\perp}{|\bka| |\bka - \gamma \bu|} + O(\gamma^2),
\end{align}
and similar for $\Gamma_{12}$, with a negative sign. Substituting into
(\ref{eq:HF1}-\ref{eq:HF2}) and then into the definition
\eqref{eq:MR7} of the scattering mean free paths, we obtain using the
identity
\begin{align*}
\int_{\RR^2} d \bu \, \widetilde \cR \big(k \bu,k \bu \cdot \nabla
\beta(\bka)\big) &= \int_{\RR^2} d \bu \int_{\RR^2} d \br
\int_{-\infty}^\infty d \zeta \, \cR(\br,\zeta) e^{-i k (\br,\zeta)
  \cdot (\bu,\bu \cdot \nabla \beta(\bka))} \\ &= \int_{\RR^2} d \br
\int_{-\infty}^\infty d \zeta\, \cR(\br,\zeta) (2 \pi)^2 \delta \big[
  k \big(\br + \zeta\nabla \beta(\bka) \big) \big] \\ &= \frac{(2
  \pi)^2}{k^2} \int_{-\infty}^\infty d \zeta \, \cR\left(\frac{\bka
  \zeta}{\beta(\bka)},\zeta\right),
\end{align*}
that
\begin{align} 
\cS(\bka) &= \frac{\gamma}{k^2} \, \frac{8 \beta^2(\bka)}{ {\alpha}^2
  \int_{-\infty}^\infty d \zeta\, \cR\left(\frac{\bka
    \zeta}{\beta(\bka)},\zeta\right)} + O(\gamma^2), \qquad
\cS^\perp(\bka) - \cS(\bka)= O(\gamma^2).
\label{eq:expandsmfp2}
\end{align}
This result shows that when $\gamma \ll 1$, the scattering mean free
paths are almost the same. They are proportional to $\gamma$ and
decrease as the negative power of $2$ with the frequency i.e., the
wave number $k$.  That the waves at higher frequency lose coherence
faster, is expected, as observed in Figure \ref{fig:SCMF}.

\textbf{Remark 3:} 
If $\cR(\vec\bx)$ depends only on $|\vec\bx|$,
i.e. $\cR(\vec\bx)=\cR_{\rm iso}(|\vec\bx|)$, then we find from
(\ref{eq:expresReQter}) in \ref{ap:RT} that ${\rm Re}\big[ {\bf
    Q}(\bka)\big]$ is proportional to the $2 \times 2$ identity
matrix, and
\begin{align}
\cS(\bka) =\cS^\perp(\bka) = \frac{\gamma}{k^2} \, \frac{4
  \beta(\bka)}{ {\alpha}^2\int_0^\infty d \zeta \, \cR_{\rm iso}
  (\zeta)} + O(\gamma^2).
\end{align}
This is in agreement with \eqref{eq:expandsmfp2} because 
\[
\cR \left(\frac{\bka \zeta}{\beta(\bka)},\zeta \right) = \cR_{\rm iso}
\left( \sqrt{\frac{|\bka|^2 \zeta^2}{\beta^2(\bka)} + \zeta^2} \right)
= \cR_{\rm iso}\left(\frac{|\zeta|}{\beta(\bka)}\right).
\]

\textbf{Remark 4:} These results hold also when $|\bka | \lesssim \gamma$, 
and we have 
\begin{align}
\label{eq:Shf}
\cS( \gamma \bka) &= \frac{\gamma}{k^2} \, \frac{8}{ {\alpha}^2
  \int_{-\infty}^\infty d \zeta\, \cR\big({\bf 0} ,\zeta\big)} + O(\gamma^2), \qquad
\cS^\perp(\gamma \bka) - \cS(\gamma \bka)= O(\gamma^2).
\end{align}
This shows that in the high-frequency regime,  in the dimensional variables,
the scattering mean free paths for small wave vectors are equal to $ 8/[\alpha^2 k^2
\ell \int_{-\infty}^\infty \cR ({\bf 0}, \zeta) d\zeta]$. These are the scattering mean free paths obtained
in the white-noise paraxial regime as described in \ref{app:par},
equation (\ref{eq:smfppar}). The agreement is expected because when
$\gamma\to 0$  we look at a narrow cone beam
propagating through a random medium, as in the paraxial regime.

\subsection{Evolution of the coherence matrix}
\label{sect:HF2}
In this section we consider the coherence matrix in the high-frequency regime, denoted by 
\begin{equation}
\boldsymbol{\mathscr{P}}_{\rm hf}(\bka,z) = \lim_{\gamma \to 0} \gamma^4 \boldsymbol{\mathscr{P}}(\gamma \bka, \gamma z)  ,
\end{equation}
and compare it to the coherence  matrix in the paraxial regime. We rescaled $z$ by $\gamma$ because we know from \eqref{eq:Shf} that the waves lose coherence over such propagation distances, 
and we see below that  energy exchange between the transverse electric and magnetic modes occurs at the same scale.
We also rescaled $\bka$ by $\gamma$ because the initial condition generates such wave vectors when $\gamma_{_\mJ}$
 is of the same order as $ \gamma$, as we assume here.
Taking the limit of (\ref{eq:T2}) as $\gamma \to 0$,
and using 
\begin{equation}
 \lim_{\gamma \to 0} \bGamma  (\gamma \bka,\gamma \bka')  =
\bGamma_{\rm hf} (\bka,\bka')  :=
\frac{1}{|\bka||\bka'|} 
\begin{pmatrix}
\bka \cdot \bka' & -\bka^\perp\cdot \bka' \\
\bka^\perp\cdot \bka'  & \bka \cdot \bka' 
\end{pmatrix}
\end{equation}
and 
\begin{equation}
  \lim_{\gamma \to 0} \gamma {\bf Q}  (\gamma \bka)  =
{\bf Q}_{\rm hf} (\bka) 
: = - \left[ \frac{k^4 \alpha^2}{8 } \int_{\RR^2}  \frac{d\bu}{(2\pi)^2} \widetilde \cR(k \bu,0) \right] {\bf I} ,
\end{equation}
we find that $\boldsymbol{\mathscr{P}}_{\rm hf}(\bka,z) $ satisfies
\begin{align}
 \partial_z \boldsymbol{\mathscr{P}}_{\rm hf}(\bka,z) =& 2{\bf Q}_{\rm hf} (\bka) 
 \boldsymbol{\mathscr{P}}_{\rm hf}(\bka,z) 
+\frac{k^2 {\alpha}^2}{4}\int_{\RR^2} \frac{d (k \bka')}{(2 \pi)^2}\,
 \bGamma_{\rm hf}(\bka,\bka') \boldsymbol{\mathscr{P}}_{\rm hf}(\bka',z)
 \bGamma_{\rm hf}(\bka',\bka) \widetilde \cR\big( k(\bka-\bka') , 0 \big),
 \label{eq:TRHF}
\end{align}
starting from $\boldsymbol{\mathscr{P}}_{\rm hf}(\bka,0) =\boldsymbol{\mathscr{A}}_{{\rm hf},o}(\bka)
\boldsymbol{\mathscr{A}}_{{\rm hf},o}(\bka)^\dagger$,
with
\begin{align}
\boldsymbol{\mathscr{A}}_{{\rm hf},o}(\bka) =
\begin{pmatrix}
 a_{{\rm hf},o} (\bka) \\
 a_{{\rm hf},o}^\perp(\bka) 
\end{pmatrix} ,
\qquad
a_{{\rm hf},o}(\bka) = 
-\frac{\bka}{2|\bka|}\cdot \what
  \bJ\left(\frac{k\bka}{\bar{\kappa}_{_\mJ} }\right)
 , 
  \qquad  a_{{\rm hf},o}^\perp(\bka) = -  \frac{\bka^\perp}{2|\bka|} \cdot
\what \bJ\left(\frac{k\bka}{\bar{\kappa}_{_\mJ} }\right) ,
\end{align}
and $\bar{\kappa}_{_\mJ} = \lim_{\gamma \to 0} \gamma_{_\mJ} / \gamma = O(1)$.
More explicitly, equation \eqref{eq:TRHF} written componentwise is 
\begin{align*}
\nonumber
 \partial_z  {\mathscr{P}}_{{\rm hf},11}(\bka,z) = &
 \frac{k^2 {\alpha}^2}{4}\int_{\RR^2} \frac{d (k \bka')}{(2 \pi)^2}\,\widetilde \cR\big( k(\bka-\bka') , 0 \big) \Big[  - {\mathscr{P}}_{{\rm hf},11}(\bka,z)  + 
\\ &  \frac{
  (\bka\cdot \bka')^2 {\mathscr{P}}_{{\rm hf},11}(\bka',z) 
 -(\bka\cdot\bka')(\bka^\perp\cdot\bka') ( {\mathscr{P}}_{{\rm hf},12}(\bka',z) + {\mathscr{P}}_{{\rm hf},21}(\bka',z) )
+  (\bka^\perp\cdot \bka')^2 {\mathscr{P}}_{{\rm hf},22}(\bka',z) }{|\bka|^2|\bka'|^2} \Big] ,\\
\nonumber
\partial_z  {\mathscr{P}}_{{\rm hf},12}(\bka,z) = &
 \frac{k^2 {\alpha}^2}{4}\int_{\RR^2}  \frac{d (k \bka')}{(2 \pi)^2}\,\widetilde \cR\big( k(\bka-\bka') , 0 \big) \Big[  - {\mathscr{P}}_{{\rm hf},12}(\bka,z)  +
\\ &  \frac{
(\bka\cdot\bka')(\bka^\perp\cdot\bka') ( {\mathscr{P}}_{{\rm hf},11}(\bka',z) - {\mathscr{P}}_{{\rm hf},22}(\bka',z) )
  + (\bka\cdot \bka')^2 {\mathscr{P}}_{{\rm hf},12}(\bka',z) 
  -  (\bka^\perp\cdot \bka')^2 {\mathscr{P}}_{{\rm hf},21}(\bka',z) }{|\bka|^2|\bka'|^2} \Big] , \\
\nonumber
\partial_z  {\mathscr{P}}_{{\rm hf},21}(\bka,z) = &
 \frac{k^2 {\alpha}^2}{4}\int_{\RR^2}  \frac{d (k \bka')}{(2 \pi)^2}\,\widetilde \cR\big( k(\bka-\bka') , 0 \big) \Big[  - {\mathscr{P}}_{{\rm hf},21}(\bka,z)  +
\\ & \frac{
(\bka\cdot\bka')(\bka^\perp\cdot\bka') ( {\mathscr{P}}_{{\rm hf},11}(\bka',z) - {\mathscr{P}}_{{\rm hf},22}(\bka',z) )
- (\bka^\perp\cdot \bka')^2 {\mathscr{P}}_{{\rm hf},12}(\bka',z) 
+  (\bka\cdot \bka')^2 {\mathscr{P}}_{{\rm hf},21}(\bka',z) }{|\bka|^2|\bka'|^2} \Big] ,\\
\nonumber
\partial_z  {\mathscr{P}}_{{\rm hf},22}(\bka,z) = &
 \frac{k^2 {\alpha}^2}{4}\int_{\RR^2}  \frac{d (k \bka')}{(2 \pi)^2}\,\widetilde \cR\big( k(\bka-\bka') , 0 \big) \Big[  - {\mathscr{P}}_{{\rm hf},22}(\bka,z)  +
\\ & \frac{
  (\bka^\perp\cdot \bka')^2 {\mathscr{P}}_{{\rm hf},11}(\bka',z) 
 +(\bka\cdot\bka')(\bka^\perp\cdot\bka') ( {\mathscr{P}}_{{\rm hf},12}(\bka',z) + {\mathscr{P}}_{{\rm hf},21}(\bka',z) )
+  (\bka \cdot \bka')^2 {\mathscr{P}}_{{\rm hf},22}(\bka',z) }{|\bka|^2|\bka'|^2} \Big] .
\end{align*}
This system  can be simplified by  changing of coordinates, as follows.  Let us introduce (for $\bka=(\kappa_1,\kappa_2)$)
\begin{align}
a_{1}^\ep(\bka,z) =  \frac{\kappa_1 a^\ep(\bka,z) -\kappa_2 a^{\ep,\perp}(\bka,z) }{|\bka|}, \qquad
a_{2}^\ep(\bka,z) =  \frac{\kappa_2 a^\ep(\bka,z) +\kappa_1 a^{\ep,\perp}(\bka,z) }{|\bka|},
\end{align}
and note from \eqref{eq:EField} that they are the Fourier coefficients of the electric field on the $x_1$ and $x_2$ axes.
Define also  the coherence matrix 
\begin{equation}
   \widetilde{\boldsymbol{\mathscr{P}}}_{\rm hf} (\bka,z) = \lim_{\gamma \to 0}
  \gamma^4   \widetilde{\boldsymbol{\mathscr{P}}}(\gamma \bka, \gamma z) , \qquad 
   \widetilde{\boldsymbol{\mathscr{P}}}(\bka,z) = \lim_{\ep \to 0} \EE
   \left[ \begin{pmatrix} a^\ep_1(\bka,z)
       \\ a_2^{\ep}(\bka,z) \end{pmatrix}
 \begin{pmatrix}
   a_1^\ep(\bka,z) \\ a_2^{\ep }(\bka,z) \end{pmatrix}^\dagger \right],
\end{equation}
which is related to ${\boldsymbol{\mathscr{P}}}_{\rm hf} (\bka,z) $ as
\begin{align*}
   \widetilde{ {\mathscr{P}}}_{{\rm hf},11} (\bka,z) &=
   \frac{ 
\kappa_1^2 {\mathscr{P}}_{{\rm hf},11} (\bka,z)  
+
\kappa_2^2 {\mathscr{P}}_{{\rm hf},22} (\bka,z)  
-\kappa_1 \kappa_2 (  {\mathscr{P}}_{{\rm hf},12}(\bka,z) +{\mathscr{P}}_{{\rm hf},21} (\bka,z) )
}{|\bka|^2} ,\\
   \widetilde{ {\mathscr{P}}}_{{\rm hf},12} (\bka,z) &=
   \frac{ 
\kappa_1\kappa_2 ( {\mathscr{P}}_{{\rm hf},11} (\bka,z)  
-  {\mathscr{P}}_{{\rm hf},22} (\bka,z)  )
+\kappa_1^2  {\mathscr{P}}_{{\rm hf},12}(\bka,z)  -\kappa_2^2 {\mathscr{P}}_{{\rm hf},21} (\bka,z) 
}{|\bka|^2} ,\\   
\widetilde{ {\mathscr{P}}}_{{\rm hf},21} (\bka,z) &=
   \frac{ 
\kappa_1\kappa_2 ( {\mathscr{P}}_{{\rm hf},11} (\bka,z)  
-  {\mathscr{P}}_{{\rm hf},22} (\bka,z)  )
-\kappa_2^2  {\mathscr{P}}_{{\rm hf},12}(\bka,z)  +\kappa_1^2 {\mathscr{P}}_{{\rm hf},21} (\bka,z) 
}{|\bka|^2} ,\\
   \widetilde{ {\mathscr{P}}}_{{\rm hf},22} (\bka,z) &=
   \frac{ 
\kappa_2^2 {\mathscr{P}}_{{\rm hf},11} (\bka,z)  
+
\kappa_1^2 {\mathscr{P}}_{{\rm hf},22} (\bka,z)  
+\kappa_1 \kappa_2 (  {\mathscr{P}}_{{\rm hf},12}(\bka,z) +{\mathscr{P}}_{{\rm hf},21} (\bka,z) )
}{|\bka|^2} .
\end{align*}
Elementary calculations then show that $\widetilde{\boldsymbol{\mathscr{P}}}_{\rm hf}(\bka,z)$ satisfies
\begin{align}
 \partial_z   \widetilde{\boldsymbol{\mathscr{P}}}_{\rm hf}(\bka,z) = 
 \frac{k^2 {\alpha}^2}{4}\int_{\RR^2} \frac{d (k \bka')}{(2 \pi)^2}\,\widetilde \cR\big( k(\bka-\bka') , 0 \big) \Big[  -  \widetilde{\boldsymbol{\mathscr{P}}}_{\rm hf}(\bka,z)  
  + \widetilde{\boldsymbol{\mathscr{P}}}_{\rm hf}(\bka',z)  \Big] ,
  \label{eq:evoltildePhf}
  \end{align}
starting from $ \widetilde{\boldsymbol{\mathscr{P}}}_{\rm hf}(\bka,0) 
= \widetilde{\boldsymbol{\boldsymbol{\mathscr{A}}}}_{{\rm hf},o}(\bka)
 \widetilde{\boldsymbol{\boldsymbol{\mathscr{A}}}}_{{\rm hf},o}(\bka)^\dagger$,
with
\begin{align}
 \widetilde{\boldsymbol{\boldsymbol{\mathscr{A}}}}_{{\rm hf},o}(\bka) = - \frac{1}{2} 
\begin{pmatrix}
\widehat  J_1\Big(\frac{k\bka}{\bar{\kappa}_{_\mJ} }\Big) \\
\widehat  J_2\Big(\frac{k\bka}{\bar{\kappa}_{_\mJ} }\Big) 
\end{pmatrix}.
\end{align}
This equation agrees to the one found in the white noise paraxial regime,
 as described in \ref{app:par}, equation (\ref{eq:rtepar}). It says that in this high frequency 
 scaling there is no polarization exchange. If the electric field is along the $x_1$ axis initially,
 it remains along the $x_1$ axis for distances that are comparable to the scattering mean free path.

\textbf{Remark 5:} The polarization conservation described here happens only in
the high-frequency regime, where the transport equation reduces to
\eqref{eq:evoltildePhf}.  It does not hold in general, and equation
(\ref{eq:T2}) exhibits polarization exchange.

\section{Summary}
\label{sect:sum}
We presented a detailed analysis of cumulative scattering effects on
electromagnetic waves that propagate long distances in random media
with small fluctuations of the wave speed. The results are derived
from first principles, by studying Maxwell's equations with random
electric permittivity, in a regime of separation of scales modeled by
two parameters. The first is $\ep$, the ratio of the wavelength and
the distance of propagation, and the second is $\gamma$, the ratio of
the wavelength and correlation length of the random fluctuations,
which is similar to the spatial support of the source.  By controlling
$\gamma$, we ensure that the source emits a cone wave beam which
propagates along a preferred direction, called range. Depending on
$\gamma$, the opening angle of the cone may be much larger than in
paraxial regimes. However, the angle is smaller than $180$ degrees, so
that evanescent and backscattered waves can be neglected.

The analysis of the solution of Maxwell's equation is in the long
range limit $\ep \to 0$. It involves the decomposition of the waves in
transverse electric and magnetic plane wave components, whose
amplitudes are random fields.  They satisfy a stochastic system of
differential equations driven by the random fluctuations of the
electric permittivity. These equations model the evolution in range of
the random amplitudes, starting from the initial values determined by
the source of excitation. The $\ep \to 0$ limit of the solution of the
system of stochastic differential equations is obtained with the
Markov limit theorem.

We study in detail the first two statistical moments of the limit
amplitudes.  Their expectation decays exponentially with range, on
length scales called scattering mean free paths. This is the
manifestation of the randomization of the waves, due to scattering in
the random medium. The second moments of the amplitudes describe the
decorrelation of the waves over directions, and the transport of
energy. Of particular interest is the energy density (mean Wigner
transform), which characterizes the state of polarization of the waves
and the diffusion of energy over directions. The transport equations
with polarization satisfied by the energy density are derived with the
Markov limit theorem, and the result is related to the radiative
transport theory.  In the high frequency limit $\gamma \to 0$ the
equations simplify, and are related to those satisfied by the Wigner
transform of the solution of the paraxial wave equation in random
media. We quantify the transfer of energy between the components of
the wave field, and illustrate the results with numerical simulations.

\section*{Acknowledgements}
Liliana Borcea's work was partially supported by grant \#339153 from
the Simons Foundation and by AFOSR Grant FA9550-15-1-0118. 
\appendix
\section{The Markov limit}
\label{ap:DL}
Let $\mathscr{O}$ be an open set in $\RR^d$ and ${\cal
  D}(\mathscr{O},\RR^p)$ the space of infinitely differentiable
functions with compact support.  We consider the process $\bX^\ep$ in
${\cal C}([0,L], {\cal D}')$, the solution of
\begin{equation}
\frac{d \bX^\ep}{dz} = \frac{1}{\sqrt{\ep}} {\cal
  F}\left(\frac{z}{\ep},\frac{z}{\ep}\right) \bX^\ep + {\cal
  G}\left(\frac{z}{\ep},\frac{z}{\ep}\right) \bX^\ep ,
\end{equation}
where ${\cal F}(\zeta,\zeta')$ and ${\cal G}(\zeta,\zeta')$ are random
linear operators from ${\cal D}'$ to ${\cal D}'$.  We assume that the
mappings $\zeta \to {\cal F}(\zeta,\zeta')$ and $\zeta\to {\cal
  G}(\zeta,\zeta')$ are stationary and possess strong ergodic
properties, and that ${\cal F}(\zeta,\zeta')$ has mean zero. Moreover,
the mappings $\zeta' \to {\cal F}(\zeta,\zeta')$ and $\zeta'\to {\cal
  G}(\zeta,\zeta')$ are periodic.

We are interested in particular in equation (\ref{eq:ivteps}), where
the process $\bX^\ep$ is in ${\cal D}'(\mathscr{O},\RR^4)$, defined by
\begin{align}
\bX^\ep(z) = \begin{pmatrix} {\rm Re}\big( {a}^\ep(\bka,z) )\\ {\rm
    Im}\big( {a}^\ep(\bka,z) )\\ {\rm Re} \big(
  {a}^{\ep,\perp}(\bka,z) \\ {\rm Im} \big( {a}^{\ep,\perp}(\bka,z)
\end{pmatrix}, \quad \mbox{for} ~ \bka \in \mathscr{O} = \{\bka \in \RR^2, \,
|\bka|<1\}.
\label{eq:A2}
\end{align}
The operator ${\cal F}(\zeta,\zeta')$ is
\begin{align}
\left< {\cal F}(\zeta,\zeta') \bX ,\bphi\right> &= \sum_{j=1}^4 \int_\mathscr{O} d
\bka\, [{\cal F}(\zeta,\zeta') \bX ]_j (\bka) \phi_j(\bka) 
= \int_\mathscr{O} d \bka \, \bphi(\bka) \cdot \int_\mathscr{O} d \bka' \,
\mathbb{F} (\bka,\bka',\zeta,\zeta') {\itbf X}(\bka'),
\label{eq:A3}
\end{align}
for $\boldsymbol{\phi} \in {\cal D}(\mathscr{O},\RR^4)$ with components
$\phi_j$, and ${\itbf X} \in {\cal D}'(\mathscr{O},\RR^4)$ with components
$X_j$. The kernel matrix $\mathbb{ F}(\bka,\bka',\zeta,\zeta')$ in
(\ref{eq:A3}) is given by
\begin{align}
{\bf \mathbb{F}}= \begin{pmatrix} \mathscr{F}_{11}^r &
  -\mathscr{F}_{11}^i& \mathscr{F}_{12}^r & -\mathscr{F}_{12}^i
  \\ \mathscr{F}_{11}^i & \mathscr{F}_{11}^r &\mathscr{F}_{12}^i &
  \mathscr{F}_{12}^r \\ \mathscr{F}_{21}^r & -\mathscr{F}_{21}^i
  &\mathscr{F}_{22}^r & -\mathscr{F}_{22}^i \\ \mathscr{F}_{21}^i &
  \mathscr{F}_{21}^r &\mathscr{F}_{22}^i &
  \mathscr{F}_{22}^r \end{pmatrix},
\label{eq:A4}
\end{align}
in terms of
\begin{align}
\mathscr{F}^{\rm r}_{jl}(\bka,\bka',\zeta,\zeta') &= {\rm Re}
\Big[\frac{i k^3 {\alpha}}{2 (2\pi)^2 \gamma^2} \widehat{\nu} \left(
  \frac{k(\bka-\bka')}{\gamma},\gamma \zeta \right) {F}^{aa}_{jl}
  \big(\bka,\bka',\zeta' \big) \Big]
,\label{eq:A5}\\ \mathscr{F}^{\rm i}_{jl}(\bka,\bka',\zeta,\zeta') &=
    {\rm Im} \Big[\frac{i k^3 {\alpha}}{2 (2\pi)^2 \gamma^2} \widehat{\nu}
      \left( \frac{k(\bka-\bka')}{\gamma},\gamma \zeta\right) {F}^{aa}_{jl}
      \big(\bka,\bka',\zeta' \big) \Big] , \label{eq:A6}
\end{align}
where we recall from (\ref{eq:Fzetaaa})-(\ref{eq:DefGamma}) the
expression of $F_{jl}^{aa}(\bka,\bka',\zeta')$.  The adjoint operator
${\cal F}^*(\zeta,\zeta')$ is defined by
$$ \left< {\cal F}(\zeta,\zeta') \bX ,\bphi\right> = \left< \bX ,{\cal
  F}^*(\zeta,\zeta') \bphi\right>
$$ for $\bphi \in {\cal D}(\mathscr{O},\RR^4)$ and $\bX \in {\cal
  D}'(\mathscr{O},\RR^4)$, and has kernel
$\mathbb{F}^*(\bka,\bka',\zeta,\zeta') = \mathbb{F}^T
(\bka',\bka,\zeta,\zeta')$.
Similar expressions hold for ${\cal G}$.

To obtain the Markov limit we use the results in \cite{PW1994}.  {The
  interested reader can find a self-contained introduction to such
  limit theorems in \cite[Chapter 6]{fouque07} or in \cite[Appendix
    A]{bal}.}  We get that $\bX^\ep(z)$ converges weakly in ${\cal
  C}([0,L],{\cal D}')$ to $\bX(z)$, which is the solution of a
martingale problem with generator ${\cal L}$ defined by:
\begin{align}
\nonumber
{\cal L} f(\left<\bX , \bphi \right>) =&
\int_0^\infty  d \zeta \lim_{Z\to \infty} \frac{1}{Z}
\int_0^Z d h \, \EE \big[ \left< \bX , {\cal F}^*(0,h)
  \bphi\right> \left< \bX , {\cal F}^*(\zeta,\zeta+h) \bphi\right>
  \big] \, f'' ( \left< \bX , \bphi\right> ) \\ 
  \nonumber&+
\int_0^\infty d \zeta \lim_{Z\to \infty} \frac{1}{Z}
\int_0^Z d h \, \EE \big[ \left< \bX , {\cal F}^*(0,h)
  {\cal F}^*(\zeta,\zeta+h) \bphi\right> \big] \, f' ( \left< \bX ,
\bphi\right> ) \\ &+ \lim_{Z\to \infty} \frac{1}{Z}
\int_0^Z dh \, \EE \big[ \left< \bX , {\cal G}^*(0,h)
  \bphi\right> \big] \, f' ( \left< \bX , \bphi\right> ) ,
  \label{eq:GEN0}
\end{align}
for any $\bX \in {\cal D}'(\mathscr{O},\RR^4)$, $\bphi\in {\cal
  D}(\mathscr{O},\RR^4)$, and smooth $f:\RR\to \RR$.  This
means that, for any $\bphi \in {\cal D}(\mathscr{O},\RR^4)$ and smooth
function $f :\RR\to \RR$, the real-valued process
$$
 f(\left< \bX(z) , \bphi \right>)- \int_0^z dz' \, {\cal L} f(\left<
\bX(z') , \bphi \right>)
$$ is a martingale.  More generally, if $n \in \NN$,
$\bphi_1,\ldots,\bphi_n\in {\cal D}(\mathscr{O},\RR^4)$, and $f:\RR^n \to \RR $
is a smooth function, then
$$ 
f(\left< \bX(z) , \bphi_1\right>,\ldots,\left< \bX(z) ,
\bphi_n\right>)- \int_0^z dz' \, {\cal L}^n f(\left< \bX(z') ,
\bphi_1\right>,\ldots,\left< \bX(z') , \bphi_n\right>) 
$$
is a martingale,
where 
\begin{align}
&  {\cal L}^n f(\left< \bX , \bphi_1\right>,\ldots,\left< \bX ,
  \bphi_n\right>)\nonumber \\ & \quad  =\left\{\sum_{j,l=1}^n
  \int_0^\infty d \zeta \lim_{Z\to \infty}
  \frac{1}{Z}\int_0^Z d h \, \EE \big[ \left< \bX ,
    {\cal F}^*(0,h) \bphi_j\right> \left< \bX , {\cal
      F}^*(\zeta,\zeta+h) \bphi_l\right> \big]\, \partial^2_{jl}
  \right.  \nonumber \\ &\qquad + \sum_{j=1}^n
  \int_0^\infty  d\zeta \lim_{Z\to \infty}
  \frac{1}{Z}\int_0^Z dh \, \EE \big[ \left< \bX ,
    {\cal F}^*(0,h) {\cal F}^*(\zeta,\zeta+h) \bphi_j\right> \big] \,
  \partial_{j}  \nonumber \\ &\left. \qquad + \sum_{j=1}^n \lim_{Z\to \infty}
  \frac{1}{Z}\int_0^Z   dh \, \EE \big[ \left< X ,
    {\cal G}^*(0,h) \phi_j\right> \big]  \partial_j \right\}f(\left< \bX ,
  \bphi_1\right>,\ldots,\left< \bX , \bphi_n\right>\big) . \label{eq:GEN}
\end{align}

To calculate the first moment of the limit process $\bX(z)$, let $n =
1$ and $f(y)=y$ in (\ref{eq:GEN}). We find that
\begin{eqnarray*}
\frac{ d \, \EE \big[ \left< {\bX}(z) , \bphi \right>\big]}{dz} = \EE
\Big[ \left< \bX(z) , \overline{\cal F}^*\bphi \right>\Big] + \EE
\Big[ \left< \bX(z) , \overline{\cal G}^*\bphi \right>\Big] ,
\end{eqnarray*}
where 
\begin{eqnarray*}
\overline{\cal F}^* &=& \int_0^\infty d \zeta
\lim_{Z\to \infty}\frac{1}{Z}\int_0^Z d h \, \EE
\big[ {\cal F}^*(0,h) {\cal F}^*(\zeta,\zeta+h) \big]
,\\ \overline{\cal G}^* &=& \lim_{Z\to \infty} \frac{1}{Z}
\int_0^Z  d h \, \EE \big[ {\cal G}^*(0,h) \big] .
\end{eqnarray*}
This shows that
$$
 \overline{\bX}(z) = \EE \big[ \bX(z) \big] ,
$$ 
satisfies a closed system of ordinary differential equations
\begin{eqnarray*}
\frac{ d \left< \overline{\bX}(z) , \bphi \right> }{dz} = \left<
\overline{\bX}(z) , \overline{\cal F}^*\bphi \right> + \left<
\overline{\bX}(z) , \overline{\cal G}^*\bphi \right> ,
\end{eqnarray*}
or, equivalently in ${\cal D}'$,
\begin{equation}
\frac{ d \overline{\bX}(z) }{dz} = \overline{\cal F}
\,\overline{\bX}(z) + \overline{\cal G} \, \overline{\bX}(z) ,
\label{eq:MeanField}
 \end{equation}
where $\overline{\cal F}$, resp. $\overline{\cal G}$, is the adjoint
of $\overline{\cal F}^*$, resp. $\overline{\cal G}^*$.

Recalling from (\ref{eq:A3}-\ref{eq:A6}) the expression of the
kernel $\mathbb{F}^T(\kappa',\bka,\zeta,\zeta')$ of ${\cal
  F}^*(\zeta,\zeta')$, we obtain
\begin{align*}
\overline{\cal F}^*_{jl}(\bka,\bka') =& \sum_{q=1}^4 \int_\mathscr{O} d \bka''
\int_0^\infty  d \zeta \lim_{Z\to \infty} \frac{1}{Z}
\int_0^Z   dh \, \EE \big[ {\mathbb
    F}_{lq}(\bka',\bka'',\zeta,\zeta+h) {\mathbb
    F}_{qj}(\bka'',\bka,0,h) \big] ,
\end{align*}
for $j,l = 1, \ldots, 4$.  For instance,
\begin{align*}
\overline{\cal F}^*_{11}(\bka,\bka') =& \int_\mathscr{O} d \bka'' \int_0^\infty
  d \zeta \lim_{Z\to \infty} \frac{1}{Z}
\int_0^Z  dh \, \EE \big[ {\mathscr F}^{\rm
    r}_{11}(\bka',\bka'',\zeta,\zeta+h) {\mathscr F}^{\rm
    r}_{11}(\bka'',\bka,0,h) \big]\\ &-
 \int_\mathscr{O} d \bka'' \int_0^\infty
 d \zeta \lim_{Z\to \infty} \frac{1}{Z}
\int_0^Z  dh \,\EE \big[
  {\mathscr F}^{\rm i}_{11}(\bka',\bka'',\zeta,\zeta+h) {\mathscr
    F}^{\rm i}_{11}(\bka'',\bka,0,h) \big]  \\ &+
\int_\mathscr{O} d \bka'' \int_0^\infty
 d \zeta \lim_{Z\to \infty} \frac{1}{Z}
\int_0^Z   dh \, \EE \big[ {\mathscr F}^{\rm
    r}_{12}(\bka',\bka'',\zeta,\zeta+h) {\mathscr F}^{\rm
    r}_{21}(\bka'',\bka,0,h) \big] \\ &- \int_\mathscr{O} d \bka''
\int_0^\infty
 d \zeta \lim_{Z\to \infty} \frac{1}{Z}
\int_0^Z dh \, \EE \big[ {\mathscr F}^{\rm
    i}_{12}(\bka',\bka'',\zeta,\zeta+h) {\mathscr F}^{\rm
    i}_{21}(\bka'',\bka,0,h) \big],
\end{align*}
and using (\ref{eq:A5}-\ref{eq:A6}), we get 
\begin{align*}
\overline{\cal F}^*_{11}(\bka,\bka') = {\rm Re} \Big\{ \left(
\frac{ik^3 {\alpha}}{2 (2\pi)^2 \gamma^2}\right)^2\int_\mathscr{O} d
\bka'' \int_0^\infty d \zeta\lim_{Z\to \infty} \frac{1}{Z} \int_0^Z d
h \, \EE\left[ \widehat{\nu}\left(\frac{k(\bka'-\bka'')}{\gamma},
  \gamma\zeta\right)\widehat{\nu}\left(\frac{k(\bka''-\bka)}{\gamma},0\right)
  \right] \times \\ \big[ {\bf F}^{aa} (\bka',\bka'',\zeta+h) {\bf
    F}^{aa}(\bka'',\bka,h) \big]_{11} \Big\} .
\end{align*}
Moreover, using the identity
\begin{eqnarray*}
\EE \left[
  \widehat{\nu}\left(\frac{k(\bka'-\bka'')}{\gamma},\gamma\zeta\right)
\widehat{\nu}\left( \frac{k(\bka''-\bka)}{\gamma},0\right)\right]
= \left( \frac{2\pi \gamma}{k}\right)^2 \delta (\bka-\bka') \widehat{\cR}\left( \frac{k(\bka-\bka'')}{\gamma},\gamma\zeta \right) ,
 \end{eqnarray*}
derived from the definition of the autocorrelation with
straightforward algebraic manipulations, and obtaining from
(\ref{eq:Fzetaaa}) that 
\begin{eqnarray*}
{\bf F}^{aa} (\bka,\bka'',\zeta+h) {\bf F}^{aa}(\bka'',\bka,h) =
\bGamma^{aa} (\bka,\bka'' ) \bGamma^{aa}(\bka'',\bka ) e^{i k
  (\beta(\bka'')-\beta(\bka)) \zeta},
\end{eqnarray*}
we get 
\begin{align*}
\overline{\cal F}^*_{11}(\bka,\bka') = -\frac{k^4 {\alpha}^2}{16 \pi^2 \gamma^2}
         {\rm Re}\Big\{ \int_\mathscr{O} d \bka'' \int_0^\infty d \zeta \,
           \widehat{\cR}\left(
           \frac{k(\bka-\bka'')}{\gamma},\gamma\zeta \right) e^{i k
             (\beta(\bka'')-\beta(\bka)) \zeta} \\ \times \big[
             \bGamma^{aa}(\bka,\bka'')\bGamma^{aa}(\bka'',\bka)
             \big]_{11} \,\delta(\bka-\bka') \Big\}.
\end{align*}

The expressions of the other components of $\overline{\cal
  F}^*_{jl}(\bka,\bka')$ are of the same type. Similarly, we
calculate $\overline{\cal G}^\star$, and substituting into
(\ref{eq:MeanField}) we obtain the explicit expression of the
differential equations satisfied by the mean wave amplitudes. This is
equation (\ref{eq:Amean}), written in complex form.

The calculation of the second moments is similar, by letting $n = 1$
and $f(y) = y^2$ in (\ref{eq:GEN}), and carrying the lengthy
calculations.

\section{Connection to the radiative transport theory}
\label{ap:RT}
To connect our transport equations (\ref{eq:T7}) (see also
(\ref{eq:T2})) to the radiative transport theory in
\cite{chandra,ryzhik96,bal00}, we adhere to the notation in
\cite{ryzhik96}.  First note that in \cite{ryzhik96} the random medium
is assumed to be statistically isotropic, i.e. $\cR(\vec\bx)$ depends
only on $|\vec\bx|$.  

For any $\vk \in \RR^3$, consider the orthonormal basis $\{\vec{\itbf
  z}^{(0)}(\vk),\vec{\itbf z}^{(1)}(\vk),\vec{\itbf z}^{(2)}(\vk)\}$
of $\mathbb{R}^3$, defined by
\begin{equation}
\vec{\itbf z}^{(0)} (\vk)= \frac{\vk}{|\vk|} = \frac{1}{|\vk|} \left(
\bmk, \mk_z \right), \quad \vec{\itbf z}^{(1)}(\vk) = \frac{1}{|\vk|}
\left( \mk_z \frac{\bmk}{|\bmk|} , - |\bmk| \right), \quad \vec{\itbf
  z}^{(2)}({\vk}) = \left( \frac{\bmk^\perp}{|\bmk|}, 0 \right),
\label{eq:rt1}
\end{equation}
and satisfying
\begin{equation}
\vec{\itbf z}^{(0)} \times \vec{\itbf z}^{(1)} = \vec{\itbf z}^{(2)},
\qquad \vec{\itbf z}^{(0)} \times \vec{\itbf z}^{(2)} = - \vec{\itbf
  z}^{(1)}.
\label{eq:RT3}
\end{equation}
In \cite{ryzhik96} the vector $\vec{\itbf z}^{(0)}$ is denoted by
$\widehat{{\itbf k}}$.  We call it $\vec{\itbf z}^{(0)}$ to avoid
confusion with the wavenumber $k$.  Note that when $\vk = k \vka$,
with $|\vka|=1$, then $\vec{\itbf z}^{(0)} = \vka$, and the vectors
$\vec{\itbf z}^{(1)}$ and $\vec{\itbf z}^{(2)}$ are the same as $\vu$
and $\vu^\perp$ defined in \eqref{eq:ORT_SR}.

Following the notation in \cite{ryzhik96}, we define 
\begin{equation}
f_j(\vk,\vx) = \frac{2 \pi}{\sqrt{2} \ep} \Big[ \vec{\itbf E}^\ep(\vx)
  \cdot \vec{\itbf z}^{(j)}(\vk) + \zeta_o \vec{\itbf H}^\ep(\vx)
  \cdot \left( \vec{\itbf z}^{(0)}(\vk) \times \vec{\itbf
    z}^{(j)}(\vk)\right)\Big], \quad j = 1, 2,
\label{eq:RT4}
\end{equation}
where we use a different constant of proportionality than in
\cite{ryzhik96}, to simplify the relation \eqref{eq:RT6}.  The Wigner
transform $\bW(\vk,\vx)$ is the $2 \times 2$ matrix with components
\begin{equation}
W_{lq}(\vk,\vx) = \int \frac{d \vy}{(2 \pi)^3} \, f_l\left(\vx -
\frac{\ep \vy}{2}, \vk \right) \overline{f_q \left( \vx + \frac{\ep
    \vy}{2},\vk \right)} e^{ i \vk \cdot \vy}, \quad l,q= 1,2,
\label{eq:RT5}
\end{equation}
where the bar denotes complex conjugate. Substituting the wave
decompositions (\ref{eq:EField}-\ref{eq:MField}) into \eqref{eq:RT4},
and keeping only the forward propagating modes, we obtain after some
algebraic manipulations that
\begin{equation}
\label{eq:RT6}
 \bW(\vk,\vx) = \frac{\delta \big[ \mk_z - k \beta(\bmk/k)
     \big]}{\beta(\bmk/k)} \bcW(\bmk/k,\bx,z),
\end{equation}
with $\bcW$ the Wigner transform in \eqref{eq:T6}, satisfying
\eqref{eq:T7}.

The transport equation in \cite{ryzhik96} is 
\begin{align}
 \vna \om(\vk) \cdot \vna_{\vx} \bW(\vk,\vx)
= \int d \vk' \sigma( \vk,\vk') [\bW(\vk',\vx)] -
\Sigma(\vk)\bW(\vk,\vx),
\label{eq:RTT}
\end{align}
for the dispersion relation
\[
\om(\vk) = c_o |\vk|.
\]
The integral kernel in the right hand side of \eqref{eq:RTT} is the
differential scattering cross section, defined by \cite{ryzhik96}
\begin{align}
\sigma( \vk,\vk') [ \bW(\vk',\vx) ] = \frac{\pi c_o^2 k^2
  {\alpha}^2}{2(2 \pi)^3\gamma^3} \widetilde
\cR \left(\frac{\vk-\vk'}{\gamma}\right) \delta[\om(\vk)-\om(\vk')]
%\nonumber \\ \times 
{\bf  T}(\vk,\vk') \bW(\vk',\vx) {\bf  T}(\vk',\vk),
\label{eq:RT9}
\end{align}
with 
the $2 \times 2$ matrix 
$$
{\bf T}(\vk,\vk') = \big( \vec{\itbf z}^{(l)}(\vk) \cdot
\vec{\itbf z}^{(q)}(\vk') \big)_{l,q = 1,2}.
$$
The scalar $\Sigma(\vk)$ is the total scattering cross section,
which is shown in \cite{ryzhik96} to satisfy
\begin{align}
\Sigma(\vk) {\bf I} = \int d\vk' \sigma(\vk,\vk') [{\bf I}] . 
\label{eq:RT9b}
\end{align}

Since \eqref{eq:RT6} gives that $\bW(\vk,\vx)$ is supported at vectors $\vk$ of
the form $\vk = k \vka$, with $\vka = (\bka,\beta(\bka))$, we see that
the operator on the left hand side of (\ref{eq:RTT}) is given by
\begin{equation}
\vna \om(\vk) \cdot \vna_{\vx} = c_o \frac{\vk}{|\vk|} \cdot
\vna_{\vx} = c_o \big[ \bka \cdot \nabla_{\bx} + \beta(\bka)
  \partial_z\Big] = c_o \beta(\bka) \Big[ \partial_z - \nabla
  \beta(\bka) \cdot \nabla_{\bx} \Big].
\label{eq:RT8}
\end{equation}
Again, using \eqref{eq:RT6}, we see that the integral kernel in the
right hand side of (\ref{eq:RTT}) is supported at vectors $\vk' = k
\vka'$, with $\vka' = (\bka',\beta(\bka'))$, so the Dirac distribution
in \eqref{eq:RT9} is
\begin{align}
\delta \Big[ \om(\vk)-\om(k \vka')\Big] &= \delta \Big[
  c_o\sqrt{|\bmk|^2 + |\mk_z|^2}-c_o k \Big] = \frac{\delta \big[
    \mk_z - k \beta(\bmk/k)\big]}{c_o \beta(\bmk/k)}.
\label{eq:RT10}
\end{align}
Thus, the integral in \eqref{eq:RTT} is supported at 
$\vk = k \vka$, with $\vka = (\bka,\beta(\bka))$. 
For such vectors, the matrix ${\bf T}$ equals 
\begin{equation}
{\bf T}(k \vka,k \vka') = \sqrt{\beta(\bka)
  \beta(\bka')}\bGamma(\bka,\bka').
\label{eq:RT7}
\end{equation}
From \eqref{eq:RT6} and (\ref{eq:RT8}), we
obtain that
\begin{align}
 \vna \om(\vk) \cdot \vna_{\vx} \bW(\vk,\vx) &= c_o \delta \big[\mk_z
   - k \beta(\bmk/k)\big]\big[ \partial_z - \nabla \beta(\bmk/k) \cdot
   \nabla_{\bx} \big] \bcW(\bmk/k,\bx,z), \label{eq:RT11}
\end{align}
from \eqref{eq:RT6}, (\ref{eq:RT9}), and (\ref{eq:RT10}-\ref{eq:RT7}), we
obtain that
\begin{align}
\int d \vk' \sigma( \vk,\vk') [\bW(\vk',\vx)] = \frac{c_o
    k^2{\alpha}^2}{4 \gamma^3} \delta \big[\mk_z - k
    \beta(\bmk/k)\big] \int_{|\bka'|\leq 1}
  \frac{d (k \bka')}{(2 \pi)^2} \, \widetilde
  \cR \left(\frac{\bmk-k\bka'}{\gamma},\frac{k(\beta(\bmk/k)-\beta(\bka'))}{\gamma}
  \right) \times \nonumber \\   \bGamma(\bmk/k,\bka') \bcW(\bka', \bx,z) \bGamma(\bka',\bmk
  / k).\label{eq:RT12}
\end{align}
We also find from (\ref{eq:RT9b}) that
\begin{align}
\Sigma(\vk) {\bf I} = \frac{ c_o^2 k^2{\alpha}^2}{4(2 \pi)^2\gamma^3}
\int d \vk' \delta \big[ \omega(\vk')-\omega(\vk)\big] \widetilde
\cR \left(\frac{\vk-\vk'}{\gamma }\right) {\bf T}(\vk,\vk')
   {\bf T}(\vk',\vk) ,
\label{eq:RT12a}
\end{align}
or equivalently, for $\vk=k(\bka,\beta(\bka))$,
\begin{align}
\Sigma(\vk) {\bf I} = \frac{ c_o k^2 \beta(\bka)
  {\alpha}^2}{4\gamma^3} \int_{|\bka'| < 1} \frac{d (k \bka')}{(2
  \pi)^2} \widetilde
\cR \left(\frac{k(\bka-\bka')}{\gamma},\frac{k(\beta(\bka)-\beta(\bka'))}{\gamma}
\right) \bGamma(\bka,\bka') \bGamma(\bka' ,\bka) ,
\label{eq:RT12aa}
\end{align}
so that
\begin{align}
\nonumber \Sigma(\vk) \bW(\vk,\vx) = \frac{ c_o
  k^2{\alpha}^2}{4\gamma^3} \delta \big[ \mk_z - k \beta(\bmk/k) \big]
\int_{|\bka'| < 1} \frac{d (k \bka')}{(2 \pi)^2} \widetilde \cR
\left(\frac{\bmk-k\bka'}{\gamma},\frac{k
  \big(\beta(\bmk/k)-\beta(\bka')\big)}{\gamma} \right) \times
\\ \bGamma(\bmk/k,\bka') \bGamma(\bka',\bmk / k )\bcW(\bmk/k,\bx,z) .
\label{eq:RT12b}
\end{align}
Finally, using the transport equation \eqref{eq:T7} satisfied by
$\bcW(\bka,\bx,z)$, and the relation \eqref{eq:RT6}, we obtain that
\begin{align}
 \vna \om(\vk) \cdot \vna_{\vx} \bW(\vk,\vx) =& \int d \vk'\,
 \sigma(\vk,\vk') \bW(\vk',\vx)
  + c_o \beta(\bmk/k) \Big[
   \mathbf{Q}(\bmk/k) \bW(\vk,\vx) + \bW(\vk,\vx)
   \mathbf{Q}(\bmk/k)^\dagger\Big].
\label{eq:RT13}
\end{align}
This is similar to the transport equation \eqref{eq:RTT}, except that
we do not have the scalar valued total scattering cross section
$\Sigma(\vk)$ multiplying $\bW(\vk,\vx)$, but a linear operator $- c_o
\beta(\bmk/k)\Big[ \mathbf{Q}(\bmk/k) \bW(\vk,\vx) + \bW(\vk,\vx)
  \mathbf{Q}(\bmk/k)^\dagger\Big]$ acting on $\bW(\vk,\vx)$. The
results would be exactly the same if $\mathbf{Q}(\bka)$ were a
multiple of the identity of the form \[- 2c_o \beta(\bka)
\mathbf{Q}(\bka) = \Sigma(\vk) {\bf I}\] for
$\vk=k(\bka,\beta(\bka))$.  This turns out to be the case in the
high-frequency limit $\gamma \to 0$, as explained in Section
\ref{sect:DEPOL}.  It is also the case for the real part of
$\mathbf{Q}(\bka)$ in general (provided that $\cR$ is isotropic, as
assumed throughout this appendix), as we now show.

For $\bka\in \RR^2$ such that $|\bka|<1$, we have by (\ref{eq:C5pa})
  \begin{align}
 {\rm Re}\big( \mathbf{Q}(\bka)\big) = -
 \frac{k^2{\alpha}^2}{8\gamma^3} \int_{|\bka'| < 1} \frac{d (k
   \bka')}{(2 \pi)^2} \widetilde \cR
   \left(\frac{k(\bka-\bka')}{\gamma},\frac{k
   \big(\beta(\bka)-\beta(\bka')\big)}{\gamma} \right)
 \bGamma(\bka,\bka') \bGamma(\bka' ,\bka) .
    \label{eq:expresReQ}
  \end{align}
By comparing with (\ref{eq:RT12aa}) and denoting $\vk=k(\bka,\beta(\bka))$, this shows that
\begin{align}
- 2 c_o \beta(\bka) {\rm Re}\big( \mathbf{Q}(\bka)\big) 
=  \Sigma(\vk)  {\bf I} ,
  \end{align}
which is a multiple of the identity matrix.  This also follows by
inspection from (\ref{eq:expresReQ}), rewritten as
\begin{align}
 {\rm Re}\big( \mathbf{Q}(\bka)\big) = -\frac{ k^2{\alpha}^2}{8(2
   \pi)^2\gamma^3 \beta(\bka)} \int_{|\vk'|=k} d S(\vk') \widetilde
 \cR  \left(\frac{\vk-\vk'}{\gamma}\right) {\bf
   T}(\vk,\vk') {\bf T}(\vk',\vk) .
  \label{eq:expresReQbis}
  \end{align}
Indeed, if we parameterize
$$
\bka = \begin{pmatrix} 
\sin \theta \cos \phi \\
\sin \theta \sin \phi 
\end{pmatrix},\quad \quad
\vk = k \begin{pmatrix}
\sin \theta \cos \phi \\
\sin \theta \sin \phi  \\
\cos \theta
\end{pmatrix},\quad \quad 
\vk' = k \begin{pmatrix}
\sin \theta' \cos \phi' \\
\sin \theta' \sin \phi'  \\
\cos \theta'
\end{pmatrix},
$$
then we have
$$
{\bf  T}(\vk,\vk') {\bf  T}(\vk',\vk)
= \begin{pmatrix}
1 - \kappa_2''^2 & -\kappa_1'' \kappa_2'' \\
-\kappa_1''\kappa_2'' & 1- \kappa_1''^2
\end{pmatrix} ,\quad
\quad
|\vk-\vk'|^2=2 k^2  (1-\kappa_3'') ,
$$
where
$$
\vka''
=
\begin{pmatrix}
\kappa_1''\\
\kappa_2''\\
\kappa_3''
\end{pmatrix}
=
\begin{pmatrix}
 \sin \theta' \sin(\phi'-\phi)  \\
\cos \theta \sin \theta' \cos(\phi'-\phi) +\sin \theta \cos \theta'\\
\sin \theta \sin \theta' \cos(\phi'-\phi) +\cos \theta \cos \theta'
\end{pmatrix} .
$$
If we introduce the rotation matrix
$$
{\bf U}_\bka = \begin{pmatrix} -\sin \phi & \cos \phi & 0\\
\cos \theta \cos \phi & \cos \theta \sin \phi & -\sin \theta\\
\sin \theta \cos \phi & \sin \theta \sin \phi & \cos \theta
\end{pmatrix} ,
$$
then 
$$
{\bf U}_\bka \vk' = k \vka'' .
$$ By carrying out the change of variable $\vka'' = {\bf U}_\bka \vk'
/ k$ in (\ref{eq:expresReQbis}), we find that since $\widetilde
\cR (\vka)=\check\cR_{\rm iso}(|\vka|)$,
\begin{align}
 {\rm Re}\big( \mathbf{Q}(\bka)\big) = -\frac{ k^4{\alpha}^2}{8(2
   \pi)^2\gamma^3 \beta(\bka)} \int_{|\vka''|=1} d S(\vka'')
 \check\cR_{\rm iso}\left(\frac{ \sqrt{2
     (1-\kappa_3'')}}{\gamma}\right) \begin{pmatrix} 1 - \kappa_2''^2
   & -\kappa_1'' \kappa_2'' \\ -\kappa_1''\kappa_2'' & 1- \kappa_1''^2
\end{pmatrix}  .
  \label{eq:expresReQter}
  \end{align}
By changing $(\kappa''_1,\kappa''_2)$ into $(-\kappa''_2,\kappa''_1)$
and carrying the integration over the unit sphere, we obtain that
${\rm Re}\big( \mathbf{Q}(\bka)\big)$ is proportional to the identity
matrix.

\section{Connection to the white-noise paraxial theory}
\label{app:par}%
The white-noise paraxial theory is valid when the four length scales
satisfy
\begin{equation}
\la   \ll \ell \sim X \ll \bar{L},
\label{eq:hypscal1par}
\end{equation}
with $\la \bar{L} \sim X^2$, which corresponds to a Fresnel number of
order one.  The separation of scales can then be characterized by a
unique dimensionless parameter $\ep \ll 1$ such that
\begin{equation}
\ep = \frac{\la}{\bar{L}}, \qquad \frac{\la}{\ell} \sim \frac{\la}{X}
\sim \ep^{\frac{1}{2}}.
\end{equation}
We again let $\bar{L}$ be the reference length scale and introduce the scaled
length variables
\begin{equation}
\bx' = \bx/(\ep^{\frac{1}{2}}\bar{L}), \quad z' = z/\bar{L}, \quad L' = L/\bar{L}, \quad
\ell' = \ell /(\ep^{\frac{1}{2}}\bar{L})  ,\quad
X' = X/(\ep^{\frac{1}{2}}\bar{L})   .
\end{equation} 
The scaled wavenumber is $k' = k \bar{L} \ep = 2\pi$. We finally
assume that the standard deviation ${\alpha}$ of the fluctuations of
the random medium in (\ref{eq:F12}) is small, of order $\ep^{3/2}$,
and we let ${\alpha}' = {\alpha} / \ep^{3/2}$ be the normalized
standard deviation of the fluctuations.  If we drop the primes and
consider the asymptotic regime $\ep \to 0$, then the electric field
$\vbE^\ep$ has the form~\cite{garnier2009paraxial}
\begin{align}
  \vbE^\ep (\vx) \approx  \begin{pmatrix}
   \bE _{\rm par}(\bx,z) \\ 0\end{pmatrix} 
   e^{i k \frac{z}{\ep}}, \quad \quad 
  \bE _{\rm par}(\bx,z) =   \int_{\RR^2} 
T_{\rm par}(\bx,\bx',z) \bE_{o}(\bx')  d\bx' ,
\end{align}
where the scalar transmission kernel satisfies the It\^o-Schr\"odinger equation
\begin{align}
d T_{\rm par}(\bx,\bx' , z) &= \frac{i}{2k} \Delta_{\bx'} T_{\rm
  par}(\bx,\bx',z)dz +\frac{ik}{2} T_{\rm par}(\bx,\bx',z) \circ
dB(\bx',z) ,
\end{align}
starting from $T_{\rm par}(\bx,\bx',z=0)=\delta(\bx-\bx')$.
Here $B(\bx,z)$ is a Brownian field with mean zero and covariance function
$$ \EE \big[ B(\bx,z) B(\bx',z')\big ] =\min(z,z') \, \alpha^2 \ell
\int_{-\infty}^\infty \cR \left(
\frac{\bx-\bx'}{\ell},\zeta\right)d\zeta ,
$$
and $\circ$ stands for the Stratonovich integral.
Using It\^o's formula,  the coherent field
 ${\itbf A}_{\rm par}(\bx,z) =  \EE [\bE_{\rm par} (\bx,z)]$
satisfies~\cite{garnier2009paraxial}
\begin{align*}
\partial_z {\itbf A}_{\rm par}(\bx,z) = \frac{i}{2k} \Delta_\bx {\itbf
  A}_{\rm par}(\bx,z)- \frac{1}{{\mathcal S}_{\rm par} } {\itbf
  A}_{\rm par}(\bx,z) ,
\end{align*}
which shows that it decays exponentially with $z$ on the scale
\begin{align}
\label{eq:smfppar}
{\mathcal S}_{\rm par} = \frac{8}{k^2 \ell \alpha^2
  \int_{-\infty}^\infty \cR({\bf 0},\zeta) d \zeta} .
\end{align}
The decay length ${\mathcal S}_{\rm par}$ corresponds to the
scattering mean free path $\cS(\bka)$ defined by
(\ref{eq:expandsmfp2}) in the high-frequency regime $\gamma \to 0$
with $|\bka|=O(\gamma)$.

The Wigner transform is 
\begin{align*}
\nonumber {\bf W}_{\rm par} (\bka,\bx,z) = \int_{\RR^2} d \by \, e^{  i
  k \bka \cdot \by } \EE \Big[ \bE_{\rm
    par}\Big(\bx-\frac{\by}{2},z\Big) \bE_{\rm
    par}\Big(\bx+\frac{\by}{2} ,z\Big)^\dagger \Big] \\ = \int_{\RR^2}
\frac{ d (k\bq)}{(2\pi)^2} \, e^{ i k\bq\cdot \bx } \EE \Big[
  \widehat{\bE}_{\rm par}\left( k\Big(\bka+\frac{\bq}{2}\Big),z\right)
  \widehat{\bE}_{\rm par}\left(k\Big( \bka-\frac{\bq}{2}\Big),z \right)^\dagger
  \Big] .
\end{align*}
Using It\^o's formula it is shown in \cite{garnier2009paraxial} to
satisfy the transport equation
\begin{align}
  \partial_z {\bf W}_{\rm par}+ \bka \cdot
\nabla_\bx {\bf W}_{\rm par} = \frac{k^2
  \ell^{3} \alpha^2}{4} \int_{\RR^2} \frac{{\rm d} (k
  \bka')}{(2\pi)^2}\, \widetilde \cR \big( k \ell (\bka-\bka') ,0
\big) \big[ {\bf W}_{\rm par}(\bka') - {\bf W}_{\rm par}(\bka) \big] ,
\label{eq:rtepar}
\end{align}
which corresponds to (\ref{eq:evoltildePhf}) in the high-frequency regime $\gamma \to 0$.

\bibliographystyle{elsarticle-num}

\bibliography{EM_POLARIZ}

\end{document}

%% file: mathmacros.tex
\DeclareMathAlphabet{\itbf}{OML}{cmm}{b}{it}
 \DeclareMathAlphabet\mathbfcal{OMS}{cmsy}{b}{n}
 
\def\EE{\mathbb{E}}

\def\RR{\mathbb{R}}
\def\NN{\mathbb{N}}
\def\eps{\varepsilon}

\def\bq{{{\itbf q}}}
\def\bu{{{\itbf u}}}

\def\by{{{\itbf y}}}
\def\bx{{{\itbf x}}}

\def\br{{{\itbf r}}}
\def\bX{{{\itbf X}}}

\newcommand{\bka}{{\boldsymbol{\kappa}}}

\newcommand{\bpsi}{{\boldsymbol{\psi}}}

\newcommand{\bGamma}{{\boldsymbol{\Gamma}}}

\newcommand{\la}{\lambda}
\newcommand{\vbE}{\vec{\itbf E}}
\newcommand{\bE}{{\itbf E}}
\newcommand{\bH}{{\itbf H}}
\newcommand{\vbH}{\vec{\itbf H}}
\newcommand{\vbJ}{\vec{{\itbf J}}}
\newcommand{\bJ}{{\itbf J}}
\newcommand{\vx}{\vec{\itbf x}}
\newcommand{\vy}{\vec{\itbf y}}
\newcommand{\vr}{\vec{\itbf r}}
\newcommand{\vna}{\vec{\nabla}}
\newcommand{\om}{\omega}
\newcommand{\mE}{E}
\newcommand{\mH}{H}
\newcommand{\mJ}{J}
\newcommand{\bU}{{\itbf U}}
\newcommand{\cR}{{\mathscr R}}
\newcommand{\ep}{\epsilon}
\newcommand{\what}{\widehat}
\newcommand{\bphi}{\boldsymbol{\phi}}
\newcommand{\cS}{{\mathcal{S}}}

\newcommand{\bcW}{\boldsymbol{\mathscr{W}}}
\newcommand{\bW}{{\bf W}}
\newcommand{\vu}{\vec {\boldsymbol{{u}}}}

\newcommand{\vka}{\vec{\bka}}
\newcommand{\vk}{\vec{\boldsymbol{\mathscr{K}}}}

\newcommand{\mk}{\mathscr{K}}
\newcommand{\bmk}{\boldsymbol{\mk}}